\documentclass[12pt]{amsart} 
\usepackage{amsfonts,graphics,amsmath,amsthm,amsfonts,amscd,amssymb,amsmath,latexsym,multicol,euscript}
\usepackage{epsfig,url}
\usepackage{flafter}
\usepackage{hyperref}
\usepackage[UglyObsolete]{diagrams}

\diagramstyle[noPostScript]

\makeatletter

\def\jobis#1{FF\fi
  \def\preedicate{#1}%
  \edef\preedicate{\expandafter\strip@prefix\meaning\preedicate}%
  \edef\job{\jobname}%
  \ifx\job\preedicate
}

\makeatother

\if\jobis{proposal}%
 \def\try{subsection}%
\else
  \def\try{section}%
\fi

\theoremstyle{plain}
\newtheorem{theorem}{Theorem}[\try]
\newtheorem{corollary}[theorem]{Corollary}
\newtheorem{lemma}[theorem]{Lemma}

\newtheorem{proposition}[theorem]{Proposition}
\newtheorem{definition-lemma}[theorem]{Definition-Lemma}

\newtheorem{definition}[theorem]{Definition}
\newtheorem{remark}[theorem]{Remark}

\newtheorem{conjecture}[theorem]{Conjecture}

\newtheorem{theorema}{Theorem}

\def\scr#1{\mathbf{\EuScript{#1}}}

\def\ideal#1.{I_{#1}}
\def\ring#1.{\mathcal {O}_{#1}}
\def\fring#1.{\hat{\mathcal {O}}_{#1}}
\def\proj#1.{\mathbb {P}(#1)}
\def\pr #1.{\mathbb {P}^{#1}}
\def\dpr #1.{\hat{\mathbb {P}}^{#1}}
\def\af #1.{\mathbb{A}^{#1}}
\def\Hz #1.{\mathbb{F}_{#1}}
\def\Hbz #1.{\overline{\mathbb {F}}_{#1}}
\def\fb#1.{\underset {#1} {\times}}
\def\rest#1.{\underset {\ \ring #1.} \to \otimes}
\def\au#1.{\operatorname {Aut}\,(#1)}
\def\deg#1.{\operatorname {deg } (#1)}
\def\pic#1.{\operatorname {Pic}\,(#1)}
\def\pico#1.{\operatorname{Pic}^0(#1)}
\def\picg#1.{\operatorname {Pic}^G(#1)}
\def\ner#1.{NS (#1)}
\def\rdown#1.{\llcorner#1\lrcorner}
\def\rfdown#1.{\lfloor{#1}\rfloor}
\def\rup#1.{\ulcorner{#1}\urcorner}
\def\rcup#1.{\lceil{#1}\rceil}
\def\cone#1.{\operatorname {NE}(#1)}
\def\mone#1.{\operatorname {NM}(#1)}
\def\ccone#1.{\overline{\operatorname {NE}}(#1)}
\def\cmone#1.{\overline{\operatorname {NM}}(#1)}
\def\none#1.{\operatorname {NF}(#1)}
\def\cnone#1.{\overline{\operatorname {NF}}(#1)}
\def\coef#1.{\frac{(#1-1)}{#1}}
\def\vit#1.{D_{\langle #1 \rangle}}
\def\mm#1.{\overline {M}_{0,#1}}
\def\H1#1.{H^1(#1,{\ring #1.})}
\def\ac#1.{\overline {\mathbb F}_{#1}}

\def\adj#1.{\frac {#1-1}{#1}}
\def\spn#1.{\overline{#1}}
\def\pek#1.#2.{\Cal P^{#1}(#2)}
\def\plk#1.#2.{\Cal P^{\leq #1}(#2)}
\def\ev#1.{\operatorname{ev_{#1}}}
\def\ilist#1.{{#1}_1,{#1}_2,\ldots}
\def\bminv#1.{(\nu_1,s_1;\nu_2,s_2;\dots ;\nu_{#1},s_{#1};\nu_{r+1})}
\def\zinv#1.{(\nu_1,s_1;\nu_2,s_2;\dots ;\nu_{#1},s_{#1};0)}
\def\iinv#1.{(\nu_1,s_1;\nu_2,s_2;\dots ;\nu_{#1},s_{#1};\infty)}
\def\scr#1.{\mathbf{\EuScript{#1}}}
\def\mg#1.{\overline {M}_{#1}}
\def\inter#1.{\underset #1{\cdot}}
\def\WDiv{\operatorname{WDiv}}

\def\llist#1.#2.{{#1}_1,{#1}_2,\dots,{#1}_{#2}}
\def\ulist#1.#2.{{#1}^1,{#1}^2,\dots,{#1}^{#2}}
\def\lomitlist#1.#2.{{#1}_1,{#1}_2,\dots,\hat {{#1}_i}, \dots, {#1}_{#2}}
\def\lomitlistz#1.#2.{{#1}_0,{#1}_1,\dots,\hat {{#1}_i}, \dots, {#1}_{#2}}
\def\loc#1.#2.{\Cal O_{#1,#2}}
\def\fderiv#1.#2.{\frac {\partial #1}{\partial #2}}
\def\deriv#1.#2.{\frac {d #1}{d #2}}
\def\map#1.#2.{#1 \longrightarrow #2}
\def\rmap#1.#2.{#1 \dasharrow #2}
\def\emb#1.#2.{#1 \hookrightarrow #2}
\def\non#1.#2.{\text {Spec }#1[\epsilon]/(\epsilon)^{#2}}
\def\Hi#1.#2.{\text {Hilb}^{#1}(#2)}
\def\sym#1.#2.{\operatorname {Sym}^{#1}(#2)}
\def\Hb#1.#2.{\text {Hilb}_{#1}(#2)}
\def\Hm#1.#2.{\Hom_{#1}(#2)}
\def\prd#1.#2.{{#1}_1\cdot {#1}_2\cdots {#1}_{#2}}
\def\Bl #1.#2.{\operatorname {Bl}_{#1}#2}
\def\pl #1.#2.{#1^{\otimes #2}}
\def\mgn#1.#2.{\overline {M}_{#1,#2}}
\def\ialist#1.#2.{{#1}_1 #2 {#1}_2, #2\dots}
\def\pair#1.#2.{\langle #1, #2\rangle}
\def\gproj#1.#2.{\mathbb{P}_{#1}(#2)}
\def\gpr #1.#2.{\mathbb{P}^{#1}_{#2}}
\def\gaf #1.#2.{\mathbb{A}^{#1}_{#2}}
\def\vandermonde#1.#2.{\left|
\begin{matrix}
1 & 1 & 1 & \dots & 1\\
{#1}_1 & {#1}_2 & {#1}_3 & \dots & {#1}_{#2}\\
{#1}_1^2 & {#1}_2^2 & {#1}_3^2 & \dots & {#1}_{#2}^2\\
\vdots & \vdots & \vdots & \ddots & \vdots\\
{#1}_1^{#2-1} & {#1}_2^{#2-1} & {#1}_2^{#2-1} & \dots & {#1}_{#2}^{#2-1}\\
\end{matrix}
\right|
}
\def\vandermondet#1.#2.{\left|
\begin{matrix}
1 & {#1}_1   & {#1}_1^2 & \dots & {#1}_1^{#2-1}\\
1 & {#1}_2   & {#1}_2^2 & \dots & {#1}_2^{#2-1}\\
1 & {#1}_3   & {#1}_3^2 & \dots & {#1}_3^{#2-1}\\
\vdots & \vdots & \vdots & \ddots & \vdots\\
1 & {#1}_{#2}& {#1}_{#2}^2 & \dots & {#1}_{#2}^{#2-1}\\
\end{matrix}
\right|
}
\def\gr#1.#2.{\mathbb{G}(#1,#2)}

\def\alist#1.#2.#3.{{#1}_1 #2 {#1}_2 #2\dots #2 {#1}_{#3}}
\def\zlist#1.#2.#3.{#1_0 #2 #1_1 #2\dots #2 #1_{#3}}
\def\lomitlist30#1.#2.#3.{{#1}_0,{#1}_1 #2 \dots #2\hat {{#1}_i} #2\dots #2 {#1}_{#3}}
\def\lmap#1.#2.#3.{#1 \overset{#2}{\longrightarrow} #3}
\def\mes#1.#2.#3.{#1 \longrightarrow #2 \longrightarrow #3}
\def\ses#1.#2.#3.{0\longrightarrow #1 \longrightarrow #2 \longrightarrow #3 \longrightarrow 0}
\def\les#1.#2.#3.{0\longrightarrow #1 \longrightarrow #2 \longrightarrow #3}
\def\res#1.#2.#3.{#1 \longrightarrow #2 \longrightarrow #3\longrightarrow 0}
\def\Hi#1.#2.#3.{\text {Hilb}^{#1}_{#2}(#3)}
\def\ten#1.#2.#3.{#1\underset {#2}{\otimes} #3}
\def\lomitlist30#1.#2.#3.{{#1}_0 #2 {#1}_1 #2 \dots #2 \hat {{#1}_i} #2 \dots #2 {#1}_{#3}}
\def\mderiv#1.#2.#3.{\frac {d^{#3} #1}{d #2^{#3}}}
\def\ggr#1.#2.#3.{\mathbb{G}_{#1}(#2,#3)}

\def\Hom{\operatorname{Hom}}

\def\Supp{\operatorname{Supp}}

\def\dim{\operatorname{dim}}

\def\deg{\operatorname{deg}}

\def\mult{\operatorname{mult}}

\def\fix{\operatorname{Fix}}

\def\rest{\operatorname{res}}

\def\adj{\operatorname{adj}}

\def\e{\Cal E}

\def\e1{E_1}
\def\e2{E_2}

\def\mapdown#1{\big\downarrow\rlap{$\vcenter
{\hbox{$\scriptstyle#1$}}$}}

\def\mapse#1{
{\vcenter{\hbox{$\mathop{\smash{\raise1pt\hbox{$\diagdown$}\!\lower7pt
\hbox{$\searrow$}}\vphantom{p}}\limits_{#1}\vphantom{\mapdown{}}$}}}}

\def\VR#1.{height#1pt&\omit&&\omit&&\omit&&\omit&&\omit&\cr}

\def\VRT#1.{height#1pt&\omit&&\omit&\cr}

\begin{document}
\title{Existence of minimal models for varieties of log general type}  
\author{Caucher Birkar}
\address{DPMMS\\
Centre for Mathematical Sciences\\
University of Cambridge\\
Wilberforce Road\\
Cambridge CB3 0WB, UK}
\email{c.birkar@dpmms.cam.ac.uk}
\author{Paolo Cascini} 
\address{Department of Mathematics\\ 
University of California at Santa Barbara\\ 
Santa Barbara, CA 93106, USA}
\email{cascini@math.ucsb.edu}
\author{Christopher D. Hacon} 
\date{\today}
\address{Department of Mathematics \\  
University of Utah\\  
155 South 1400 East\\
JWB 233\\
Salt Lake City, UT 84112, USA}
\email{hacon@math.utah.edu}
\author{James M\textsuperscript{c}Kernan} 
\address{Department of Mathematics\\ 
University of California at Santa Barbara\\ 
Santa Barbara, CA 93106, USA} 
\email{mckernan@math.ucsb.edu}
\address{Department of Mathematics\\ 
MIT\\ 
77 Massachusetts Avenue\\
Cambridge, MA 02139, USA} 
\email{mckernan@math.mit.edu}

\thanks{The first author was partially supported by EPSRC grant GR/S92854/02, the second
  author was partially supported by NSF research grant no: 0801258, the third author was
  partially supported by NSF research grant no: 0456363 and an AMS Centennial fellowship
  and the fourth author was partially supported by NSA grant no: H98230-06-1-0059 and NSF
  grant no: 0701101.  The last author would like to thank Sogang University and Professor
  Yongnam Lee for their generous hospitality, where some of the work for this paper was
  completed.  We would like to thank Dan Abramovich, Valery Alexeev, Florin Ambro, Tommaso
  deFernex, Stephane Dreul, Se\'an Keel, Kalle Karu, J\'anos Koll\'ar, S\'andor Kov\'acs,
  Michael McQuillan, Shigefumi Mori, Martin Olsson, Genia Tevelev, Burt Totaro, Angelo
  Vistoli and Chengyang Xu for answering many of our questions and pointing out some
  errors in an earlier version of this paper.}

\begin{abstract} We prove that the canonical ring of a smooth projective variety is 
finitely generated.  
\end{abstract}

\maketitle

\tableofcontents

\section{Introduction}
\label{s_introduction}

The purpose of this paper is to prove the following result in birational algebraic
geometry:
\begin{theorem} Let $(X,\Delta )$ be a projective kawamata log terminal pair.  

If $\Delta$ is big and $K_X+\Delta$ is pseudo-effective then $K_X+\Delta$ has a log
terminal model.
\end{theorem} 

In particular, it follows that if $K_X+\Delta$ is big then it has a log canonical model
and the canonical ring is finitely generated.  It also follows that if $X$ is a smooth
projective variety, then the ring
$$
R(X,K_X)=\bigoplus_{m\in \mathbb{N}}H^0(X,\ring X.(mK_X)),
$$ 
is finitely generated.

The birational classification of complex projective surfaces was understood by the Italian
Algebraic Geometers in the early 20th century: If $X$ is a smooth complex projective
surface of non-negative Kodaira dimension, that is $\kappa (X,K_X)\geq 0$, then there is a
unique smooth surface $Y$ birational to $X$ such that the canonical class $K_Y$ is nef
(that is $K_Y\cdot C\geq 0$ for any curve $C\subset Y$).  $Y$ is obtained from $X$ simply
by contracting all $-1$-curves, that is all smooth rational curves $E$ with $K_X\cdot
E=-1$.  If, on the other hand, $\kappa (X,K_X)=-\infty$, then $X$ is birational to either
$\pr 2.$ or a ruled surface over a curve of genus $g>0$.

The minimal model program aims to generalise the classification of complex projective
surfaces to higher dimensional varieties.  The main goal of this program is to show that
given any $n$-dimensional smooth complex projective variety $X$, we have:
\begin{itemize}
\item If $\kappa (X,K_X)\geq 0$, then there exists a \textbf{minimal model}, that
is a variety $Y$ birational to $X$ such that $K_Y$ is nef.
\item If $\kappa (X,K_X)=-\infty$, then there is a variety $Y$ birational to $X$ which
admits a \textbf{Fano fibration}, that is a morphism $\map Y.Z.$ whose general fibres $F$
have ample anticanonical class $-K_F$.
\end{itemize}

It is possible to exhibit 3-folds which have no smooth minimal model, see for example
(16.17) of \cite{Ueno82}, and so one must allow varieties $X$ with singularities.
However, these singularities cannot be arbitrary.  At the very minimum, we must still be
able to compute $K_X\cdot C$ for any curve $C\subset X$.  So, we insist that $K_X$ is
$\mathbb{Q}$-Cartier (or sometimes we require the stronger property that $X$ is
$\mathbb{Q}$-factorial).  We also require that $X$ and the minimal model $Y$ have the same
pluricanonical forms.  This condition is essentially equivalent to requiring that the
induced birational map $\phi\colon\rmap X.Y.$ is $K_X$-negative.

There are two natural ways to construct the minimal model (it turns out that if one can
construct a minimal model for pseudo-effective $K_X$, then we can construct Mori fibre
spaces whenever $K_X$ is not pseudo-effective).  Since one of the main ideas of this paper
is to blend the techniques of both methods, we describe both methods.

 The first method is to use the ideas behind finite generation.  If the canonical ring
$$
R(X,K_X)=\bigoplus_{m\in \mathbb{N}}H^0(X,\ring X.(mK_X)),
$$
is finitely generated and $K_X$ is big, then the canonical model $Y$ is nothing more
than the Proj of $R(X,K_X)$.  It is then automatic that the induced rational map
$\phi\colon\rmap X.Y.$ is $K_X$-negative.

The other natural way to ensure that $\phi$ is $K_X$-negative is to factor $\phi$ into a
sequence of elementary steps all of which are $K_X$-negative.  We now explain one way to
achieve this factorisation.

If $K_X$ is not nef, then, by the cone Theorem, there is a rational curve $C\subset X$
such that $K_X\cdot C<0$ and a morphism $f\colon\map X.Z.$ which is surjective, with
connected fibres, on to a normal projective variety and which contracts an irreducible
curve $D$ if and only if $[D]\in \mathbb{R}^+[C] \subset N_1(X)$.  Note that $\rho
(X/Z)=1$ and $-K_X$ is $f$-ample.  We have the following possibilities:
\begin{itemize}
\item If $\dim Z<\dim X$, this is the required Fano fibration.
\item If $\dim Z=\dim X$ and $f$ contracts a divisor, then we say that $f$ is a divisorial
contraction and we replace $X$ by $Z$.
\item If $\dim Z=\dim X$ and $f$ does not contract a divisor, then we say that $f$ is a
small contraction.  In this case $K_Z$ is not $\mathbb{Q}$-Cartier, so that we can not
replace $X$ by $Z$.  Instead, we would like to replace $f\colon\map X.Z.$ by its flip
$f^+\colon\map X^+.Z.$ where $X^+$ is isomorphic to $X$ in codimension $1$ and $K_{X^+}$
is $f^+$-ample.  In other words, we wish to replace some $K_X$-negative curves by
$K_{X^+}$-positive curves.
\end{itemize}
The idea is to simply repeat the above procedure until we obtain either a minimal model or
a Fano fibration.  For this procedure to succeed, we must show that flips always exist and
that they eventually terminate.  Since the Picard number $\rho (X)$ drops by one after
each divisorial contraction and is unchanged after each flip, there can be at most finitely many divisorial contractions.
So we must show that there is no infinite sequence of flips.

This program was successfully completed for $3$-folds in the 1980's by the work of
Kawamata, Koll\'ar, Mori, Reid, Shokurov and others.  In particular, the existence of
$3$-fold flips was proved by Mori in \cite{Mori88}.

Naturally, one would hope to extend these results to dimension $4$ and higher by induction
on the dimension.

Recently, Shokurov has shown the existence of flips in dimension $4$ \cite{Shokurov03} and
Hacon and M\textsuperscript{c}Kernan \cite{HM05c} have shown that assuming the minimal
model program in dimension $n-1$ (or even better simply finiteness of minimal models in
dimension $n-1$), then flips exist in dimension $n$.  Thus we get an inductive approach to
finite generation.

Unfortunately the problem of showing termination of an arbitrary sequence of flips seems
to be a very difficult problem and in dimension $\geq 4$ only some partial answers are
available.  Kawamata, Matsuda and Matsuki proved \cite{KMM87} the termination of terminal
$4$-fold flips, Matsuki has shown \cite{Matsuki91} the termination of terminal $4$-fold
flops and Fujino has shown \cite{Fujino04} the termination of canonical $4$-fold (log)
flips.  Alexeev, Hacon and Kawamata \cite{AHK06} have shown the termination of kawamata
log terminal $4$-fold flips when the Kodaira dimension of $-(K_X+\Delta )$ is non-negative
and the existence of minimal models of kawamata log terminal $4$-folds when either
$\Delta$ or $K_X+\Delta$ is big by showing the termination of a certain sequence of flips
(those that appear in the MMP with scaling).  However, it is known that termination of
flips follows from two natural conjectures on the behaviour of the log discrepancies of
$n$-dimensional pairs (namely the ascending chain condition for minimal log discrepancies
and semicontinuity of log discrepancies, cf.  \cite{Shokurov04}).  Moreover, if $\kappa
(X,K_X+\Delta)\geq 0$, Birkar has shown \cite{Birkar05} that it suffices to establish acc
for log canonical thresholds and the MMP in dimension one less.

We now turn to the main result of the paper:
\begin{theorem}\label{t_main} Let $(X,\Delta)$ be a kawamata log terminal pair, where 
$K_X+\Delta$ is $\mathbb{R}$-Cartier.  Let $\pi\colon\map X.U.$ be a projective morphism
of quasi-projective varieties.

If either $\Delta$ is $\pi$-big and $K_X+\Delta$ is $\pi$-pseudo-effective or $K_X+\Delta$ is 
$\pi$-big, then
\begin{enumerate} 
\item $K_X+\Delta$ has a log terminal model over $U$,
\item if $K_X+\Delta$ is $\pi$-big then $K_X+\Delta$ has a log canonical model over $U$, and
\item if $K_X+\Delta$ is $\mathbb{Q}$-Cartier, then the $\ring U.$-algebra
$$
\mathfrak{R}(\pi ,\Delta)=\bigoplus_{m\in\mathbb{N}} \pi_*\ring X.(\rdown m(K_X+\Delta).),
$$
is finitely generated.  
\end{enumerate} 
\end{theorem}

We now present some consequences of \eqref{t_main}, most of which are known to follow from the
MMP.  Even though we do not prove termination of flips, we are able to derive many of the
consequences of the existence of the MMP.  In many cases we do not state the strongest
results possible; anyone interested in further applications is directed to the references.
We group these consequences under different headings.

\makeatletter
\renewcommand{\thetheorem}{\thesubsection.\arabic{theorem}}
\@addtoreset{theorem}{subsection}
\makeatother

\subsection{Minimal models}

An immediate consequence of \eqref{t_main} is:
\begin{corollary}\label{c_general} Let $X$ be a smooth projective variety of general 
type.

Then
\begin{enumerate} 
\item $X$ has a minimal model,
\item $X$ has a canonical model, 
\item the ring
$$
R(X,K_X)=\bigoplus_{m\in\mathbb{N}} H^0(X,\ring X.(mK_X)),
$$
is finitely generated, and
\item $X$ has a model with a K\"ahler-Einstein metric.
\end{enumerate} 
\end{corollary}

Note that (4) follows from (2) and Theorem D of \cite{EGZ06}.  Note that Siu has announced
a proof of finite generation for varieties of general type, using analytic methods, see
\cite{Siu06}.
 
\begin{corollary}\label{c_generated} Let $(X,\Delta)$ be a projective kawamata log terminal pair,
where $K_X+\Delta$ is $\mathbb{Q}$-Cartier.
 
 Then the ring
$$
R(X,K_X+\Delta)=\bigoplus_{m\in\mathbb{N}} H^0(X,\ring X.(\rdown m(K_X+\Delta).)),
$$
is finitely generated.  
\end{corollary}

Let us emphasize that in \eqref{c_generated} we make no assumption about $K_X+\Delta$ or
$\Delta$ being big.  Indeed Fujino and Mori, \cite{FM02}, proved that \eqref{c_generated}
follows from the case when $K_X+\Delta$ is big.

We will now turn our attention to the geography of minimal models.  It is well known that
log terminal models are not unique.  The first natural question about log terminal models
is to understand how any two are related.  In fact there is a very simple connection:

\begin{corollary}\label{c_flops} Let $\pi\colon\map X.U.$ be a projective morphism of 
normal quasi-projective varieties.  Suppose that $K_X+\Delta$ is kawamata log terminal and
$\Delta$ is big over $U$.  Let $\phi_i\colon\rmap X.Y_i.$, $i=1$ and $2$ be two log
terminal models of $(X,\Delta)$ over $U$.  Let $\Gamma_i=\phi_{i*}\Delta$.

Then the birational map $\rmap Y_1.Y_2.$ is the composition of a sequence of
$(K_{Y_1}+\Gamma_1)$-flops over $U$.
\end{corollary}

Note that \eqref{c_flops} has been generalised recently to the case when $\Delta$ is not
assumed big, \cite{Kawamata07}.  The next natural problem is to understand how many
different models there are.  Even if log terminal models are not unique, in many important
contexts, there are only finitely many.  In fact Shokurov realised that much more is true.
He realised that the dependence on $\Delta$ is well-behaved.  To explain this, we need
some definitions:
\begin{definition}\label{d_cones} Let $\pi\colon\map X.U.$ be a projective morphism of 
normal quasi-projective varieties, and let $V$ be a finite dimensional affine subspace of
the real vector space $\WDiv_{\mathbb{R}}(X)$ of Weil divisors on $X$.  Fix an
$\mathbb{R}$-divisor $A\geq 0$ and define
\begin{align*} 
V_A &=\{\,\Delta \,|\, \Delta=A+B, B\in V \,\}, \\
\mathcal{L}_A(V)&=\{\,\Delta=A+B\in V_A \,|\, \text{$K_X+\Delta$ is log canonical and $B\geq 0$} \,\}, \\
\mathcal{E}_{A,\pi}(V) &=\{\,\Delta\in \mathcal{L}_A(V) \,|\, \text{$K_X+\Delta$ is pseudo-effective over $U$} \,\}, \\
\mathcal{N}_{A,\pi}(V)&=\{\,\Delta\in \mathcal{L}_A(V) \,|\, \text{$K_X+\Delta$ is nef over $U$} \,\}.
\end{align*}
Given a birational contraction $\phi\colon\rmap X.Y.$ over $U$, define
$$
\mathcal{W}_{\phi,A,\pi}(V)=\{\, \Delta\in \mathcal E_{A,\pi}(V)  \,|\, \text{$\phi$ is a
  weak log canonical model for $(X,\Delta)$ over $U$} \,\},
$$
and given a rational map $\psi\colon\rmap X.Z.$ over $U$, define
$$
\mathcal{A}_{\psi,A,\pi}(V)=\{\, \Delta\in \mathcal E_{A,\pi}(V) \,|\, \text{$\psi$ is the
  ample model for $(X,\Delta)$ over $U$} \,\},
$$
(cf. \eqref{d_log} and \eqref{d_ample} for the definitions of weak log canonical model and
ample model for $(X,\Delta)$ over $U$).
\end{definition}

We will adopt the convention that $\mathcal{L}(V)=\mathcal{L}_0(V)$.  If the support of
$A$ has no components in common with any element of $V$ then the condition that $B\geq 0$
is vacuous.  In many applications, $A$ will be an ample $\mathbb{Q}$-divisor over $U$.  In
this case, we often assume that $A$ is \textit{general} in the sense that we fix a
positive integer such that $kA$ is very ample over $U$, and we assume that $A=\frac 1k A'$, where
$A'\sim _{U} kA$ is very general.  With this choice of $A$, we have
$$
\mathcal{N}_{A,\pi}(V)\subset \mathcal{E}_{A,\pi}(V)\subset\mathcal{L}_A(V)=\mathcal{L}(V)+A \subset V_A=V+A,
$$
and the condition that the support of $A$ has no common components with any element of $V$
is then automatic.  The following result was first proved by Shokurov \cite{Shokurov96}
assuming the existence and termination of flips:

\begin{corollary}\label{c_finite} Let $\pi\colon\map X.U.$ be a projective morphism of 
normal quasi-projective varieties.  Let $V$ be a finite dimensional affine subspace of
$\WDiv _{\mathbb{R}}(X)$ which is defined over the rationals.  Suppose there is a divisor
$\Delta_0\in V$ such that $K_X+\Delta_0$ is kawamata log terminal.  Let $A$ be a general ample
$\mathbb{Q}$-divisor over $U$, which has no components in common with any element of $V$.
\begin{enumerate} 
\item\label{i_model} There are finitely many birational contractions $\phi_i\colon\rmap
X.Y_i.$ over $U$, $1\leq i\leq p$ such that
$$
\mathcal{E}_{A,\pi}(V)=\bigcup_{i=1}^p\mathcal{W}_i,
$$
where each $\mathcal{W}_i=\mathcal{W}_{\phi_i,A,\pi}(V)$ is a rational polytope.  Moreover,
if $\phi\colon\rmap X.Y.$ is a log terminal model of $(X,\Delta)$ over $U$, for some $\Delta\in
\mathcal{E}_{A,\pi}(V)$, then $\phi=\phi_i$, for some $1\leq i\leq p$.
\item \label{i_cmodel} There are finitely many rational maps $\psi_j\colon\rmap X.Z_j.$ over $U$,
$1\leq j\leq q$ which partition $\mathcal{E}_{A,\pi}(V)$ into the subsets
$\mathcal{A}_j=\mathcal{A}_{\psi_j,A,\pi}(V)$.
\item \label{i_subset} For every $1\leq i\leq p$ there is a $1\leq j\leq q$ and a morphism
$f_{i,j}\colon\map Y_i.Z_j.$ such that $\mathcal{W}_i\subset\bar{\mathcal{A}}_j$.
\end{enumerate} 

In particular $\mathcal{E}_{A,\pi}(V)$ and each $\bar{\mathcal{A}}_j$ are rational
polytopes.
\end{corollary}

\begin{definition}\label{d_kodaira} Let $(X,\Delta)$ be a kawamata log terminal pair and
let $D$ be a big divisor.  Suppose that $K_X+\Delta$ is not pseudo-effective.  The
\textbf{effective log threshold} is
$$
\sigma(X,\Delta,D)=\sup \{\, t\in\mathbb{R} \,|\, \text{$D+t(K_X+\Delta)$ is pseudo-effective} \,\}.
$$
The \textbf{Kodaira energy} is the reciprocal of the effective log threshold.  
\end{definition}

Following ideas of Batyrev, one can easily show that:
\begin{corollary}\label{c_kodaira}  Let $(X,\Delta)$ be a projective kawamata log terminal 
pair and let $D$ be an ample divisor.  Suppose that $K_X+\Delta$ is not pseudo-effective.

If both $K_X+\Delta$ and $D$ are $\mathbb{Q}$-Cartier then the effective log threshold and
the Kodaira energy are rational.
\end{corollary}

\begin{definition}\label{d_cox} Let $\pi\colon\map X.U.$ be a projective morphism of 
normal quasi-projective varieties.  Let $D^{\bullet}=(\llist D.k.)$ be a sequence of
$\mathbb{Q}$-divisors on $X$.  The sheaf of $\ring U.$-modules
$$
\mathfrak{R}(\pi,D^{\bullet})=\bigoplus_{m\in\mathbb{N}^k} \pi_*\ring X.(\rdown \sum m_iD_i.),
$$
is called the \textbf{Cox ring} associated to $D^{\bullet}$.  
\end{definition}

Using \eqref{c_finite} one can show that adjoint Cox rings are finitely generated:
\begin{corollary}\label{c_cox} Let $\pi\colon\map X.U.$ be a projective morphism of 
normal quasi-projective varieties.  Fix $A\geq 0$ an ample $\mathbb{Q}$-divisor over $U$.
Let $\Delta_i=A+B_i$, for some $\mathbb{Q}$-divisors $\llist B.k.\geq 0$.  Assume that
$D_i=K_X+\Delta_i$ is divisorially log terminal and $\mathbb{Q}$-Cartier.  Then the Cox
ring
$$
\mathfrak{R}(\pi,D^{\bullet})=\bigoplus_{m\in\mathbb{N}^k} \pi_*\ring X.(\rdown \sum m_iD_i.),
$$
is a finitely generated $\ring U.$-module.  
\end{corollary}

\subsection{Moduli spaces}

At first sight \eqref{c_finite} might seem a hard result to digest.  For this reason, we
would like to give a concrete, but non-trivial example.  The moduli spaces $\mgn g.n.$ of
$n$-pointed stable curves of genus $g$ are probably the most intensively studied moduli
spaces.  In particular the problem of trying to understand the related log canonical
models via the theory of moduli has attracted a lot of attention (e.g. see \cite{GKM00},
\cite{Manin04} and \cite{HH06}).
\begin{corollary}\label{c_curves} Let $X=\mgn g.n.$ the moduli space of 
stable curves of genus $g$ with $n$ marked points and let $\Delta_i$, $1\leq i\leq k$
denote the boundary divisors.

Let $\Delta=\sum _i a_i\Delta_i$ be a boundary.  Then $K_X+\Delta$ is log canonical and if
$K_X+\Delta$ is big then there is a log canonical model $\rmap X.Y.$.  Moreover if we fix
a positive rational number $\delta$ and require that the coefficient $a_i$ of $\Delta_i$
is at least $\delta$ for each $i$ then the set of all log canonical models obtained this
way is finite.
\end{corollary}

\subsection{Fano varieties}

The next set of applications is to Fano varieties.  The key observation is that given any 
divisor $D$, a small multiple of $D$ is linearly equivalent to a divisor of the 
form $K_X+\Delta$, where $\Delta$ is big and $K_X+\Delta$ is kawamata log terminal.  
 
 Using this observation we can show:
\begin{corollary}\label{c_dream} Let $\pi\colon\map X.U.$ be a projective morphism 
of normal varieties, where $U$ is affine.  Suppose that $X$ is $\mathbb{Q}$-factorial,
$K_X+\Delta$ is divisorially log terminal and $-(K_X+\Delta)$ is ample over $U$.  

 Then $X$ is a Mori dream space.  
\end{corollary}

There are many reasons why Mori dream spaces (see \cite{HK00} for the definition) are
interesting.  As the name might suggest, they behave very well with respect to the minimal
model program.  Given any divisor $D$, one can run the $D$-MMP, and this ends with either
a nef model, or a fibration, for which $-D$ is relatively ample, and in fact any sequence
of $D$-flips terminates.

\eqref{c_dream} was conjectured in \cite{HK00} where it is also shown that Mori dream
spaces are GIT quotients of affine varieties by a torus.  Moreover the decomposition given
in \eqref{c_finite} is induced by all the possible ways of taking GIT quotients, as one
varies the linearisation.

Finally, it was shown in \cite{HK00} that if one has a Mori dream space, then the Cox Ring
is finitely generated.

We next prove a result that naturally complements \eqref{t_main}.  We show that if
$K_X+\Delta $ is not pseudo-effective, then we can run the MMP with scaling to get a Mori
fibre space:
\begin{corollary}\label{c_fibre} Let $(X,\Delta)$ be a $\mathbb{Q}$-factorial 
kawamata log terminal pair.  Let $\pi\colon\map X.U.$ be a projective morphism of normal
quasi-projective varieties.  Suppose that $K_X+\Delta$ is not $\pi$-pseudo-effective.

Then we may run $f\colon\rmap X.Y.$ a $(K_X+\Delta)$-MMP over $U$ and end with a Mori
fibre space $g\colon\map Y.W.$ over $U$.
\end{corollary}

Note that we do not claim in \eqref{c_fibre} that however we run the $(K_X+\Delta)$-MMP
over $U$, we always end with a Mori fibre space, that is we do not claim that every
sequence of flips terminates.  

Finally we are able to prove a conjecture of Batyrev on the closed cone of nef 
curves for a Fano pair.  

\begin{definition}\label{d_nef-cone} Let $X$ be a projective variety.  A curve $\Sigma$ is 
called \textbf{nef} if $B\cdot \Sigma\geq 0$ for all Cartier divisors $B\geq 0$.  
$\none X.$ denotes the cone of nef curves sitting inside $H_2(X,\mathbb{R})$ and $\cnone
X.$ denotes its closure.

Now suppose that $(X,\Delta)$ is a log pair.  A $(K_X+\Delta)$-\textbf{co-extremal ray} is
an extremal ray $F$ of the closed cone of nef curves $\cnone X.$ on which $K_X+\Delta$ is
negative.
\end{definition}

\begin{corollary}\label{c_batyrev} Let $(X,\Delta)$ be a projective $\mathbb{Q}$-factorial 
kawamata log terminal pair such that $-(K_X+\Delta)$ is ample.  

Then $\cnone X.$ is a rational polyhedron.  If $F=F_i$ is a co-extremal ray then there
exists a $\mathbb{R}$-divisor $\Theta$ such that the pair $(X,\Theta)$ is kawamata log
terminal and the $(K_X+\Theta)$-MMP $\pi\colon\rmap X.Y.$ ends with a Mori fibre space
$f\colon\map Y.Z.$ such that $F$ is spanned by the pullback to $X$ of the class of any
curve $\Sigma$ which is contracted by $f$.
\end{corollary}

\subsection{Birational geometry}

Another immediate consequence of \eqref{t_main} is the existence of flips:
\begin{corollary}\label{c_existence} Let $(X,\Delta)$ be a kawamata log terminal 
pair and let $\pi\colon\map X.Z.$ be a small $(K_X+\Delta)$-extremal contraction.  

 Then the flip of $\pi$ exists.  
\end{corollary}

As already noted, we are unable to prove the termination of flips in general.  However,
using \eqref{c_finite}, we can show that any sequence 
of flips for the MMP with scaling terminates:
\begin{corollary}\label{c_scaling} Let $\pi\colon\map X.U.$ be a projective morphism of 
normal quasi-projective varieties.  Let $(X,\Delta)$ be a $\mathbb{Q}$-factorial kawamata
log terminal pair, where $K_X+\Delta$ is $\mathbb{R}$-Cartier and $\Delta$ is $\pi$-big.
Let $C\geq 0$ be an $\mathbb{R}$-divisor.

If $K_X+\Delta+C$ is kawamata log terminal and $\pi$-nef, then we may run the
$(K_X+\Delta)$-MMP over $U$ with scaling of $C$.
\end{corollary}

Another application of \eqref{t_main} is the existence of log terminal models which 
extract certain divisors: 
\begin{corollary}\label{c_extract} Let $(X,\Delta)$ be a log canonical pair and let 
$f\colon\map W.X.$ be a log resolution.  Suppose that there is a divisor $\Delta_0$ such
that $K_X+\Delta_0$ is kawamata log terminal.  Let $\mathfrak{E}$ be any set of valuations
of exceptional divisors which satisfies the following two properties:
\begin{enumerate} 
\item $\mathfrak{E}$ contains only valuations of log discrepancy at most one, and
\item the centre of every valuation of log discrepancy one in $\mathfrak{E}$ does
not contain any non kawamata log terminal centres.
\end{enumerate} 

Then we may find a birational morphism $\pi\colon\map Y.X.$, such that $Y$ is
$\mathbb{Q}$-factorial and the exceptional divisors of $\pi$ correspond to the elements of
$\mathfrak{E}$.
\end{corollary}

For example, if we assume that $(X,\Delta)$ is kawamata log terminal and we let
$\mathfrak{E}$ be the set of all exceptional divisors with log discrepancy at most one,
then the birational morphism $\pi\colon\map Y.X.$ defined in \eqref{c_extract} above is a
\textbf{terminal model} of $(X,\Delta)$.  In particular there is an $\mathbb{R}$-divisor
$\Gamma\geq 0$ on $W$ such that $K_Y+\Gamma=\pi ^*(K_X+\Delta)$ and the pair $(Y,\Gamma)$
is terminal.

If instead we assume that $(X,\Delta)$ is kawamata log terminal but $\mathfrak{E}$ is
empty, then the birational morphism $\pi\colon\map Y.X.$ defined in \eqref{c_extract}
above is a log terminal model.  In particular $\pi$ is small, $Y$ is
$\mathbb{Q}$-factorial and there is an $\mathbb{R}$-divisor $\Gamma\geq 0$ on $Y$ such
that $K_Y+\Gamma=\pi^*(K_X+\Delta)$.

We are able to prove that every log pair admits a birational model with
$\mathbb{Q}$-factorial such that the non kawamata log terminal locus is 
a divisor:
\begin{corollary}\label{c_model} Let $(X,\Delta)$ be a log pair.   

Then there is a birational morphism $\pi\colon\map Y.X.$, where $Y$ is
$\mathbb{Q}$-factorial, such that if we write
$$
K_Y+\Gamma=K_Y+\Gamma_1+\Gamma_2=\pi^*(K_X+\Delta),
$$
where every component of $\Gamma_1$ has coefficient less than one and every component of
$\Gamma_2$ has coefficient at least one, then $K_Y+\Gamma_1$ is kawamata log terminal and
nef over $X$ and no component of $\Gamma_1$ is exceptional.  
\end{corollary}

Even thought the result in \eqref{c_extract} is not optimal as it does not fully address
the log canonical case, nevertheless, we are able to prove the following result (cf.
\cite{Shokurov93}, \cite{Kollaretal}, \cite{Kawakita05}):
\begin{corollary}[Inversion of adjunction]\label{c_inversion} Let $(X,\Delta)$ be a 
log pair and let $\nu\colon\map S.S'.$ be the normalisation of a component $S'$ of
$\Delta$ of coefficient one.

If we define $\Theta$ by adjunction,
$$
\nu^*(K_X+\Delta)=K_S+\Theta, 
$$
then the log discrepancy of $K_S+\Theta$ is equal to the minimum of the log
discrepancy with respect to $K_X+\Delta$ of any valuation whose centre on $X$ is of
codimension at least two.
\end{corollary}

One of the most compelling reasons to enlarge the category of varieties to the category of
algebraic spaces (equivalently Moishezon spaces, at least in the proper case) is to allow
the possibility of cut and paste operations, such as one can perform in topology.
Unfortunately, it is then all too easy to construct proper smooth algebraic spaces over
$\mathbb{C}$, which are not projective.  In fact the appendix to \cite{Hartshorne77} has
two very well known examples due to Hironaka.  In both examples, one exploits the fact
that for two curves in a threefold which intersect in a node, the order in which one blows
up the curves is important (in fact the resulting threefolds are connected by a flop).

It is then natural to wonder if this is the only way to construct such examples, in the
sense that if a proper algebraic space is not projective then it must contain a rational
curve.  Koll\'ar dealt with the case when $X$ is a terminal threefold with Picard number
one, see \cite{Kollar91b}.  In a slightly different but related direction, it is
conjectured that if a complex K\"ahler manifold $M$ does not contain any rational curves
then $K_M$ is nef (see for example \cite{Peternell98}), which would extend some of Mori's
famous results from the projective case.  Koll\'ar also has some unpublished proofs of
some related results.  

The following result, which was proved by Shokurov assuming the existence and termination
of flips, cf. \cite{Shokurov97}, gives an affirmative answer to the first conjecture and
at the same time connects the two conjectures:
\begin{corollary}\label{c_moishezon} Let $\pi\colon\map X.U.$ be a proper map 
of normal algebraic spaces, where $X$ is analytically $\mathbb{Q}$-factorial.   

If $K_X+\Delta$ is divisorially log terminal and $\pi$ does not contract any rational
curves then $\pi$ is a log terminal model.  In particular $\pi$ is projective and
$K_X+\Delta$ is $\pi$-nef.
\end{corollary}

\makeatletter
\renewcommand{\thetheorem}{\thesection.\arabic{theorem}}
\@addtoreset{theorem}{section} 
\makeatother

\section{Description of the proof}

\begin{theorema}[Existence of pl-flips]\label{t_existence} Let $f\colon\map X.Z.$ be
a pl-flipping contraction for an $n$-dimensional purely log terminal pair $(X,\Delta)$.

Then the flip $f^+\colon\map X^+.Z.$ of $f$ exists.
\end{theorema}

\begin{theorema}[Special finiteness]\label{t_special} Let $\pi\colon\map X.U.$ be a 
projective morphism of normal quasi-projective varieties, where $X$ is $\mathbb
Q$-factorial of dimension $n$.  Let $V$ be a finite dimensional affine subspace of
$\WDiv_{\mathbb{R}}(X)$, which is defined over the rationals, let $S$ be the sum of
finitely many prime divisors and let $A$ be a general ample $\mathbb{Q}$-divisor over $U$.
Let $(X,\Delta_0)$ be a divisorially log terminal pair such that $S\leq \Delta_0$.  Fix a
finite set $\mathfrak{E}$ of prime divisors on $X$.

Then there are finitely $1\leq i\leq k$ many birational maps $\phi_i\colon\rmap X.Y_i.$  over $U$
such that if $\phi\colon\rmap X.Y.$ is any $\mathbb{Q}$-factorial weak log
canonical model over $U$ of $K_{X}+\Delta$ where
$\Delta\in \mathcal{L}_{S+A}(V)$ which only contracts elements of
$\mathfrak{E}$ and which does not contract every component of $S$, then there 
is an index $1\leq i\leq k$ such that the induced birational map $\xi\colon\rmap Y_i.Y.$ 
is an isomorphism in a neighbourhood of the strict transforms of $S$.  
\end{theorema}

\begin{theorema}[Existence of log terminal models]\label{t_model} Let $\pi\colon\map X.U.$ 
be a projective morphism of normal quasi-projective varieties, where $X$ has dimension
$n$.  Suppose that $K_X+\Delta$ is kawamata log terminal, where $\Delta$ is big over $U$.

If there exists an $\mathbb{R}$-divisor $D$ such that $K_X+\Delta\sim_{\mathbb{R},U}
D\geq 0$, then $K_X+\Delta$ has a log terminal model over $U$.
\end{theorema} 

\begin{theorema}[Non-vanishing theorem]\label{t_effective} Let $\pi\colon\map X.U.$ be 
a projective morphism of normal quasi-projective varieties, where $X$ has dimension $n$.
Suppose that $K_X+\Delta$ is kawamata log terminal, where $\Delta$ is big over $U$.

If $K_X+\Delta$ is $\pi$-pseudo-effective, then there exists an $\mathbb{R}$-divisor $D$
such that $K_X+\Delta\sim_{\mathbb{R},U} D\geq 0$.
\end{theorema}

\begin{theorema}[Finiteness of models]\label{t_finite}  Let $\pi\colon\map X.U.$ be 
a projective morphism of normal quasi-projective varieties, where $X$ has dimension $n$.
Fix a general ample $\mathbb{Q}$-divisor $A\geq 0$ over $U$.  Let $V$ be a finite
dimensional affine subspace of $\WDiv _{\mathbb{R}}(X)$ which is defined over the
rationals.  Suppose that there is a kawamata log terminal pair $(X,\Delta_0)$.

Then there are finitely many birational maps $\psi_j\colon\rmap X.Z_j.$ over $U$, $1\leq
j\leq l$ such that if $\psi\colon\rmap X.Z.$ is a weak log canonical model of $K_X+\Delta$
over $U$, for some $\Delta\in\mathcal{L}_A(V)$ then there is an index $1\leq j\leq l$ and
an isomorphism $\xi\colon\map Z_j.Z.$ such that $\psi=\xi\circ\psi_j$.
\end{theorema}

\begin{theorema}[Finite generation]\label{t_ezd} Let $\pi\colon\map X.Z.$ 
be a projective morphism to a normal affine variety.  Let $(X,\Delta=A+B)$ be a kawamata
log terminal pair of dimension $n$, where $A\geq 0$ is an ample $\mathbb{Q}$-divisor and
$B\geq 0$.  If $K_X+\Delta$ is pseudo-effective, then
\begin{enumerate} 
\item The pair $(X,\Delta)$ has a log terminal model $\mu\colon\rmap X.Y.$.  
In particular if $K_X+\Delta$ is $\mathbb{Q}$-Cartier then the log canonical ring 
$$
R(X,K_X+\Delta)=\bigoplus_{m\in\mathbb{N}}H^0(X,\ring X.(\rdown m(K_X+\Delta).)),
$$
is finitely generated. 
\item Let $V\subset \WDiv_{\mathbb{R}}(X)$ be the vector space spanned by the components
of $\Delta$.  Then there is a constant $\delta>0$ such that if $G$ is a prime divisor
contained in the stable base locus of $K_X+\Delta$ and $\Xi\in \mathcal L _{A}(V)$ such
that $\|\Xi-\Delta\|<\delta$, then $G$ is contained in the stable base locus of $K_X+\Xi$.
\item Let $W\subset V$ be the smallest affine subspace of $\WDiv_{\mathbb{R}}(X)$
containing $\Delta$, which is defined over the rationals.  Then there is a constant
$\eta>0$ and a positive integer $r>0$ such that if $\Xi\in W$ is any divisor and $k$ is
any positive integer such that $\|\Xi-\Delta\|<\eta$ and $k(K_X+\Xi)/r$ is Cartier, then
every component of $\fix (k(K_X+\Xi))$ is a component of the stable base locus of
$K_X+\Delta$.
\end{enumerate} 
\end{theorema}

The proof of Theorem~\ref{t_existence}, Theorem~\ref{t_special}, Theorem~\ref{t_model},
Theorem~\ref{t_effective}, Theorem~\ref{t_finite} and Theorem~\ref{t_ezd} proceeds by
induction:
\begin{itemize} 
\item Theorem~\ref{t_ezd}$_{n-1}$ implies Theorem~\ref{t_existence}$_n$, see the main result of \cite{HM07}.  
\item Theorem~\ref{t_finite}$_{n-1}$ implies Theorem~\ref{t_special}$_n$, cf. \eqref{l_A-to-B}.
\item Theorem~\ref{t_existence}$_{n}$ and Theorem~\ref{t_special}$_n$
imply Theorem~\ref{t_model}$_n$, cf. \eqref{l_B-to-C}.
\item Theorem~\ref{t_effective}$_{n-1}$, Theorem~\ref{t_special}$_{n}$ and
Theorem~\ref{t_model}$_n$ imply Theorem~\ref{t_effective}$_n$, cf.  \eqref{l_C-to-D}.
\item Theorem~\ref{t_model}$_n$ and Theorem~\ref{t_effective}$_n$ imply Theorem~\ref{t_finite}$_n$, cf. \eqref{l_D-to-E}.  
\item Theorem~\ref{t_model}$_n$, Theorem~\ref{t_effective}$_n$ and
Theorem~\ref{t_finite}$_n$ imply Theorem~\ref{t_ezd}$_n$, cf. \eqref{l_E-to-F}.
\end{itemize}

\subsection{Sketch of the proof}

To help the reader navigate through the technical problems which naturally arise when
trying to prove \eqref{t_main}, we review a natural approach to proving that the canonical
ring
$$
R(X,K_X)=\bigoplus_{m\in\mathbb{N}} H^0(X,\ring X.(mK_X)),
$$
of a smooth projective variety $X$ of general type is finitely generated.  Even though we
do not directly follow this method to prove the existence of log terminal models, instead
using ideas from the MMP, many of the difficulties which arise in our approach are
mirrored in trying to prove finite generation directly.

A very natural way to proceed is to pick a divisor $D\in |kK_X|$, whose existence is
guaranteed as we are assuming that $K_X$ is big, and then to restrict to $D$.  One obtains
an exact sequence
$$
\les {H^0(X,\ring X.((l-k)K_X))}.{H^0(X,\ring X.(lK_X))}.{H^0(D,\ring D.(mK_D))}.,
$$
where $l=m(1+k)$ is divisible by $k+1$, and it is easy to see that it suffices to prove
that the restricted algebra, given by the image of the maps
$$
\map {H^0(X,\ring X.(m(1+k)K_X))}.{H^0(D,\ring D.(mK_D))}.,
$$
is finitely generated.  Various problems arise at this point.  First $D$ is neither
smooth nor even reduced (which, for example, means that the symbol $K_D$ is only formally
defined; strictly speaking we ought to work with the dualising sheaf $\omega_D$).  It is
natural then to pass to a log resolution, so that the support $S$ of $D$ has simple normal
crossings, and to replace $D$ by $S$.  The second problem is that the kernel of the map
$$
\map {H^0(X,\ring X.(m(1+k)K_X))}.{H^0(S,\ring S.(mK_S))}.,
$$
no longer has any obvious connection with 
$$
H^0(X,\ring X.((m(1+k)-k)K_X)),
$$ 
so that even if we knew that the new restricted algebra were finitely generated, it is not
immediate that this is enough.  Another significant problem is to identify the restricted
algebra as a subalgebra of
$$
\bigoplus_{m\in\mathbb{N}}H^0(S,\ring S.(mK_S)),
$$ 
since it is only the latter that we can handle by induction.  Yet another problem is that
if $C$ is a component of $S$, it is no longer the case that $C$ is of general type, so
that we need a more general induction.  In this case the most significant problem to deal
with is that even if $K_C$ is pseudo-effective, it is not clear that the linear system
$|kK_C|$ is non-empty for any $k>0$.  Finally, even though this aspect of the problem may
not be apparent from the description above, in practise it seems as though we need to work
with infinitely many different values of $k$ and hence $D=D_k$, which entails working with
infinitely many different birational models of $X$ (since for every different value of
$k$, one needs to resolve the singularities of $D$).

Let us consider one special case of the considerations above, which will hopefully throw
some more light on the problem of finite generation.  Suppose that to resolve the
singularities of $D$ we need to blow up a subvariety $V$.  The corresponding divisor $C$
will typically fibre over $V$ and if $V$ has codimension two then $C$ will be close to a
$\pr 1.$-bundle over $V$.  In the best case, the projection $\pi\colon\map C.V.$ will be a
$\pr 1.$-bundle with two disjoint sections (this is the toroidal case) and sections of
tensor powers of a line bundle on $C$ will give sections of an algebra on $V$ which is
graded by $\mathbb{N}^2$, rather than just $\mathbb{N}$.  Let us consider then the
simplest possible algebras over $\mathbb{C}$ which are graded by $\mathbb{N}^2$.  If we
are given a submonoid $M\subset \mathbb{N}^2$ (that is a subset of $\mathbb{N}^2$ which
contains the origin and is closed under addition) then we get a subalgebra
$R\subset\mathbb{C}^2[x,y]$ spanned by the monomials
$$
\{\, x^iy^j \,|\, (i,j)\in M \,\}.  
$$
The basic observation is that $R$ is finitely generated iff $M$ is a finitely generated
monoid.  There are two obvious cases when $M$ is not finitely generated, 
$$
M=\{\, (i,j)\in\mathbb{N}^2 \,|\, j>0 \,\}\cup \{(0,0)\} \quad \text{and} \quad \{\,(i,j)\in\mathbb{N}^2 \,|\, i\geq \sqrt{2}j \,\}.  
$$
In fact, if $\mathcal{C}\subset\mathbb{R}^2$ is the convex hull of the set $M$, then $M$
is finitely generated iff $\mathcal{C}$ is a rational polytope.  In the general case, we
will be given a convex subset $\mathcal{C}$ of a finite dimensional vector space of Weil
divisors on $X$ and a key part of the proof is to show that the set $\mathcal{C}$ is in
fact a rational polytope.  As naive as these examples are, hopefully they indicate why it
is central to the proof of finite generation to
\begin{itemize} 
\item consider divisors with real coefficients, and
\item prove a non-vanishing result.  
\end{itemize} 

We now review our approach to the proof of \eqref{t_main}.  As is clear from the plan of
the proof given in the previous subsection, the proof of \eqref{t_main} is by induction on
the dimension and the proof is split into various parts.  Instead of proving directly that
the canonical ring is finitely generated, we try to construct a log terminal model for
$X$.  The first part is to prove the existence of pl-flips.  This is proved by induction
in \cite{HM05c}, and we will not talk about the proof of this result here, since the
methods used to prove this result are very different from the methods we use here.
Granted the existence of pl-flips, the main issue is to prove that some MMP terminates,
which means that we must show that we only need finitely many flips.

As in the scheme of the proof of finite generation sketched above, the first step is to
pick $D\in |kK_X|$, and to pass to a log resolution of the support $S$ of $D$.  By way of
induction we want to work with $K_X+S$ rather than $K_X$.  As before this is tricky since
a log terminal model for $K_X+S$ is not the same as a log terminal model for $K_X$.  In
other words having added $S$, we really want to subtract it as well.  The trick however is
to first add $S$, construct a log terminal model for $K_X+S$ and then subtract $S$ (almost
literally component by component).  This is one of the key steps, to show that
Theorem~\ref{t_existence}$_n$ and Theorem~\ref{t_special}$_n$ imply
Theorem~\ref{t_model}$_n$.  This part of the proof splits naturally into two parts.  First
we have to prove that we may run the relevant minimal model programs, see \S
\ref{s_special} and the beginning of \S \ref{s_models}, then we have to prove this does
indeed construct a log terminal model for $K_X$, see \S \ref{s_models}.

To gain intuition for how this part of the proof works, let us first consider a simplified
case.  Suppose that $D=S$ is irreducible.  In this case it is clear that $S$ is of general
type and $K_X$ is nef if and only if $K_X+S$ is nef and in fact a log terminal model for
$K_X$ is the same as a log terminal model for $K_X+S$.  Consider running the
$(K_X+S)$-MMP.  Then every step of this MMP is a step of the $K_X$-MMP and vice-versa.
Suppose that we have a $(K_X+S)$-extremal ray $R$.  Let $\pi\colon\map X.Z.$ be the
corresponding contraction.  Then $S\cdot R<0$, so that every curve $\Sigma$ contracted by
$\pi$ must be contained in $S$.  If $\pi$ is a divisorial contraction, then $\pi$ must
contract $S$ and $K_Z\sim _{\mathbb{Q}} 0$, so that $Z$ is a log terminal model.
Otherwise $\pi$ is a pl-flip and by Theorem~\ref{t_existence}$_n$, we can construct the
flip of $\pi$, $\phi\colon\rmap X.Y.$.  Consider the restriction $\psi\colon\rmap S.T.$ of
$\phi$ to $S$, where $T$ is the strict transform of $S$.  Since log discrepancies increase
under flips and $S$ is irreducible, $\psi$ is a birational contraction.  After finitely
many flips, we may therefore assume that $\psi$ does not contract any divisors, since the
Picard number of $S$ cannot keep dropping.  Consider what happens if we restrict to $S$.
By adjunction, we have
$$
(K_X+S)|_S=K_S.
$$
Thus $\psi\colon\rmap S.T.$ is $K_S$-negative.  We have to show that this cannot happen
infinitely often.  If we knew that every sequence of flips on $S$ terminates, then we
would be done.  In fact this is how special termination works.  Unfortunately we cannot
prove that every sequence of flips terminates on $S$, so that we have to do something
slightly different.  Instead we throw in an auxiliary ample divisor $H$ on $X$, and
consider $K_X+S+tH$, where $t$ is a positive real number.  If $t$ is large enough then
$K_X+S+tH$ is ample.  Decreasing $t$, we may assume that there is an extremal ray $R$ such
that $(K_X+S+tH)\cdot R=0$.  If $t=0$, then $K_X+S$ is nef and we are done.  Otherwise
$(K_X+S)\cdot R<0$, so that we are still running a $(K_X+S)$-MMP, but with the additional
restriction that $K_X+S+tH$ is nef and trivial on any ray we contract.  This is the
$(K_X+S)$-MMP with scaling of $H$.  Let $G=H|_S$.  Then $K_S+tG$ is nef and so is
$K_T+tG'$, where $G'=\psi_*G$.  In this case $K_T+tG'$ is a weak log canonical model for
$K_S+tG$ (it is not a log terminal model, both because $\phi$ might contract divisors on
which $K_S+tG$ is trivial and more importantly because $T$ need not be
$\mathbb{Q}$-factorial).  In this case we are then done, by finiteness of weak log
canonical models for $(S,tG)$, where $t\in [0,1]$ (cf. Theorem~\ref{t_finite}$_{n-1}$).

We now turn to the general case.  The idea is similar.  First we want to use finiteness of
log terminal models on $S$ to conclude that there are only finitely many log terminal
models in a neighbourhood of $S$.  Secondly we use this to prove the existence of a very
special MMP and construct log terminal models using this MMP.  The intuitive idea is that
if $\phi\colon\rmap X.Y.$ is $K_X$-negative then $K_X$ is bigger than $K_Y$ (the
difference is an effective divisor on a common resolution) so that we can never return to
the same neighbourhood of $S$.  As already pointed out, in the general case we need to
work with $\mathbb{R}$-divisors.  This poses no significant problem at this stage of the
proof, but it does make some of the proofs a little more technical.  By way of induction,
suppose that we have a log pair $K_X+\Delta=K_X+S+A+B$, where $S$ is a sum of prime
divisors, $A$ is an ample divisor (with rational coefficients) and the coefficients of
$\Delta$ are real numbers between zero and one.  We are also given a divisor $D\geq 0$
such that $K_X+\Delta \sim_{\mathbb{R}} D$.  The construction of log terminal models is
similar to the one sketched above and breaks into two parts.

In the first part, for simplicitly of exposition we assume that $S$ is a prime divisor and
$K_X+\Delta$ is purely log terminal.  We fix $S$ and $A$ but we allow $B$ to vary and we
want to show that finiteness of log terminal models for $S$ implies finiteness of log
terminal models in a neighbourhood of $S$.  We are free to pass to a log resolution, so we
may assume that $(X,\Delta)$ is log smooth and if $B=\sum b_iB_i$ the coefficients
$(\llist b.k.)$ of $B$ lie in $[0,1]^k$.  Let $\Theta=(\Delta-S)|_S$ so that
$(K_X+\Delta)|_S=K_S+\Theta$.

Suppose that $f\colon\rmap X.Y.$ is a log terminal model of $(X,\Delta)$.  There are three
problems that arise, two of which are quite closely related.  Suppose that $g\colon\rmap
S.T.$ is the restriction of $f$ to $S$, where $T$ is the strict transform of $S$.  The
first problem is that $g$ need not be a birational contraction.  For example, suppose that
$X$ is a threefold and $f$ flips a curve $\Sigma$ intersecting $S$, which is not contained
in $S$.  Then $S\cdot \Sigma>0$ so that $T\cdot E<0$, where $E$ is the flipped curve.  In
this case $E\subset T$ so that the induced birational map $\rmap S.T.$ extracts the curve
$E$.  The basic observation is that $E$ must have log discrepancy less than one with
respect to $(S,\Theta)$.  Since the pair $(X,\Delta)$ is purely log terminal if we replace
$(X,\Delta)$ by a fixed model which is high enough then we can ensure that the pair
$(S,\Theta)$ is terminal, so that there are no such divisors $E$ and $g$ is then always a
birational contraction.  The second problem is that if $E$ is a divisor intersecting $S$
which is contracted to a divisor lying in $T$, then $E\cap S$ is not contracted by $S$.
For this reason, $g$ is not necessarily a weak log canonical model of $(S,\Theta)$.
However we can construct a divisor $0\leq \Xi\leq \Theta$ such that $g$ is a weak log
canonical model for $(S,\Xi)$.  Suppose that we start with a smooth threefold $Y$ and a
smooth surface $T\subset Y$ which contains a $-2$-curve $\Sigma$, such that $K_Y+T$ is
nef.  Let $f\colon\map X.Y.$ be the blow up of $Y$ along $\Sigma$ with exceptional divisor
$E$ and let $S$ be the strict transform of $T$.  Then $f$ is a step of the
$(K_X+S+eE)$-MMP for any $e>0$ and $f$ is a log terminal model of $K_X+S+eE$.  The
restriction of $f$ to $S$, $g\colon\map S.T.$ is the identity, but $g$ is not a log
terminal model for $K_S+e\Sigma$, since $K_S+e\Sigma$ is negative along $\Sigma$.  It is a
weak log canonical model for $K_T$, so that in this case $\Xi=0$.  The details of the
construction of $\Xi$ are contained in \eqref{l_terminal-negative}.

The third problem is that the birational contraction $g$ does not determine $f$.  This is
most transparent in the case when $X$ is a surface and $S$ is a curve, since in this case
$g$ is always an isomorphism.  To remedy this particular part of the third problem we use
the different, which is defined by adjunction,
$$
(K_Y+T)|_T=K_T+\Phi.
$$
The other parts of the third problem only occurs in dimension three or more.  For example,
suppose that $Z$ is the cone over a smooth quadric in $\pr 3.$ and $p\colon\map X.Z.$ and
$q\colon\map Y.Z.$ are the two small resolutions, so that the induced birational map
$f\colon\rmap X.Y.$ is the standard flop.  Let $\pi\colon\map W.Z.$ blow up the maximal
ideal, so that the exceptional divisor $E$ is a copy of $\pr 1.\times\pr 1.$.  Pick a
surface $R$ which intersects $E$ along a diagonal curve $\Sigma$.  If $S$ and $T$ are the
strict transforms of $R$ in $X$ and $Y$ then the induced birational map $g\colon\map S.T.$
is an isomorphism (both $S$ and $T$ are isomorphic to $R$).  To get around this problem,
one can perturb $\Delta$ so that $f$ is the ample model, and one can distinguish between
$X$ and $Y$ by using the fact that $g$ is the ample model of $(S,\Xi)$.  Finally it is
not hard to write down examples of flops which fix $\Xi$, but switch the individual
components of $\Xi$.  In this case one needs to keep track not only of $\Xi$ but the
individual pieces $(g_*B_i)|_T$, $1\leq i\leq k$.  We prove that an ample model $f$ is
determined in a neighbourhood of $T$ by $g$, the different $\Phi$ and $(f_*B_i)|_T$, see
\eqref{l_determine-neighbourhood}.  To finish this part, by induction we assume that there
are finitely many possibilities for $g$ and it is easy to see that there are then only
finitely many possibilities for the different $\Phi$ and the divisors $(f_*B_i)|_T$ and
this shows that there are only finitely many possibilities for $f$.  The details are
contained in \S \ref{s_special}.

The second part consists of using finiteness of models in a neighbourhood of $S$ to run a
sequence of minimal model programs to construct a log terminal model.  We may assume that
$X$ is smooth and the support of $\Delta+D$ has normal crossings.  

Suppose that there is a divisor $C$ such that
\[
K_X+\Delta \sim_{\mathbb{R},U} D+\alpha C, \tag{$\ast$}   
\]
where $K_X+\Delta+C$ is divisorially log terminal and nef and the support of $D$ is
contained in $S$.  If $R$ is an extremal ray which is $(K_X+\Delta)$-negative, then
$D\cdot R<0$, so that $S_i\cdot R<0$ for some component $S_i$ of $S$.  As before this
guarantees the existence of flips.  It is easy to see that the corresponding step of the
$(K_X+\Delta)$-MMP is not an isomorphism in a neighbourhood of $S$.  Therefore the
$(K_X+\Delta)$-MMP with scaling of $C$ must terminate with a log terminal model for
$K_X+\Delta$.  To summarise, whenever the conditions above hold, we can always construct a
log terminal model of $K_X+\Delta$.  

We now explain how to construct log terminal models in the general case.  We may write
$D=D_1+D_2$, where every component of $D_1$ is a component of $S$ and no component of
$D_2$ is a component of $S$.  If $D_2$ is empty, that is every component of $D$ is a
component of $S$ then we take $C$ to be a sufficiently ample divisor, and the argument in
the previous paragraph implies that $K_X+\Delta$ has a log terminal model.  If $D_2\neq 0$
then instead of constructing a log terminal model, we argue that we can construct a
neutral model, which is exactly the same as a log terminal model, except that we drop the
hypothesis on negativity.  Consider $(X,\Theta=\Delta+\lambda D_2)$, where $\lambda$ is
the largest real number so that the coefficients of $\Theta$ are at most one.  Then more
components of $\rdown \Theta.$ are components of $D$.  By induction $(X,\Theta)$ has a
neutral model, $f\colon\rmap X.Y.$.  It is then easy to check the conditions in the
paragraph above apply, and we can construct a log terminal model $g\colon\rmap Y.Z.$ for
$K_X+g_*\Delta$.  It is then automatic that the composition $h=g\circ f\colon\rmap X.Z.$
is a neutral model of $K_X+\Delta$ (since $f$ is not $(K_X+\Delta)$-negative it is not true
in general that $h$ is a log terminal model of $K_X+\Delta$).  However $g$ is
automatically a log terminal model model provided we only contract components of the
stable base locus of $K_X+\Delta$.  For this reason, we pick $D$ so that we may write
$D=M+F$ where every component of $M$ is semiample and every component of $F$ is a
component of the stable base locus.  This explains the implication
Theorem~\ref{t_finite}$_{n-1}$ implies Theorem~\ref{t_special}$_n$.  The details are
contained in \S \ref{s_special} and \S \ref{s_models}.

Now we explain how prove that if $K_X+\Delta=K_X+A+B$ is pseudo-effective, then
$K_X+\Delta\sim_{\mathbb{R}}D\geq 0$.  The idea is to mimic the proof of the non-vanishing
theorem.  As in the proof of the non-vanishing theorem and following the work of Nakayama,
there are two cases.  In the first case, for any ample divisor $H$,
$$
h^0(X,\ring X.(\rdown m(K_X+\Delta).+H)),
$$
is a bounded function of $m$. In this case it follows that $K_X+\Delta$ is numerically
equivalent to the divisor $N_{\sigma}(K_X+\Delta)\geq 0$.  It is then not hard to prove
that Theorem~\ref{t_model}$_n$ implies that $K_X+\Delta$ has a log terminal model and we
are done by the base point free theorem.

In the second case we construct a non kawamata log terminal centre for 
$$
m(K_X+\Delta)+H,
$$ 
when $m$ is sufficiently large.  Passing to a log resolution, and using standard
arguments, we are reduced to the case when 
$$
K_X+\Delta=K_X+S+A+B,
$$ 
where $S$ is irreducible and $(K_X+\Delta)|_S$ is pseudo-effective, and the support of
$\Delta$ has global normal crossings.  Suppose first that $K_X+\Delta$ is
$\mathbb{Q}$-Cartier.  We may write
$$
(K_X+S+A+B)|_S=K_S+C+D,
$$
where $C$ is ample and $D\geq 0$.  By induction we know that there is a positive integer
$m$ such that $h^0(S,\ring S.(m(K_S+C+D)))>0$.  To lift sections, we need to know that
$h^1(X,\ring X.(m(K_X+S+A+B)-S))=0$.  Now
\begin{align*}
m(K_X+\Delta)-(K_X+B)-S&=(m-1)(K_X+\Delta)+A \\
                      &=(m-1)(K_X+\Delta+\frac 1{m-1}A).
\end{align*}
As $K_X+\Delta+A/(m-1)$ is big, we can construct a log terminal model $\phi\colon\rmap
X.Y.$ for $K_X+\Delta+A/(m-1)$, and running this argument on $Y$, the required vanishing
holds by Kawamata-Viehweg vanishing.  In the general case, $K_X+S+A+B$ is an
$\mathbb{R}$-divisor.  The argument is now a little more delicate as $h^0(S,\ring
S.(m(K_S+C+D)))$ does not make sense.  We need to approximate $K_S+C+D$ by rational
divisors, which we can do by induction.  But then it is not so clear how to choose $m$.
In practise we need to prove that the log terminal model $Y$ constructed above does not
depend on $m$, at least locally in a neighbourhood of $T$, the strict transform of $S$,
and then the result follows by Diophantine approximation.  This explains the implication
Theorem~\ref{t_effective}$_{n-1}$, Theorem~\ref{t_finite}$_{n-1}$ and
Theorem~\ref{t_model}$_n$ imply Theorem~\ref{t_effective}$_n$.  The details are in \S
\ref{s_effective}.

Finally, in terms of induction, we need to prove finiteness
of weak log canonical models.  We fix an ample divisor $A$ and work with divisors of the
form $K_X+\Delta=K_X+A+B$, where the coefficients of $B$ are variable.  For ease of
exposition, we assume that the support of $A$ and $B$ have global normal crossings, so
that $K_X+\Delta=K_X+A+\sum b_iB_i$ is log canonical if and only if $0\leq b_i\leq 1$ for
all $i$.  The key point is that we allow the coefficients of $B$ to be real numbers, so
that the set of all possible choices of coefficients $[0,1]^k$ is a compact subset of
$\mathbb{R}^k$.  Thus we may check finiteness locally.  In fact since $A$ is ample, we can
always perturb the coefficients of $B$ so that none of the coefficients are equal to one
or zero and so we may even assume that $K_X+\Delta$ is kawamata log terminal.

Observe that we are certainly free to add components to $B$ (formally we add components
with coefficient zero and then perturb so that their coefficients are non-zero).  In
particular we may assume that $B$ is the support of an ample divisor and so working on the
weak log canonical model, we may assume that we have a log canonical model for a perturbed
divisor.  Thus it suffices to prove that there are only finitely many log canonical
models.  Since the log canonical model is determined by any log terminal model, it
suffices to prove that we can find a cover of $[0,1]^k$ by finitely many log terminal
models.  By compactness, it suffices to do this locally.

So pick $b\in [0,1]^k$.  There are two cases.  If $K_X+\Delta$ is not pseudo-effective,
then $K_X+A+B'$ is not pseudo-effective, for $B'$ in a neighbourhood of $B$, and there are
no weak log canonical models at all.  Otherwise we may assume that $K_X+\Delta$ is
pseudo-effective.  By induction we know that $K_X+\Delta\sim_{\mathbb{R}} D\geq 0$.  Then
we know that there is a log terminal model $\phi\colon\rmap X.Y.$.  Replacing $(X,\Delta)$
by $(Y,\Gamma=\phi_*\Delta)$, we may assume that $K_X+\Delta$ is nef.  By the base point
free theorem it is semiample.  Let $\map X.Z.$ be the corresponding morphism.  The key
observation is that locally about $\Delta$, any log terminal model over $Z$ is an absolute
log terminal model.  Working over $Z$, we may assume that $K_X+\Delta$ is numerically
trivial.  In this case the problem of finding a log terminal model for $K_X+\Delta'$ only
depends on the line segment spanned by $\Delta$ and $\Delta'$.  Working in a small box
about $\Delta$, we are then reduced to finding a log terminal model on the boundary of the
box and we are done by induction on the dimension of the affine space containing $B$.
Note that in practise, we need to work in slightly more generality than we have indicated;
first we need to work in the relative setting and secondly we need to work with an
arbitrary affine space containing $B$ (and not just the space spanned by the components of
$B$).  This poses no significant problem.  This explains the implication
Theorem~\ref{t_model}$_n$ and Theorem~\ref{t_effective}$_n$ imply
Theorem~\ref{t_finite}$_n$.  The details are contained in \S \ref{s_finite}.

Let us end the sketch of the proof by pointing out some of the technical advantages with
working with kawamata log terminal pairs $(X,\Delta)$, where $\Delta$ is big.  The first
observation is that since the kawamata log terminal condition is open, it is
straightforward to show that $\Delta$ is $\mathbb{Q}$-linearly equivalent to $A+B$ where
$A$ is an ample $\mathbb{Q}$-divisor, $B\geq 0$ and $K_X+A+B$.  The presence of the ample
divisor $A$ is very convenient for a number of reasons, two of which we have already seen
in the sketch of the proof.

Firstly the restriction of an ample divisor to any divisor $S$ is ample, so that if $B$
does not contain $S$ in its support then the restriction of $A+B$ to $S$ is big.  
This is very useful for induction.  

Secondly, as we vary the coefficients of $B$, the closure of the set of kawamata log
terminal pairs is the set of log canonical pairs.  However, we can use a small piece of
$A$ to perturb the coefficients of $B$ so that they are bounded away from zero and
$K_X+A+B$ is always kawamata log terminal. 

Finally, if $(X,\Delta)$ is divisorially log terminal and $f\colon\map X.Y.$ is a
$(K_X+\Delta)$-trivial contraction then $K_Y+\Gamma=K_X+f_*\Delta$ is not necessarily
divisorially log terminal, only log canonical.  For example, suppose that $Y$ is a surface
with a simple elliptic singularity and $f\colon\map X.Y.$ is the blow up with exceptional
divisor $E$.  Then $f$ is a weak log canonical model of $K_X+E$, but $Y$ is not log
terminal as it does not have rational singularities.  On the other hand, if $\Delta=A+B$ where
$A$ is ample then $K_Y+\Gamma$ is always divisorially log terminal.

\subsection{Standard conjectures of the MMP}

Having sketched the proof of \eqref{t_main}, we should point out the main obstruction to
extending these ideas to the case when $X$ is not of general type.  The main issue seems
to be the implication $K_X$ pseudo-effective implies $\kappa(X,K_X)\geq 0$.  In other words 
we need:

\begin{conjecture}\label{c_pseudo} Let $(X,\Delta)$ be a projective kawamata log terminal 
pair.  

 If $K_X+\Delta$ is pseudo-effective then $\kappa(X,K_X+\Delta)\geq 0$.
\end{conjecture}

 We also probably need 
\begin{conjecture}\label{c_non} Let $(X,\Delta)$ be a projective kawamata log terminal 
pair.  

 If $K_X+\Delta$ is pseudo-effective and
$$
h^0(X,\ring X.(\rdown m(K_X+\Delta).+H)),
$$
is not a bounded function of $m$, for some ample divisor $H$, then $\kappa(X,K_X+\Delta)\geq 1$.
\end{conjecture}

In fact, using the methods of this paper, together with some results of
Kawamata (cf. \cite{Kawamata87} and \cite{Kawamata88}), \eqref{c_pseudo} and \eqref{c_non}
would seem to imply one of the main outstanding conjectures of higher dimensional 
geometry:
\begin{conjecture}[Abundance]\label{c_abundance} Let $(X,\Delta)$ be a projective kawamata log terminal pair. 

 If $K_X+\Delta$ is nef then it is semiample.  
\end{conjecture}
We remark that the following seemingly innocuous generalisation of
\eqref{t_main} (in dimension $n$) would seem to imply \eqref{c_abundance}
(in dimension $n+1$).
\begin{conjecture}\label{c_canonical} Let $(X,\Delta)$ be a projective log canonical pair
of dimension $n$.  

If $K_X+\Delta$ is big, then $(X,\Delta)$ has a log canonical model.  
\end{conjecture}

 It also seems worth pointing out that the other remaining conjecture is:
\begin{conjecture}[Borisov-Alexeev-Borisov]\label{c_ab} Fix a positive integer 
$n$ and a positive real number $\epsilon>0$.  

Then the set of varieties $X$ such that $K_X+\Delta$ has log discrepancy at least 
$\epsilon$ and $-(K_X+\Delta)$ is ample, forms a bounded family.   
\end{conjecture}

\section{Preliminary results}
\label{s_basics}

In this section we collect together some definitions and results that will be needed for
the proof of \eqref{c_scaling}.

\makeatletter
\renewcommand{\thetheorem}{\thesubsection.\arabic{theorem}}
\@addtoreset{theorem}{subsection} 
\makeatother

\subsection{Notation and conventions}

We work over the field of complex numbers $\mathbb{C}$.  We say that two
$\mathbb{Q}$-divisors $D_1$, $D_2$ are \textit{$\mathbb{Q}$-linearly equivalent} ($D_1\sim
_{\mathbb{Q}} D_2$) if there exists an integer $m>0$ such that $mD_i$ are linearly
equivalent.  We say that a $\mathbb{Q}$-divisor $D$ is \textit{$\mathbb{Q}$-Cartier} if
some integral multiple is Cartier.  We say that $X$ is \textit{$\mathbb{Q}$-factorial} if
every Weil divisor is $\mathbb{Q}$-Cartier.  We say that $X$ is \textit{analytically
  $\mathbb{Q}$-factorial} if every analytic Weil divisor (that is an analytic subset of
codimension one) is analytically $\mathbb{Q}$-Cartier (i.e. some multiple is locally
defined by a single analytic function).  We recall some definitions involving divisors
with real coefficients.
   
\begin{definition}\label{d_divisor} Let $\pi\colon\map X.U.$ be a proper morphism of 
normal algebraic spaces.
\begin{enumerate} 
\item An \textbf{$\mathbb{R}$-Weil divisor} (frequently abbreviated to
$\mathbb{R}$-divisor) $D$ on $X$ is an $\mathbb{R}$-linear combination of prime divisors.
\item An \textbf{$\mathbb{R}$-Cartier divisor} $D$ is an $\mathbb{R}$-linear combination of
Cartier divisors.
\item Two $\mathbb{R}$-divisors $D$ and $D'$ are \textbf{$\mathbb{R}$-linearly equivalent
  over $U$}, denoted $D\sim_{\mathbb{R},U} D'$, if their difference is an
$\mathbb{R}$-linear combination of principal divisors and an $\mathbb{R}$-Cartier divisor
pulled back from $U$.
\item Two $\mathbb{R}$-divisors $D$ and $D'$ are \textbf{numerically equivalent over $U$},
denoted $D\equiv_{U} D'$, if their difference is an $\mathbb{R}$-Cartier divisor such that
$(D-D')\cdot C=0$ for any curve $C$ contained in a fibre of $\pi$.
\item An $\mathbb{R}$-Cartier divisor $D$ is \textbf{ample over $U$} (or $\pi$-ample) if
it is $\mathbb{R}$-linearly equivalent to a positive linear combination of ample (in the
usual sense) Cartier divisors over $U$.  
\item An $\mathbb{R}$-Cartier divisor $D$ on $X$ is \textbf{nef over $U$} (or $\pi$-nef) if
$D\cdot C\geq 0$ for any curve $C\subset X$, contracted by $\pi$.
\item An $\mathbb{R}$-divisor $D$ is \textbf{big over $U$} (or $\pi$-big) if
$$
\limsup \frac {h^0(F,\ring F.(\rdown mD.))}{m^{\dim F}}>0,
$$
for the fibre $F$ over any generic point of $U$.  Equivalently $D$ is big over $U$ if
$D\sim_{\mathbb{R},U}A+B$ where $A$ is ample over $U$ and $B\geq 0$ (cf.~\cite[II
3.16]{Nakayama04}).
\item An $\mathbb{R}$-Cartier divisor $D$ is \textbf{semiample over $U$} (or
$\pi$-semiample) if there is a morphism $f\colon\map X.Y.$ over $U$ such that $D$ is
$\mathbb{R}$-linearly equivalent to the pullback of an ample $\mathbb{R}$-divisor  over $U$.
\item An $\mathbb{R}$-divisor $D$ is \textbf{pseudo-effective} if the numerical class of
$D$ over $U$ is a limit of divisors $D_i\geq 0$.
\item An $\mathbb{R}$-divisor $D$ is \textbf{$\pi$-pseudo-effective} if the restriction of 
$D$ to the generic fibre is pseudo-effective.  
\end{enumerate} 
\end{definition}

Note that the group of Weil divisors with rational coefficients $\WDiv_{\mathbb{Q}}(X)$,
or with real coefficients $\WDiv_{\mathbb{R}}(X)$, forms a vector space, with a canonical
basis given by the prime divisors.  Given an $\mathbb{R}$-divisor, $\|D\|$ denotes the sup
norm with respect to this basis.  If $A=\sum a_iC_i$ and $B=\sum b_iC_i$ are two
$\mathbb{R}$-divisors,
$$
A \wedge B= \sum \min(a_i,b_i) C_i.
$$

Given an $\mathbb{R}$-divisor $D$ and a subvariety $Z$ which is not contained in the
singular locus of $X$, $\mult_Z D$ denotes the multiplicity of $D$ at the generic point of
$Z$.  If $Z=E$ is a prime divisor this is the coefficient of $E$ in $D$.

A \textit{log pair} $(X,\Delta)$ (sometimes abbreviated by $K_X+\Delta$) is a normal
variety $X$ and an $\mathbb{R}$-divisor $\Delta\geq 0$ such that $K_X+\Delta$ is
$\mathbb{R}$-Cartier.  We say that a log pair $(X,\Delta)$ is \textit{log smooth}, if $X$
is smooth and the support of $\Delta$ is a divisor with global normal crossings.  A
birational morphism $g\colon\map Y.X.$ is a \textit{log resolution} of the pair $(X,\Delta
)$ if $g$ is projective, $Y$ is smooth, the exceptional locus is a divisor and
$g^{-1}(\Delta)$ union the exceptional set of $g$ is a divisor with global normal
crossings support.  By Hironaka's Theorem we may, and often will, assume that the
exceptional locus supports an ample divisor over $X$.  If we write
$$
K_Y+\Gamma=K_Y+\sum b_i\Gamma _i=g^*(K_X+\Delta),
$$
where $\Gamma _i$ are distinct prime divisors then the log discrepancy
$a(\Gamma_i,X,\Delta)$ of $\Gamma_i$ is $1-b_i$.  The log discrepancy of $(X,\Delta)$ is
then the infimum of the log discrepancy for every $\Gamma_i$ and for every resolution.
The image of any component of $\Gamma$ of coefficient at least one (equivalently log
discrepancy at most zero) is a \textit{non kawamata log terminal centre of the pair $(X,\Delta)$}.
The pair $(X,\Delta )$ is \textit{kawamata log terminal} if for every (equivalently for
one) log resolution $g\colon\map Y.X.$ as above, the coefficients of $\Gamma$ are strictly
less than one, that is $b_i<1$ for all $i$.  Equivalently, the pair $(X,\Delta)$ is
kawamata log terminal if there are no non kawamata log terminal centres.  The \textit{non kawamata log
  terminal locus} of $(X,\Delta)$ is the union of the non kawamata log terminal centres.  We say that
the pair $(X,\Delta)$ is \textit{purely log terminal} if the log discrepancy of any
exceptional divisor is greater than zero.  We say that the pair
$(X,\Delta=\sum\delta_i\Delta_i)$, where $\delta_i\in (0,1]$, is \textit{divisorially log
  terminal} if there is a log resolution such that the log discrepancy of every
exceptional divisor is greater than zero.  By \cite[(2.40)]{KM98} $(X,\Delta)$ is
divisorially log terminal if and only if there is a closed subset $Z\subset X$ such that
\begin{itemize} 
\item $(X\backslash Z,\Delta |_{X\backslash Z})$ is log smooth, and
\item if $f\colon\map Y.X.$ is a projective birational morphism and $E\subset Y$ is an
irreducible divisor with centre contained in $Z$ then $a(E,X,\Delta)>0$.
\end{itemize} 

We will also often write
$$
K_Y+\Gamma=g^*(K_X+\Delta)+E,
$$
where $\Gamma\geq 0$ and $E\geq 0$ have no common components, $g_*\Gamma=\Delta$ and $E$
is $g$-exceptional.  Note that this decomposition is unique.

We say that a birational map $\phi\colon\rmap X.Y.$ is a \textit{birational contraction}
if $\phi$ is proper and $\phi^{-1}$ does not contract any divisors.  If in addition
$\phi^{-1}$ is also a birational contraction, we say that $\phi$ is a \textit{small
  birational map}.  

\subsection{Preliminaries}

\begin{lemma}\label{l_ridiculous} Let $\pi\colon\map X.U.$ be a projective morphism 
of normal quasi-projective varieties.  Let $D$ be an $\mathbb{R}$-Cartier divisor on $X$
and let $D'$ be its restriction to the generic fibre of $\pi$.

If $D'\sim_{\mathbb{R}}B'\geq 0$ for some $\mathbb{R}$-divisor $B'$ on the generic fibre
of $\pi$, then there is a divisor $B$ on $X$ such that $D\sim_{\mathbb{R},U} B\geq 0$.
\end{lemma}
\begin{proof} Taking the closure of the generic points of $B'$, we may assume that there
is an $\mathbb{R}$-divisor $B_1\geq 0$ on $X$ such that the restriction of $B_1$ to the
generic fibre is $B'$.  As
$$
D'-B'\sim_{\mathbb{R}} 0,
$$
it follows that there is an open subset $U_1$ of $U$, such that
$$
(D-B_1)|_{V_1}\sim_{\mathbb{R}} 0,
$$
where $V_1$ is the inverse image of $U_1$.  But then there is a divisor $G$ on $X$ such
that
$$
D-B_1\sim_{\mathbb{R}} G,
$$
where $Z=\pi(\Supp G)$ is a proper closed subset.  As $U$ is quasi-projective, there is an
ample divisor $H\geq 0$ on $U$ which contains $Z$.  Possibly rescaling, we may assume that
$F=\pi^* H\geq -G$.  But then
$$
D \sim_{\mathbb{R}} (B_1+F+G)-F,
$$
so that 
\[
D \sim_{\mathbb{R},U} (B_1+F+G) \geq 0. \qedhere
\]
\end{proof}

\subsection{Nakayama-Zariski decomposition}

We will need some definitions and results from \cite{Nakayama04}.
\begin{definition-lemma}\label{d_sigma} Let $X$ be a smooth projective variety, $B$ be a 
big $\mathbb{R}$-divisor and let $C$ be a prime divisor.  Let
$$
\sigma_C(B)=\inf \{\, \mult_C(B') \,|\, B'\sim_\mathbb{Q}B, B'\geq 0\,\}.
$$
Then $\sigma_C$ is a continuous function on the cone of big divisors. 

Now let $D$ be any pseudo-effective $\mathbb{R}$-divisor and let $A$ be any ample
$\mathbb{Q}$-divisor.  Let 
$$
\sigma_C(D)=\lim_{\epsilon\to 0} \sigma_C(D+\epsilon A).
$$
Then $\sigma_C(D)$ exists and is independent of the choice of $A$.  

There are only finitely many prime divisors $C$ such that $\sigma _{C}(D)>0$ and the
$\mathbb{R}$-divisor $N_\sigma(D)=\sum _C\sigma_C(D)C$ is determined by the numerical
equivalence class of $D$. Moreover $D-N_\sigma(D)$ is pseudo-effective and
$N_\sigma(D-N_\sigma(D))=0$.
\end{definition-lemma}
\begin{proof} See \S III.1 of \cite{Nakayama04}.  \end{proof}

\begin{proposition}\label{p_linear} Let $X$ be a smooth projective variety and let 
$D$ be a pseudo-effective $\mathbb{R}$-divisor.  Let $B$ be any big $\mathbb{R}$-divisor. 

If $D$ is not numerically equivalent to $N_{\sigma}(D)$, then there is a positive integer
$k$ and a positive rational number $\beta$ such that
$$
h^0(X,\ring X.(\rdown mD.+\rdown kB.))> \beta m, \qquad \text{for all} \qquad m\gg 0.
$$
\end{proposition}
\begin{proof} Let $A$ be any integral divisor.  Then we may find a positive integer $k$
such that
$$
h^0(X,\ring X. (\rdown kB.-A))\geq 0.
$$
Thus it suffices to exhibit an ample divisor $A$ and a positive rational number $\beta$
such that
$$
h^0(X,\ring X.(\rdown mD.+A)) > \beta m \qquad \text{for all} \qquad m\gg 0. 
$$
Replacing $D$ by $D-N_{\sigma}(D)$, we may assume that $N_{\sigma}(D)=0$.  Now apply
(V.1.12) of \cite{Nakayama04}.  \end{proof}

\subsection{Adjunction}

We recall some basic facts about adjunction, see \cite[\S16, \S 17]{Kollaretal} for more
details.
\begin{definition-lemma}\label{d_adjunction} Let $(X,\Delta)$ be a log canonical pair, and let $S$ be a normal
component of $\rdown\Delta.$ of coefficient one.  Then there is a divisor $\Theta$ on $S$
such that
$$
(K_X+\Delta)|_S=K_S+\Theta. 
$$
\begin{enumerate}
\item If $(X,\Delta)$ is divisorially log terminal then so is $K_S+\Theta$. 
\item If $(X,\Delta)$ is purely log terminal then $K_S+\Theta$ is kawamata log terminal.  
\item If $(X,\Delta=S)$ is purely log terminal then the coefficients of $\Theta$ have the
form $(r-1)/r$, where $r$ is the index of $S$ at $P$ the generic point of the
corresponding divisor $D$ on $S$ (equivalently $r$ is the index of $K_X+S$ at $P$ or $r$
is the order of the cyclic group $\text{Weil}(\ring {X,P}.)$).  In particular if $B$ is a
Weil divisor on $X$, then the coefficient of $B|_S$ in $D$ is an integer multiple of
$1/r$.
\item If $(X,\Delta)$ is purely log terminal, $f\colon \map Y.X.$ is a projective
birational morphism and $T$ is the strict transform of $S$, then $(f|_T)_*\Psi=\Theta$,
where $K_Y+\Gamma =f^*(K_X+\Delta )$ and $\Psi$ is defined by adjunction,
$$
(K_Y+\Gamma)|_T=K_T+\Psi.
$$ 
\end{enumerate}
\end{definition-lemma}

\subsection{Stable base locus}

We need to extend the definition of the stable base locus to the case of a real divisor.

\begin{definition}\label{d_stable} Let $\pi\colon\map X.U.$ be a projective morphism of 
normal varieties.

Let $D$ be an $\mathbb{R}$-divisor on $X$.  The \textbf{real linear system}
associated to $D$ over $U$ is
$$
|D/U|_{\mathbb{R}}=\{\, C\geq 0\,|\, C\sim_{\mathbb{R},U} D \,\}.
$$

The \textbf{stable base locus} of $D$ over $U$ is the Zariski closed set $\mathbf{B}(D/U)$
given by the intersection of the support of the elements of the real linear system
$|D/U|_{\mathbb{R}}$.  If $|D/U|_{\mathbb{R}}=\emptyset$, then we let $\mathbf{B}(D/U)=X$.
The \textbf{stable fixed divisor} is the divisorial support of the stable base locus.  The
\textbf{augmented base locus} of $D$ over $U$ is the Zariski closed set
$$
\mathbf{B}_{+}(D/U)=\mathbf{B}((D-\epsilon A)/U),
$$
for any ample divisor $A$ over $U$ and any sufficiently small rational number $\epsilon
>0$ (compare \cite[Definition 10.3.2]{Lazarsfeld04b})
\end{definition}

\begin{remark}\label{r_no} The stable base locus, the stable fixed divisor and the
augmented base locus are only defined as closed subsets, they do not have any scheme
structure.
\end{remark}

\begin{lemma}\label{l_same} Let $\pi\colon\map X.U.$ be a projective morphism of 
normal varieties and let $D$ be an integral Weil divisor on $X$.

Then the stable base locus as defined in \eqref{d_stable} coincides with the usual
definition of the stable base locus.
\end{lemma}
\begin{proof} Let 
$$
|D/U|_{\mathbb{Q}}=\{\, C\geq 0\,|\, C\sim_{\mathbb{Q},U} D \,\}.
$$
Let $R$ be the intersection of the elements of $|D|_{\mathbb{R},U}$ and let $Q$ be the
intersection of the elements of $|D|_{\mathbb{Q},U}$.  It suffices to prove that $Q=R$.
As $|D|_{\mathbb{Q}}\subset |D|_{\mathbb{R}}$, it is clear that $R\subset Q$.

Suppose that $x\notin R$.  We want to show that $x\notin Q$.  We may find $D'\in
|D|_{\mathbb{R},U}$ such that $\mult_xD'=0$.  But then 
$$
D'=D+\sum r_i(f_i)+\pi^*E,
$$
where $f_i$ are rational functions on $X$, $E$ is an $\mathbb{R}$-Cartier divisor on $U$
and $r_i$ are real numbers.  Let $V$ be the subspace of $\WDiv _{\mathbb{R}}(X)$ spanned
by the components of $D$, $D'$, $\pi^*E$ and $(f_i)$.  We may write $E=\sum e_jE_j$, where
$E_i$ are Cartier divisors.  Let $W$ be the span of the $(f_i)$ and the components of
$\pi^*E_i$.  Then $W\subset V$ are defined over the rationals.  Set
$$
\mathcal{P}=\{\, D''\in V \,|\, D''\geq 0,\ \mult_xD''=0, \ D''-D\in
W\,\}\subset |D|_\mathbb{R}.
$$
Then $\mathcal{P}$ is a rational polyhedron.  As $D'\in \mathcal{P}$, $\mathcal{P}$ is
non-empty, and so it must contain a rational point $D''$.  We may write
$$
D''=D+\sum s_i(f_i)+\sum f_j\pi^*E_j,
$$
where $s_i$ are real numbers.  Since $D''$ and $D$ have rational coefficients, it follows
that we may find $s_i$ and $f_j$ which are rational.  But then $D''\in |D|_{\mathbb{Q}}$,
and so $x\notin Q$.
\end{proof}

\begin{proposition}\label{p_mobile} Let $\pi\colon\map X.U.$ be projective morphism of 
normal varieties and let $D\geq 0$ be an $\mathbb{R}$-divisor.  Then we may find
$\mathbb{R}$-divisors $M$ and $F$ such that
\begin{enumerate} 
\item $M\geq 0$ and $F\geq 0$, 
\item $D\sim_{\mathbb{R},U}M+F$, 
\item every component of $F$ is a component of $\mathbf{B}(D/U)$, and
\item if $B$ is a component of $M$ then some multiple of $B$ is mobile.  
\end{enumerate} 
\end{proposition}

 We need two basic results.
\begin{lemma}\label{l_equiv} Let $X$ be a normal variety and let $D$ and $D'$ 
be two $\mathbb{R}$-divisors such that $D\sim_{\mathbb{R}} D'$.  

 Then we may find rational functions $\llist f.k.$ and real numbers $\llist r.k.$ which are independent 
over the rationals such that 
$$
D=D'+\sum _i r_i(f_i).  
$$
In particular every component of $(f_i)$ is either a component of $D$ or of $D'$.  
\end{lemma}
\begin{proof} By assumption we may find rational functions $\llist f.k.$ and real numbers $\llist r.k.$
such that 
$$
D=D'+\sum _{i=1}^k r_i(f_i).
$$
Pick $k$ minimal with this property.  Suppose that the real numbers $r_i$ are not
independent over $\mathbb{Q}$.  Then we can find rational numbers $d_i$, not all zero,
such that
$$
\sum_i d_ir_i=0.
$$
Possibly re-ordering we may assume that $d_k\neq 0$.  Multiplying through by an integer we
may assume that $d_i\in\mathbb{Z}$.  Possibly replacing $f_i$ by $f_i^{-1}$, we may assume
that $d_i\geq 0$.  Let $d$ be the least common multiple of the non-zero $d_i$.  If
$d_{i}\neq 0$, we replace $f_i$ by $f_i^{d/d_i}$ (and hence $r_{i}$ by $d_{i}r_{i}/d$) so
that we may assume that either $d_i=0$ or $1$.  For $1\leq i < k$, set
$$
g_i=\begin{cases} f_i/f_k & \text{if $d_i=1$} \\
                  f_i     & \text{if $d_i=0$.}
\end{cases}
$$
Then 
$$
D=D'+\sum _{i=1}^{k-1} r_i(g_i),
$$
which contradicts our choice of $k$.  

 Now suppose that $B$ is a component of $(f_i)$.  Then 
$$
\mult_B(D)=\mult_B(D')+\sum r_j n_j,
$$
where $n_j=\mult_B(f_j)$ is an integer and $n_i\neq 0$.  But then
$\mult_B(D)-\mult_B(D')\neq 0$, so that one of $\mult_B(D)$ and $\mult_B(D')\neq 0$ must
be non-zero.  \end{proof}

\begin{lemma}\label{l_one} Let $\pi\colon\map X.U.$ be a projective morphism of normal 
varieties and let
$$
D' \sim_{\mathbb{R},U} D, \qquad D\geq 0, \qquad D'\geq 0,
$$
be two $\mathbb{R}$-divisors on $X$ with no common components.

Then we may find $D''\in |D/U|_{\mathbb{R}}$ such that a multiple of every component of
$D''$ is mobile.
\end{lemma}
\begin{proof} Pick ample $\mathbb{R}$-divisors on $U$, $H$ and $H'$ such that
$D+\pi^*H\sim_{\mathbb{R}} D'+\pi^*H'$ and $D+\pi^*H$ and $D'+\pi^*H'$ have no common
components.  Replacing $D$ by $D+\pi^*H$ and $D'$ by $D'+\pi^*H'$, we may assume that $D' \sim_{\mathbb{R}} D$.

We may write
$$
D'=D+\sum r_i(f_i)=D+R,
$$
where $r_i\in\mathbb{R}$ and $f_i$ are rational functions on $X$.  By \eqref{l_equiv} we
may assume that every component of $R$ is a component of $D+D'$.

We proceed by induction on the number of components of $D+D'$.  If $\llist q.k.$ are
any positive rational numbers then we may always write
$$
C'=C+Q=C+\sum q_i(f_i),  
$$
where $C\geq 0$ and $C'\geq 0$ have no common components.  But now if we suppose that
$q_i$ is sufficiently close to $r_i$ then $C$ is supported on $D$ and $C'$ is supported on
$D'$.  We have that $mC\sim mC'$ for some integer $m>0$.  By Bertini we may find $C''\sim
_{\mathbb{Q}}C$ such that every component of $C''$ has a multiple which is mobile.  Pick
$\lambda>0$ maximal such that $D_1=D-\lambda C\geq 0$ and $D_1'=D'-\lambda C'\geq 0$.
Note that
$$
D_1\sim_{\mathbb{R}} D_1', \qquad D_1\geq 0, \qquad D_1'\geq 0,
$$ 
are two $\mathbb{R}$-divisors on $X$ with no common components, and that $D_1+D_1'$ has
fewer components than $D+D'$.  By induction we may then find
$$
D_1''\in |D_1|_{\mathbb{R}},
$$
such that a multiple of every component of $D_1''$ is mobile.  But then
$$
D''=\lambda C''+D_1''\in |D|_{\mathbb{R}},
$$
and every component of $D''$ has a multiple which is mobile.   \end{proof}

\begin{proof}[Proof of \eqref{p_mobile}] We may write $D=M+F$, where every component of
$F$ is contained in $\mathbf{B}(D/U)$ and no component of $M$ is contained in
$\mathbf{B}(D/U)$.  A prime divisor is \textbf{bad} if none of its multiples are mobile.

We proceed by induction on the number of bad components of $M$.  We may assume that $M$
has at least one bad component $B$.  As $B$ is a component of $M$, we may find $D_1\in
|D/U|_{\mathbb{R}}$ such that $B$ is not a component of $D_1$.  If $E=D\wedge D_1$ then
$D'=D-E\geq 0$ and $D'_1=D_1-E\geq 0$, $D'$ and $D'_1$ have no common components and $D'
\sim_{\mathbb{R},U} D'_1$.  By \eqref{l_one} there is a divisor $D''\in
|D'/U|_{\mathbb{R}}$ with no bad components.  But then $D''+E\in |D|_{\mathbb{R}}$, $B$ is
not a component of $D''+E$ and the only bad components of $D''+E$ are components of $E$,
which are also components of $D$.  Therefore $D''+E$ has fewer bad components than $D$ and
we are done by induction.
\end{proof}

\subsection{Types of model}

\begin{definition}\label{d_model} Let $\phi\colon\rmap X.Y.$ be a proper birational
contraction of normal quasi-projective varieties and let $D$ be an $\mathbb{R}$-Cartier
divisor on $X$ such that $D'=\phi_*D$ is also $\mathbb{R}$-Cartier.  We say that $\phi$ is
\textbf{$D$-non-positive} (respectively \textbf{$D$-negative}) if for some common
resolution $p\colon\map W.X.$ and $q\colon\map W.Y.$, we may write
$$
p^*D=q^*D'+E,
$$
where $E\geq 0$ is $q$-exceptional (respectively $E\geq 0$ is $q$-exceptional and the
support of $E$ contains the strict transform of the $\phi$-exceptional divisors).
\end{definition}

We will often use the following well-known lemma.
\begin{lemma}[Negativity of contraction]\label{l_negativity} Let $\pi\colon\map Y.X.$ 
be a projective birational morphism of normal quasi-projective varieties.
\begin{enumerate} 
\item If $E>0$ is an exceptional $\mathbb{R}$-Cartier divisor then there is a component
$F$ of $E$ which is covered by curves $\Sigma$ such that $E\cdot\Sigma<0$.
\item If $\pi^*L\equiv M+G+E$, where $L$ is an $\mathbb{R}$-Cartier divisor on $X$, $M$ is a
$\pi$-nef $\mathbb{R}$-Cartier divisor on $Y$, $G\geq 0$, $E$ is $\pi$-exceptional, and
$G$ and $E$ have no common components, then $E\geq 0$.  Further if $F_1$ is an exceptional
divisor such that there is an exceptional divisor $F_2$ with the same centre on $X$ as
$F_1$, with the restriction of $M$ to $F_2$ not numerically $\pi$-trivial, then $F_1$ is a
component of $E$.  
\item If $X$ is $\mathbb{Q}$-factorial then there is a $\pi$-exceptional divisor 
$E\geq 0$ such that $-E$ is ample over $X$.  In particular the exceptional locus 
of $\pi$ is a divisor.  
\end{enumerate} 
\end{lemma}
\begin{proof} Cutting by hyperplanes in $X$, we reduce to the case when $X$ is a surface,
in which case (1) reduces to the Hodge Index Theorem.  (2) follows easily from (1), see
for example (2.19) of \cite{Kollaretal}.  Let $H$ be a general ample $\mathbb{Q}$-divisor 
over $X$.  If $X$ is $\mathbb{Q}$-factorial then 
$$
E=\pi^*\pi_*H-H\geq 0,
$$ 
is $\pi$-exceptional and $-E$ is ample over $X$.  This is (3). 
\end{proof}

\begin{lemma}\label{l_enough} Let $\map X.U.$ and $\map Y.U.$ be two projective 
morphisms of normal quasi-projective varieties.  Let $\phi\colon\rmap X.Y.$ be a
birational contraction over $U$ and let $D$ and $D'$ be $\mathbb{R}$-Cartier divisors such
that $D'=\phi_*D$ is nef over $U$.

Then $\phi$ is $D$-non-positive (respectively $D$-negative) if given a common resolution
$p\colon\map W.X.$ and $q\colon\map W.Y.$, we may write
$$
p^*D=q^*D'+E,
$$
where $p_*E\geq 0$ (respectively $p_*E\geq 0$ and the support of $p_*E$ contains the union of
all $\phi$-exceptional divisors).
 
Further if $D=K_X+\Delta$ and $D'=K_Y+\phi_*\Delta$ then this is equivalent to requiring
\begin{displaymath}
a(F,X,\Delta)\leq a(F,Y,\phi_*\Delta) \qquad \text{(respectively $a(F,X,\Delta )< a(F,Y,\phi_*\Delta)$),}
\end{displaymath} 
for all $\phi$-exceptional divisors $F\subset X$.
\end{lemma}
\begin{proof} Easy consequence of \eqref{l_negativity}.  \end{proof}

\begin{definition}\label{d_ample} Let $\pi\colon\map X.U.$ be a projective morphism of 
normal quasi-projective varieties and let $D$ be an $\mathbb{R}$-Cartier divisor on $X$.

We say that a birational contraction $f\colon\rmap X.Y.$ over $U$ is a \textbf{semiample
  model} of $D$ over $U$, if $f$ is $D$-non-positive, $Y$ is normal and projective over
$U$ and $H=f_*D$ is semiample over $U$.

We say that $g\colon\rmap X.Z.$ is the \textbf{ample model} of $D$ over $U$, if $g$ is a
rational map over $U$, $Z$ is normal and projective over $U$ and there is an ample divisor
$H$ over $U$ on $Z$ such that if $p\colon\map W.X.$ and $q\colon\map W.Z.$ resolve $g$
then $q$ is a contraction morphism and we may write $p^*D \sim_{\mathbb{R},U} q^*H+E$,
where $E\geq 0$ and for every $B\in |p^*D/U|_{\mathbb{R}}$ then $B\geq E$.
\end{definition}

\begin{lemma}\label{l_ample} Let $\pi\colon\map X.U.$ be a projective morphism of normal
quasi-projective varieties and let $D$ be an $\mathbb{R}$-Cartier divisor on $X$.

\begin{enumerate} 
\item If $g_i\colon\rmap X.X_i.$, $i=1$, $2$ are two ample models of $D$ over $U$ then
there is an isomorphism $\chi\colon\map X_1.X_2.$ such that $g_2=\chi\circ g_1$.
\item Suppose that $g\colon\rmap X.Z.$ is the ample model of $D$ over $U$ and let $H$ be
the corresponding ample divisor on $Z$.  If $p\colon\map W.X.$ and $q\colon\map W.Z.$
resolve $g$ then we may write
$$
p^*D \sim_{\mathbb{R},U} q^*H+E,
$$
where $E\geq 0$ and if $F$ is any $p$-exceptional divisor whose centre lies in the
indeterminancy locus of $g$ then $F$ is contained in the support of $E$.
\item If $f\colon\rmap X.Y.$ is a semiample model of $D$ over $U$ then the ample model
$g\colon\rmap X.Z.$ of $D$ over $U$ exists and $g=h\circ f$, where $h\colon\map Y.Z.$ is a
contraction morphism and $f_*D \sim_{\mathbb{R},U} h^*H$.  If $B$ is a prime divisor contained in the stable fixed
divisor of $D$ over $U$ then $B$ is contracted by $f$.
\item If $f\colon\rmap X.Y.$ is a birational map over $U$ then $f$ is the ample model of
$D$ over $U$ if and only if $f$ is a semiample model of $D$ over $U$ and $f_*D$ is ample
over $U$.
\end{enumerate} 
\end{lemma}
\begin{proof} Let $g\colon\map Y.X.$ resolve the indeterminacy of $g_i$ and let
$f_i=g_i\circ g\colon\map Y.X_i.$ be the induced contraction morphisms.  By assumption
$g^*D \sim_{\mathbb{R},U} f_i^*H_i+E_i$, for some divisor $H_i$ on $X_i$ ample over $U$.
Since the stable fixed divisor of $f_{1}^*H_1$ over $U$ is empty, $E_1\geq E_2$.  By
symmetry $E_1=E_2$ and so $f_1^*H_1 \sim_{\mathbb{R},U} f_2^*H_2$.  But then $f_1$ and
$f_2$ contract the same curves.  This is (1).

Suppose that $g\colon\rmap X.Z.$ is the ample model of $D$.  By assumption this means that
we may write
$$
p^*D \sim_{\mathbb{R},U} q^*H+E,
$$
where $E\geq 0$.  We may write $E=E_1+E_2$, where every component of $E_2$ is exceptional
for $p$ but no component of $E_1$ is $p$-exceptional.  Let $V=p(F)$.  Possibly blowing up
more we may assume that $p^{-1}(V)$ is a divisor.  Since $V$ is contained in the
indeterminancy locus of $g$, there is an exceptional divisor $F'$ with centre $V$ such
that $\dim q(F')>0$.  But then $q^*H$ is not numerically trivial on $F'$ and we may apply
\eqref{l_negativity}.  This is (2).

Now suppose that $f\colon\rmap X.Y.$ is a semiample model of $D$ over $U$.  As $f_*D$ is
semiample over $U$, there is a contraction morphism $h\colon\map Y.Z.$ over $U$ and an
ample divisor $H$ over $U$ on $Z$ such that $f_*D \sim_{\mathbb{R},U} h^*H$.  If
$p\colon\map W.X.$ and $q\colon\map W.Y.$ resolve the indeterminacy of $f$ then $p^*D
\sim_{\mathbb{R},U} r^*H+E$, where $E\geq 0$ is $q$-exceptional and $r=h\circ q\colon\map
W.Z.$.  If $B\in |p^*D/U|_{\mathbb{R}}$ then $B\geq E$.  But then $g=h\circ f\colon\rmap
X.Z.$ is the ample model of $D$ over $U$.  This is (3).

Now suppose that $f\colon\rmap X.Y.$ is birational over $U$.  If $f$ is a semiample model
of $D$ over $U$ then (3) implies that the ample model of $D$ over $U$ exists $g\colon\rmap
X.Z.$ and there is a contraction morphism $h\colon\map Y.Z.$, such that $f_*D
\sim_{\mathbb{R},U} h^*H$ where $H$ on $Z$ is ample over $U$.  If $f_*D$ is ample over $U$
then $h$ must be the identity.

Conversely suppose that $f$ is the ample model.  Suppose that $p\colon\map W.X.$ and
$q\colon\map W.Y.$ are projective birational morphisms which resolve $g$.  By assumption
we may write $p^*D \sim_{\mathbb{R},U} q^*H+E$, where $H$ is ample over $U$.  We may
assume that there is a $q$-exceptional $\mathbb Q$-divisor $F\geq 0$ such that $q^*H-F$ is
ample over $U$.  Then there is a constant $\delta>0$ such that $q^*H-F+\delta E$ is ample over $U$.
Suppose $B$ is a component of $E$.  As $B$ does not belong to the stable base locus of
$q^*H-F+\delta E$ over $U$, $B$ must be a component of $F$.  It follows that $E$ is
$q$-exceptional. If $G$ is a divisor on $W$ and $C$ is a curve in $G$ contracted by $p$ but not contracted by $q$, then $0=C\cdot p^*D=C\cdot q^{*}H+
C\cdot E$ and so $C$ and hence $G$ are contained in the support of $E$.
Therefore $f$ is a birational contraction and $f$ is a semiample model.  Further $f_*D=H$ is ample over $U$.
This is (4).
\end{proof}

\begin{definition}\label{d_log} Let $\pi\colon\map X.U.$ be a projective morphism of 
normal quasi-projective varieties.  Suppose that $K_X+\Delta$ is log canonical and let
$\phi\colon\rmap X.Y.$ be a birational contraction of normal quasi-projective varieties
over $U$, where $Y$ is projective over $U$.  Set $\Gamma=\phi_*\Delta$.
\begin{itemize} 
\item $Y$ is a \textbf{weak log canonical model} for $K_X+\Delta$ over $U$ if $\phi$ is
$(K_X+\Delta)$-non-positive and $K_Y+\Gamma$ is nef over $U$.
\item $Y$ is the \textbf{log canonical model} for $K_X+\Delta$ over $U$ if $\phi$ is
the ample model of $K_X+\Delta$ over $U$.  
\item $Y$ is a \textbf{log terminal model} for $K_X+\Delta$ over $U$ if $\phi$ is
$(K_X+\Delta)$-negative, $K_Y+\Gamma$ is divisorially log terminal and nef over $U$, and
$Y$ is $\mathbb{Q}$-factorial.
\end{itemize} 
\end{definition}

\begin{remark} Note that there is no consensus on the definitions given in \eqref{d_log}.
\end{remark}

\begin{lemma}\label{l_equiv-models} Let $\pi\colon\map X.U.$ be a projective morphism of
normal quasi-projective varieties.  Let $\phi\colon\rmap X.Y.$ be a birational contraction
over $U$.  Let $(X,\Delta)$ and $(X,\Delta')$ be two log pairs and set
$\Gamma=\phi_*\Delta$, $\Gamma'=\phi_*\Delta'$.  Let $\mu>0$ be a positive real number.
\begin{itemize} 
\item If both $K_X+\Delta$ and $K_X+\Delta'$ are log canonical and $K_X+\Delta'
\sim_{\mathbb{R} ,U} \mu(K_X+\Delta)$ then $\phi$ is a weak log canonical
model for $K_X+\Delta$ over $U$ if and only if $\phi$ is a weak log canonical model for
$K_X+\Delta'$ over $U$. 
\item If both $K_X+\Delta$ and $K_X+\Delta'$ are kawamata log terminal and $K_X+\Delta
\equiv_U \mu (K_X+\Delta')$ then $\phi$ is a log terminal model for $K_X+\Delta$ over $U$
if and only if $\phi$ is a log terminal model for $K_X+\Delta'$ over $U$.
\end{itemize} 
\end{lemma}
\begin{proof} Note first that either $K_Y+\Gamma' \sim_{\mathbb{R},U} \mu(K_Y+\Gamma)$ or
$Y$ is $\mathbb{Q}$-factorial.  In particular $K_Y+\Gamma$ is $\mathbb{R}$-Cartier if and
only if $K_Y+\Gamma'$ is $\mathbb{R}$-Cartier.  If $p\colon\map W.X.$ and $q\colon\map
W.Y.$ resolve the indeterminacy of $\phi=q\circ p^{-1}$, then we may write
$$
p^*(K_X+\Delta)=q^*(K_Y+\Gamma)+E \quad \text{and} \quad p^*(K_X+\Delta')=q^*(K_Y+\Gamma')+E'.
$$  
Since $\mu E-E'\equiv_Y 0$ is $q$-exceptional, $\mu E=E'$ by \eqref{l_negativity}.
Therefore, $\phi$ is $(K_X+\Delta)$-non-positive (respectively $(K_X+\Delta)$-negative) if
and only if $\phi$ is $(K_X+\Delta')$-non-positive (respectively
$(K_X+\Delta')$-negative).

Since $K_X+\Delta'\equiv_U \mu(K_X+\Delta)$ and $\mu E-E'\equiv_Y 0$, it follows that
$K_Y+\Gamma'\equiv_U \mu(K_Y+\Gamma)$, so that $K_Y+\Gamma$ is nef over $U$ if and only if
$K_Y+\Gamma'$ is nef over $U$.  \end{proof}

\begin{lemma}\label{l_modelXY} Let $\pi\colon\map X.U.$ be a projective morphism of 
normal quasi-projective varieties.  Suppose that $K_X+\Delta$ is divisorially log terminal
and let $\phi\colon\rmap X.Y.$ be a birational contraction over $U$, $K_Y+\phi _*\Delta$
is divisorially log terminal and $a(F,X,\Delta)<a(F,Y,\phi _*\Delta)$ for all
$\phi$-exceptional divisors $F\subset X$.

If $\varphi\colon\rmap Y.Z.$ is a log terminal model of $(Y,\phi_*\Delta)$ over $U$, then
$\eta=\varphi\circ \phi\colon\rmap X.Z.$ is a log terminal model of $K_X+\Delta$ over $U$.
\end{lemma}
\begin{proof} Clearly $\eta$ is a birational contraction, $Z$ is $\mathbb{Q}$-factorial
and $K_Z+\eta_*\Delta$ is divisorially log terminal and nef over $U$.
  
Let $p\colon\map W.X.$, $q\colon\map W.Y.$ and $r\colon\map W.Z.$ be a common resolution.
As $\varphi$ is a log terminal model of $(Y,\phi _*\Delta)$ we have that $q^*(K_Y+\phi
_*\Delta)-r^*(K_Z+\eta _*\Delta)=E\geq 0$ and the support of $E$ contains the exceptional
divisors of $\varphi$.  By assumption $K_X+\Delta-p_*q^*(K_Y+\phi _*\Delta)$ is an
effective divisor whose support is the set of all $\phi$ exceptional divisors.  But then
$$
(K_X+\Delta)-p_*r^*(K_Z+\eta_*\Delta)=K_X+\Delta-p_*q^*(K_Y+\phi _*\Delta)+p_*E
$$ 
contains all the $\eta$-exceptional divisors and \eqref{l_enough} implies that $\eta$ is a
log terminal model of $K_X+\Delta$ over $U$.  
\end{proof}

\begin{lemma}\label{l_high} Let $\pi\colon\map X.U.$ be a projective morphism of 
normal quasi-projective varieties.  Let $(X,\Delta)$ be a kawamata log terminal pair, 
where $\Delta$ is big over $U$.  Let
$f\colon\map Z.X.$ be any log resolution of $(X,\Delta)$ and suppose that we write
$$
K_Z+\Phi_0=f^*(K_X+\Delta)+E,
$$
where $\Phi_0\geq 0$ and $E\geq 0$ have no common components, $f_*\Phi_0=\Delta$ and $E$
is exceptional.  Let $F\geq 0$ be any divisor whose support is equal to the exceptional
locus of $f$.

If $\eta>0$ is sufficiently small and $\Phi=\Phi_0+\eta F$ then $K_Z+\Phi$ is kawamata log
terminal and $\Phi$ is big over $U$.  Moreover if $\phi\colon\rmap Z.W.$ is a log terminal
model of $K_Z+\Phi$ over $U$ then the induced birational map $\psi\colon\rmap X.W.$ is in fact a
log terminal model of $K_X+\Delta$ over $U$.
\end{lemma}
\begin{proof} Set $\Psi=\phi_*\Phi$.  By \eqref{l_modelXY}, possibly blowing up more, we
may assume that $\phi$ is a morphism.  By assumption if we write
$$
K_Z+\Phi=\phi^*(K_W+\Psi)+G,
$$
then $G\geq 0$ and the support of $G$ is the union of all the $\phi$-exceptional divisors.
Thus
$$
f^*(K_X+\Delta)+E+\eta F=\phi^*(K_W+\Psi)+G.
$$
By negativity of contraction, \eqref{l_negativity}, applied to $f$, $G-E-\eta F\geq 0$.
In particular $\phi$ must contract every $f$-exceptional divisor and so $\psi$ is a
birational contraction.  But then $\psi$ is a log terminal model over $U$ by
\eqref{l_enough}.  \end{proof}

\begin{lemma}\label{l_ample-trick} Let $\pi\colon\map X.U.$ be a projective 
morphism of normal quasi-projective varieties, where $X$ is $\mathbb{Q}$-factorial.  Let
$\phi\colon\rmap X.Y.$ be a weak log canonical model over $U$ of a divisorially log
terminal pair $(X,\Delta=S+B)$ where $S=\rdown\Delta.$ is a sum of prime divisors and
$B\geq 0$.  Suppose that the components of $B$ span $\WDiv_{\mathbb{R}}(X)$ modulo
numerical equivalence over $U$.

If $V$ is any finite dimensional affine subspace of $\WDiv_{\mathbb{R}}(X)$ which contains the subspace generated by 
the components of $B$ then
$$
\mathcal{W}_{\phi,S,\pi}(V)=\bar{\mathcal{A}}_{\phi,S,\pi}(V).
$$
\end{lemma}
\begin{proof} By (4) of \eqref{l_ample} 
$\mathcal{W}_{\phi,S,\pi}(V)\supset {\mathcal{A}}_{\phi,S,\pi}(V).$.  Since $\mathcal{W}_{\phi,S,\pi}(V)$ is closed, it follows that
$$
\mathcal{W}_{\phi,S,\pi}(V)\supset \bar{\mathcal{A}}_{\phi,S,\pi}(V).
$$
To prove the reverse inclusion, it suffices to prove that a dense subset of
$\mathcal{W}_{\phi,S,\pi}(V)$ is contained in $\mathcal{A}_{\phi,S,\pi}(V)$.

Let $H$ be a general ample $\mathbb{Q}$-divisor over $U$ on $Y$.  Let $p\colon\map W.X.$
and $q\colon \map W.Y.$ resolve the indeterminacy locus of $\phi$ and let $H'=p_*q^*H$.
It follows that $\phi$ is $H'$-negative.  Pick $B'$ numerically
equivalent to $H'$ over $U$ such that the support of $B'$ is contained in the support of
$B$.

If $\epsilon>0$ is any positive real number then $\phi$ is $(K_X+\Delta+\epsilon
B')$-negative and $\phi_*(K_X+\Delta+\epsilon B')$ is ample over $U$.  If we pick
$\epsilon>0$ such that $\Delta+\epsilon B'\in \mathcal{L}_S(V)$ then $\phi$ is the ample
model of $K_X+\Delta+\epsilon B'$ over $U$ and $\Delta+\epsilon B'\in
\mathcal{A}_{\phi,S,\pi}(V)$.  \end{proof}

\subsection{Convex geometry and Diophantine approximation}

\begin{definition}\label{d_polytope} Let $V$ be a finite dimensional real affine space.  
If $\mathcal{C}$ is a convex subset of $V$ and $F$ is a convex subset of $\mathcal{C}$,
then we say that $F$ is a \textbf{face} of $\mathcal{C}$ if whenever $\sum r_i v_i \in F$,
where $\llist r.k.$ are real numbers such that $\sum r_i=1$, $r_i\geq 0$ and $\llist v.k.$
belong to $\mathcal{C}$, then $v_i\in F$ for some $i$.  We say that $v\in \mathcal{C}$ is
an \textbf{extreme point} if $F=\{v\}$ is a face of $\mathcal{C}$.

A \textbf{polyhedron} $\mathcal{P}$ in $V$ is the intersection of finitely many half
spaces.  The \textbf{interior} $\mathcal{P}^{\circ}$ of $\mathcal{P}$ is the complement of
the proper faces.  A \textbf{polytope} $\mathcal{P}$ in $V$ is a compact polyhedron.

We say that a real vector space $V_0$ is \textbf{defined over the rationals}, if $V_0=\ten
V'.{\mathbb{Q}}.\mathbb{R}.$, where $V'$ is a rational vector space.  We say that an
affine subspace $V$ of a real vector space $V_0$, which is defined over the rationals, is
\textbf{defined over the rationals}, if $V$ is spanned by a set of rational vectors of
$V_0$.  We say that a polyhedron $\mathcal{P}$ is \textbf{rational} if it is defined 
by rational half spaces.  
\end{definition}

Note that a polytope is the convex hull of any finite set of points and the polytope is
rational if those points can be chosen to be rational.

\begin{lemma}\label{l_polytopal} Let $X$ be a normal quasi-projective variety, and let 
$V$ be a finite dimensional affine subspace of $\WDiv _{\mathbb{R}}(X)$, 
which is defined over the rationals.

Then $\mathcal{L}(V)$ (cf. \eqref{d_cones} for the definition) is a rational polytope.
\end{lemma}
\begin{proof} Note that the set of divisors $\Delta$ such that $K_X+\Delta$ is
$\mathbb{R}$-Cartier forms an affine subspace $W$ of $V$, which is defined over the
rationals, so that, replacing $V$ by $W$, we may assume that $K_X+\Delta$ is
$\mathbb{R}$-Cartier for every $\Delta\in V$.

Let $\pi\colon\map Y.X.$ be a resolution of $X$, which is a log resolution of the support
of any element of $V$.  Given any divisor $\Delta\in V$, if we write
$$
K_Y+\Gamma=\pi^*(K_X+\Delta),
$$
then the coefficients of $\Gamma$ are rational affine linear functions of the coefficients
of $\Delta$.  On the other hand the condition that $K_X+\Delta$ is log canonical is
equivalent to the condition that the coefficient of every component of $\Gamma$ is at most
one and the coefficient of every component of $\Delta$ is at least zero.  \end{proof}

\begin{lemma}\label{l_linear-ample} Let $\pi\colon\map X.U.$ be a projective morphism 
of normal quasi-projective varieties.  Let $V$ be a finite dimensional affine subspace of
$\WDiv_{\mathbb{R}}(X)$ and let $A\geq 0$ be a big $\mathbb{R}$-divisor over $U$.  Let
$\mathcal{C}\subset \mathcal{L}_A(V)$ be a polytope.

If $\mathbf{B}_+(A/U)$ does not contain any non kawamata log terminal centres of $(X,\Delta)$, for
every $\Delta\in \mathcal{C}$, then we may find a general ample $\mathbb{Q}$-divisor $A'$
over $U$, a finite dimensional affine subspace $V'$ of $\WDiv_{\mathbb{R}}(X)$ and a
translation
$$
L\colon\map \WDiv_{\mathbb{R}}(X).\WDiv _{\mathbb{R}}(X).,
$$
by a $\mathbb{R}$-divisor $T$ $\mathbb{R}$-linearly equivalent to zero over $U$ such that
$L(\mathcal{C})\subset \mathcal{L}_{A'}(V')$ and $(X,\Delta-A)$ and $(X,L(\Delta))$ have
the same non kawamata log terminal centres.  Further, if $A$ is a $\mathbb{Q}$-divisor then we may
choose $T$ $\mathbb{Q}$-linearly equivalent to zero over $U$.
\end{lemma}
\begin{proof} Let $\llist \Delta.k.$ be the vertices of the polytope $\mathcal{C}$.  Let
$\mathfrak{Z}$ be the set of non kawamata log terminal centres of $(X,\Delta_i)$ for $1\leq i \leq k$.  Note that
if $\Delta\in \mathcal{C}$ then any non kawamata log terminal centre of $(X,\Delta)$ is an element of
$\mathfrak{Z}$.

By assumption, we may write $A \sim_{\mathbb{R},U}C+D$, where $C$ is a general ample
$\mathbb{Q}$-divisor over $U$ and $D\geq 0$ does not contain any element of $\mathfrak{Z}$.  Further
\eqref{l_same} implies that if $A$ is a $\mathbb{Q}$-divisor then we may assume that $A
\sim_{\mathbb{Q},U} C+D$.

Given any rational number $\delta>0$, let
$$
L\colon\map \WDiv _{\mathbb{R}}(X).\WDiv _{\mathbb{R}}(X). \quad \text{given by} \quad L(\Delta)=\Delta+\delta (C+D-A),
$$
be the translation by the divisor $T=\delta(C+D-A)\sim_{\mathbb{R},U} 0$.  Note that
$T\sim_{\mathbb{Q},U} 0$ if $A$ is a $\mathbb{Q}$-divisor.  As $C+D$ does not contain any
element of $\mathfrak{Z}$, if $\delta$ is sufficiently small then
$$
K_X+L(\Delta_i)=K_X+\Delta_i+\delta(C+D-A)=K_X+\delta C+(\Delta_i-\delta A+\delta D),
$$
is log canonical for every $1\leq i\leq k$ and has the same non kawamata log terminal centres as
$(X,\Delta_i-A)$.  But then $L(\mathcal{C})\subset \mathcal{L}_{A'}(V')$, where $A'=\delta
C$ and $V'=V_{(1-\delta)A+\delta D}$ and $(X,\Delta-A)$ and $(X,L(\Delta))$ have the same
non kawamata log terminal centres.  \end{proof}

\begin{lemma}\label{l_linear-dlt} Let $\pi\colon\map X.U.$ be a projective morphism 
of normal quasi-projective varieties.  Let $V$ be a finite dimensional affine subspace of
$\WDiv_{\mathbb{R}}(X)$, which is defined over the rationals, and let $A$ be a general
ample $\mathbb{Q}$-divisor over $U$.  Let $S$ be a sum of prime divisors.  Suppose that
there is a divisorially log terminal pair $(X,\Delta_0)$, where $S=\rdown\Delta_0.$, and
let $G\geq 0$ be any divisor whose support does not contain any non kawamata log terminal
centres of $(X,\Delta_0)$.  

Then we may find a general ample $\mathbb{Q}$-divisor $A'$ over $U$, an affine subspace
$V'$ of $\WDiv_{\mathbb{R}}(X)$, which is defined over the rationals, a divisor
$\Delta_0'\in \mathcal{L}_{S+A'}(V')$ and a rational affine linear isomorphism
$$
L\colon\map V_{S+A}.V'_{S+A'}.,
$$
such that 
\begin{itemize} 
\item $L$ preserves $\mathbb{Q}$-linear equivalence over $U$, 
\item $L(\mathcal{L}_{S+A}(V))$ is contained in the interior of $\mathcal{L}_{S+A'}(V')$, 
\item $K_X+\Delta_0'$ is divisorially log terminal and $\rdown\Delta'_0.=S$, and 
\item for any $\Delta\in L(\mathcal{L}_{S+A}(V))$, the support of $\Delta$ contains the
support of $G$.
\end{itemize} 
\end{lemma}
\begin{proof} Let $W$ be the vector space spanned by the components of $\Delta_0$.  Then
$\Delta_0\in \mathcal{L}_S(W)$ and \eqref{l_polytopal} implies that $\mathcal{L}_S(W)$ is
a non-empty rational polytope.  But then $\mathcal{L}_S(W)$ contains a rational point and
so, possibly replacing $\Delta_0$, we may assume that $K_X+\Delta_0$ is
$\mathbb{Q}$-Cartier.

We first prove the result in the case that $K_X+\Delta$ is $\mathbb{R}$-Cartier for every
$\Delta\in V_{S+A}$.  By compactness, we may pick $\mathbb{Q}$-divisors $\llist \Delta.l.\in V_{S+A}$
such that $\mathcal{L}_{S+A}(V)$ is contained in the simplex spanned by $\llist \Delta.l.$
(we do not assume that $\Delta_i\geq 0$).  Pick a rational number $\epsilon\in (0,1/4]$
such that
$$
\epsilon(\Delta_i-\Delta_0)+(1-2\epsilon)A,
$$
is an ample $\mathbb{Q}$-divisor over $U$, for $1\leq i\leq l$.  Pick
$$
A_i \sim_{\mathbb{Q},U} \epsilon(\Delta_i-\Delta_0)+(1-2\epsilon)A,
$$
general ample $\mathbb{Q}$-divisors over $U$.  Pick $A' \sim_{\mathbb{Q},U} \epsilon A$ a
general ample $\mathbb{Q}$-divisor over $U$.  If we define $L\colon\map V_{S+A}.\WDiv
_{\mathbb R}(X).$ by
$$
L(\Delta_i)=(1-\epsilon)\Delta_i+A_i+\epsilon \Delta_0+A'-(1-\epsilon) A \sim_{\mathbb{Q},U} \Delta_i,
$$
and extend to the whole of $V_{S+A}$ by linearity, then $L$ is an injective rational
linear map which preserves $\mathbb{Q}$-linear equivalence over $U$.
We let $V'$ be the rational affine subspace of $\WDiv _{\mathbb R}(X)$
defined by $V'_{S+A'}=L(V_{S+A})$.
  Note that
$K_X+\Delta_0'$ is divisorially log terminal, where $\Delta'_0=A'+\Delta_0\in
\mathcal{L}_{S+A'}(V')$.  Note also that $L$ is the composition of
$$
L_1(\Delta_i)=\Delta_i+A_i/(1-\epsilon)+A'-A \quad \text{and} \quad L_2(\Delta)=(1-\epsilon)\Delta+\epsilon\Delta'_0.
$$
If $\Delta\in \mathcal{L}_{S+A}(V)$ then $K_X+\Delta+A'-A$ is log canonical, and as $A_i$
is a general ample $\mathbb{Q}$-divisor over $U$ it follows that $K_X+\Delta+2A_i+A'-A$ is
log canonical as well.  As $1/(1-\epsilon)<2$, it follows that if $\Delta\in
\mathcal{L}_{S+A}(V)$ then $K_X+L_1(\Delta)$ is log canonical.  Therefore, if $\Delta\in
\mathcal{L}_{S+A}(V)$ then $K_X+L(\Delta)$ is divisorially log terminal and $\rdown
L(\Delta).=S$.

Pick a divisor $G'$ such that $S+A'+G'$ belongs to the interior of $L_{S+A'}(V')$.  As
$G+G'$ contains no log canonical centres of $(X,\Delta_0)$ and $X$ is smooth at the
generic point of every log canonical centre of $(X,\Delta_0)$, we may pick a
$\mathbb{Q}$-Cartier divisor $H\geq G+G'$ which contains no log canonical centres of
$(X,\Delta_0)$.  Pick a rational number $\eta>0$ such that $A'-\eta H$ is ample over $U$.
Pick $A''\sim_{\mathbb{Q},U} A'-\eta H$ a general ample $\mathbb{Q}$-divisor over $U$.
Let $\delta>0$ by any rational number and let
$$
T\colon\map \WDiv _{\mathbb{R}}(X).\WDiv _{\mathbb{R}}(X).,
$$
be translation by $\delta (\eta H+A''-A') \sim_{\mathbb{Q},U} 0$.  If $V''$ is the span of
$V'$, $A'$ and $H$ and $\delta>0$ is sufficiently small then $T(L(\mathcal{L}_{S+A}(V)))$
is contained in the interior of $\mathcal{L}_{\delta A''+S}(V'')$, $K_X+T(\Delta'_0)$ is
divisorially log terminal and the support of $T(\Delta'_0)$ contains the support of $G$.
If we replace $L$ by $T\circ L$, $V'_{S+A'}$ by $T(L(V_{S+A}))$, $A'$ by $\delta A''$ and $\Delta'_0$
by $T(\Delta'_0)$ then this finishes the case when $K_X+\Delta$ is $\mathbb{R}$-Cartier
for every $\Delta\in V_{S+A}$.

We now turn to the general case.  If 
$$
W_0=\{\, B\in V \,|\, \text{$K_X+S+A+B$ is $\mathbb{R}$-Cartier} \,\},
$$
then $W_0\subset V$ is an affine subspace of $V$, which is defined over the rationals.
Note that $\mathcal{L}_{S+A}(V)=\mathcal{L}_{S+A}(W_0)$.  By what we have already proved,
there is a rational affine linear isomorphism $L_0\colon\map W_0.W'_0.$, which preserves
$\mathbb{Q}$-linear equivalence over $U$, a general ample $\mathbb{Q}$-divisor $A'$ over
$U$, such that $L_0(\mathcal{L}_{S+A}(W_0))$ is contained in the interior of
$\mathcal{L}_{S+A'}(W'_0)$, and there is a divisor $\Delta_0'\in \mathcal{L}_{S+A'}(W'_0)$
such that $K_X+\Delta'_0$ is divisorially log terminal and the support of $\Delta'_0$
contains the support of $G$.

Let $W_1$ be any subspace of $\WDiv_{\mathbb{R}}(X)$, which is defined over the rationals,
such that $V=W_0+W_1$ and $W_0\cap W_1\subset \{0\}$.  Let $V'=W'_0+W_1$.  Since $L_0$
preserves $\mathbb{Q}$-linear equivalence over $U$, $W'_0\cap W_1=W_0\cap W_1$ and
$\mathcal{L}_{S+A'}(V')=\mathcal{L}_{S+A'}(W'_0)$.  If we define $L\colon\map
V_A.V'_{A'}.$, by sending $A+B_0+B_1$ to $L_0(A+B_0)+B_1$, where $B_i\in W_i$, then $L$ is
a rational affine linear isomorphism, which preserves $\mathbb{Q}$-linear equivalence over
$U$, $L(\mathcal{L}_{S+A}(V))$ is contained in the interior of $\mathcal{L}_{S+A'}(V')$
and $\Delta'_0\in \mathcal{L}_{S+A'}(V')$.
\end{proof}

\begin{lemma}\label{l_big-to-ample} Let $\pi\colon\map X.U.$ be a projective morphism 
of normal quasi-projective varieties.  Let $(X,\Delta=A+B)$ be a log canonical pair, where
$A\geq 0$ and $B\geq 0$.

If $\mathbf{B}_+(A/U)$ does not contain any non kawamata log terminal centres of $(X,\Delta)$ and
there is a kawamata log terminal pair $(X,\Delta_0)$ then we may find a kawamata log
terminal pair $(X,\Delta'=A'+B')$, where $A'\geq 0$ is a general ample
$\mathbb{Q}$-divisor over $U$, $B'\geq 0$ and $K_X+\Delta' \sim_{\mathbb{R},U}
K_X+\Delta$.  If in addition $A$ is a $\mathbb{Q}$-divisor then $K_X+\Delta'
\sim_{\mathbb{Q},U} K_X+\Delta$.
\end{lemma}
\begin{proof} By \eqref{l_linear-ample} we may assume that $A$ is a general ample
$\mathbb{Q}$-divisor over $U$.  If $V$ is the vector space spanned by the components of
$\Delta$ then $\Delta\in \mathcal{L}_A(V)$ and the result follows by \eqref{l_linear-dlt}.
\end{proof}

\begin{lemma}\label{l_dense} Let $V$ be a finite dimensional real vector space, which is 
defined over the rationals.  Let $\Lambda\subset V$ be a lattice spanned by rational
vectors.  Suppose that $v\in V$ is a vector which is not contained in any proper affine
subspace $W\subset V$ which is defined over the rationals.
 
 Then the set
$$
X=\{\, mv+\lambda \,|\, m\in\mathbb{N}, \lambda\in\Lambda \,\},
$$
is dense in $V$. 
\end{lemma}
\begin{proof} Let
$$
q\colon\map V.V/\Lambda.,
$$
be the quotient map and let $G$ be the closure of the image of $X$.  As $G$ 
is infinite and $V/\Lambda$ is compact, $G$ has an accumulation point.  It 
then follows that zero is also an accumulation point and that $G$ is a closed 
subgroup.  

The connected component $G_0$ of $G$ containing the identity is a Lie subgroup of
$V/\Lambda$ and so by Theorem 15.2 of \cite{Bump04}, $G_0$ is a torus.  Thus
$G_0=W/\Lambda_0$, where
$$
W=H_1(G_0,\mathbb{R})=\ten \Lambda_0.\mathbb{Z}.\mathbb{R}.=\ten H_1(G_0,\mathbb{Z}).\mathbb{Z}.\mathbb{R}.\subset \ten H_1(G,\mathbb{Z}).\mathbb{Z}.\mathbb{R}.=H_1(G,\mathbb{R}),
$$
is a subspace of $V$ which is defined over the rationals.  On the other hand,
$G/G_0$ is finite as it is discrete and compact.  Thus a translate of $v$ by a rational
vector is contained in $W$ and so $W=V$.  \end{proof}

\begin{lemma}\label{l_diophantine} Let $\mathcal{C}$ be a rational polytope contained in a 
real vector space $V$ of dimension $n$, which is defined over the rationals.  Fix a
positive integer $k$ and a positive real number $\alpha$.

If $v\in\mathcal{C}$ then we may find vectors $\llist v.p.\in\mathcal{C}$ and positive
integers $\llist m.p.$, which are divisible by $k$, such that $v$ is a convex linear
combination of the vectors $\llist v.p.$ and
$$
\|v_i-v\|< \frac{\alpha}{m_i} \qquad \text{where} \qquad \text{$\frac {m_i v_i} k$ is integral.}
$$
\end{lemma}
\begin{proof} Rescaling by $k$, we may assume that $k=1$.  We may assume that $v$ is not
contained in any proper affine linear subspace which is defined over the rationals.  In
particular $v$ is contained in the interior of $\mathcal{C}$ since the faces of
$\mathcal{C}$ are rational.

After translating by a rational vector, we may assume that $0\in \mathcal{C}$.  After
fixing a suitable basis for $V$ and possibly shrinking $\mathcal{C}$, we may assume that
$\mathcal{C}=[0,1]^n\subseteq \mathbb{R}^n$ and $v=(\llist x.n.)\in (0,1)^n$.  By
\eqref{l_dense}, for each subset $I\subseteq \{1,2,\dots,n\}$, we may find 
$$
v_I=(\llist s.n.)\in (0,1)^n\cap \mathbb{Q}^n,
$$ 
and an integer $m_I$ such that $m_I v_I$ is integral, such that 
$$
\qquad \|v-v_I\|< \frac \alpha {m_I} \qquad\text{and } s_j <x_j \quad\text{if and only if }j\in I. 
$$
In particular $v$ is contained inside the rational polytope $\mathcal{B}\subseteq
\mathcal{C}$ generated by the $v_I$.  Thus $v$ is a convex linear combination of a subset
$\llist v.p.$ of the extreme points of $\mathcal{B}$.  \end{proof}

\subsection{Rational curves of low degree}

We will need the following generalisation of a result of Kawamata, see Theorem 1 of
\cite{Kawamata91}, which is proved by Shokurov in the appendix to \cite{Nikulin94}.

\begin{theorem}\label{t_low} Let $\pi\colon\map X.U.$ be a projective morphism 
of normal quasi-projective varieties.  Suppose that $(X,\Delta)$ is a log canonical pair of
dimension $n$, where $K_X+\Delta$ is $\mathbb{R}$-Cartier.  Suppose that there is a
divisor $\Delta_0$ such that $K_X+\Delta_0$ is kawamata log terminal.

If $R$ is an extremal ray of $\ccone X/U.$ that is $(K_X+\Delta)$-negative, then there is
a rational curve $\Sigma$ spanning $R$, such that
$$
0<-(K_X+\Delta)\cdot \Sigma \leq 2n.
$$
\end{theorem}
\begin{proof} Passing to an open subset of $U$, we may assume that $U$ is affine.  Let $V$
be the vector pace spanned by the components of $\Delta+\Delta_0$.  By \eqref{l_polytopal}
the space $\mathcal{L}(V)$ of log canonical divisors is a rational polytope.  Since
$\Delta_0\in \mathcal{L}(V)$, we may find $\mathbb{Q}$-divisors $\Delta_i\in V$ with limit
$\Delta$, such that $K_X+\Delta_i$ is kawamata log terminal.  In particular we may assume
that $(K_X+\Delta_0)\cdot R<0$.  Replacing $\pi$ by the contraction defined by the
extremal ray $R$, we may assume that $-(K_X+\Delta)$ is $\pi$-ample.

Theorem 1 of \cite{Kawamata91} implies that we can find a rational curve $\Sigma_i$
contracted by $\pi$ such that
$$
-(K_X+\Delta_i)\cdot \Sigma_i\leq 2n.
$$
Pick a $\pi$-ample $\mathbb{Q}$-divisor $A$ such that $-(K_X+\Delta+A)$ is also
$\pi$-ample.  In particular $-(K_X+\Delta_i+A)$ is $\pi$-ample for $i\gg 0$.  Now
$$
A\cdot\Sigma_i=(K_X+\Delta_{i}+A)\cdot \Sigma_i-(K_X+\Delta_i)\cdot \Sigma_i<2n.
$$
It follows that the curves $\Sigma_i$ belong to a bounded family.  Thus, possibly passing
to a subsequence, we may assume that $\Sigma=\Sigma_i$ is constant.  In this case
\[
-(K_X+\Delta)\cdot \Sigma=\lim_i -(K_X+\Delta_i)\cdot \Sigma\leq 2n. \qedhere
\]
\end{proof}

\begin{corollary}\label{c_low} Let $\pi\colon\map X.U.$ be a projective morphism of 
normal quasi-projective varieties.  Suppose the pair $(X,\Delta=A+B)$ has log canonical
singularities, where $A\geq 0$ is an ample $\mathbb{R}$-divisor over $U$ and $B\geq 0$.
Suppose that there is a divisor $\Delta_0$ such that $K_X+\Delta_0$ is kawamata log
terminal.

Then there are only finitely many $(K_X+\Delta)$-negative extremal rays $\llist R.k.$ of
$\ccone X/U.$.
\end{corollary}
\begin{proof} We may assume that $A$ is a $\mathbb{Q}$-divisor.  Let $R$ be a
$(K_X+\Delta)$-negative extremal ray of $\ccone X/U.$.  Then
$$
-(K_X+B)\cdot R=-(K_X+\Delta)\cdot R+A\cdot R>0.
$$
By \eqref{t_low} $R$ is spanned by a curve $\Sigma$ such that
$$
-(K_X+B)\cdot \Sigma\leq 2n.
$$
But then 
\[
A\cdot \Sigma=-(K_X+B)\cdot \Sigma+(K_X+\Delta)\cdot \Sigma\leq 2n. \qedhere
\]
Therefore the curve $\Sigma$ belongs to a bounded family.
\end{proof}

\subsection{Effective base point free theorem}

\begin{theorem}[Effective Base Point Free Theorem]\label{t_base} Fix a positive integer $n$.  
Then there is a positive integer $m>0$ with the following property:

Let $f\colon\map X.U.$ be a projective morphism of normal quasi-projective varieties, and
let $D$ be a nef $\mathbb{R}$-divisor over $U$, such that $aD-(K_X+\Delta)$ is nef and big
over $U$, for some positive real number $a$, where $(X,\Delta)$ is kawamata log terminal
and $X$ has dimension $n$.

Then $D$ is semiample over $U$ and if $aD$ is Cartier then $maD$ is globally generated over
$U$.
\end{theorem}
\begin{proof} Replacing $D$ by $aD$ we may assume that $a=1$.  As the property that $D$ is
either semiample or globally generated over $U$ is local over $U$, we may assume that $U$
is affine.

By assumption we may write $D-(K_X+\Delta) \sim_{\mathbb{R},U} A+B$, where $A$ is an ample
$\mathbb{Q}$-divisor and $B\geq 0$.  Pick $\epsilon\in (0,1)$ such that
$(X,\Delta+\epsilon B)$ is kawamata log terminal.  Then
$$
D-(K_X+\Delta+\epsilon B)=(1-\epsilon)(D-(K_X+\Delta))+\epsilon (D-(K_X+\Delta+B)),
$$
is ample.  Replacing $(X,\Delta)$ by $(X,\Delta+\epsilon B)$ we may therefore assume that
$D-(K_X+\Delta)$ is ample.  Let $V$ be the subspace of $\WDiv_{\mathbb{R}}(X)$ spanned by
the components of $\Delta$.  As $\mathcal{L}(V)$ is a rational polytope, cf.
\eqref{l_polytopal}, which contains $\Delta$, we may find $\Delta'$ such that
$K_X+\Delta'$ is $\mathbb Q$-Cartier and kawamata log terminal, sufficiently close to
$\Delta$ so that $D-(K_X+\Delta')$ is ample.  Replacing $(X,\Delta)$ by $(X,\Delta')$ we
may therefore assume that $K_X+\Delta$ is $\mathbb{Q}$-Cartier.

The existence of the integer $m$ is Koll\'ar's effective version of the base point free theorem
\cite{Kollar93b}.

Pick a general ample $\mathbb{Q}$-divisor $A$ such that $D-(K_X+\Delta+A)$ is ample.
Replacing $\Delta$ by $\Delta+A$, we may assume that $\Delta=A+B$ where $A$ is ample and
$B\geq 0$.  By \eqref{c_low} there are finitely many $(K_X+\Delta)$-negative extremal rays
$\llist R.k.$ of $\ccone X.$.  Let
$$
F=\{\, \alpha\in \ccone X.\,|\, D\cdot \alpha=0 \,\}.
$$
Then $F$ is a face of $\ccone X.$ and if $\alpha\in F$ then $(K_X+\Delta)\cdot \alpha <0$,
and so $F$ is spanned by a subset of the extremal rays $\llist R.k.$.  Let $V$ be the
smallest affine subspace of $\WDiv_{\mathbb{R}}(X)$, which is defined over the rationals
and contains $D$.  Then
$$
\mathcal{C}=\{\, B\in V \,|\, B\cdot \alpha=0, \forall\alpha\in F \,\},
$$
is a rational polyhedron.  It follows that we may find positive real numbers $\llist r.q.$
and nef $\mathbb{Q}$-Cartier divisors $\llist D.q.$ such that $D=\sum r_pD_p$.  Possibly
re-choosing $\llist D.q.$ we may assume that $D_p-(K_X+\Delta)$ is ample.  By the usual
base point free theorem, $D_p$ is semiample and so $D$ is semiample.  \end{proof}

\begin{corollary}\label{c_base} Fix a positive integer $n$.  Then there is a constant
$m>0$ with the following property:

Let $f\colon\map X.U.$ be a projective morphism of normal quasi-projective varieties such
that $K_X+\Delta$ is nef over $U$ and $\Delta$ is big over $U$, where $(X,\Delta)$ is
kawamata log terminal and $X$ has dimension $n$.

Then $K_X+\Delta$ is semiample over $U$ and if $r$ is a positive constant such that
$r(K_X+\Delta)$ is Cartier, then $mr(K_X+\Delta)$ is globally generated over $U$.
\end{corollary}
\begin{proof} By \eqref{l_big-to-ample} we may find a kawamata log terminal pair
$$
K_X+A+B\sim_{\mathbb{R},U} K_X+\Delta,
$$ 
where $A\geq 0$ is a general ample $\mathbb{Q}$-divisor over $U$ and $B\geq 0$.  As
$$
(K_X+\Delta)-(K_X+B)\sim_{\mathbb{R},U} A,
$$
is ample over $U$ and $K_X+B$ is kawamata log terminal, the result follows by
\eqref{t_base}.  \end{proof}

\begin{lemma}\label{l_weak} Let $\pi\colon\map X.U$ be a projective morphism
of quasi-projective varieties.  Suppose that $(X,\Delta)$ is a kawamata log terminal 
pair, where $\Delta$ is big over $U$.

If $\phi\colon\rmap X.Y.$ is a weak log canonical model of $K_X+\Delta$ over $U$ then
\begin{itemize} 
\item $\phi$ is a semiample model over $U$,
\item the ample model $\psi\colon\rmap X.Z.$ of $K_X+\Delta$ over $U$ exists, and 
\item there is a contraction morphism $h\colon\map Y.Z.$ such that
$K_Y+\Gamma\sim _{\mathbb R, U}h^*H$, for some ample divisor $H$, where $\Gamma=\phi_*\Delta$.  
\end{itemize} 
\end{lemma}
\begin{proof} $K_Y+\Gamma$ is semiample over $U$ by \eqref{c_base} and (3) of
\eqref{l_ample} implies the rest.
\end{proof}

\subsection{The MMP with scaling}
In order to run a minimal model program, two kinds of operations, known as flips and
divisorial contractions are required.  We begin by recalling their definitions.

\begin{definition}\label{d_flip} Let $(X,\Delta)$ be a log canonical pair and 
$f\colon\map X.Z.$ be a projective morphism of normal varieties.  Then $f$ is a
\textbf{flipping contraction} if
\begin{enumerate}
\item $X$ is $\mathbb{Q}$-factorial and $\Delta$ is an $\mathbb{R}$-divisor,
\item $f$ is a small birational morphism of relative Picard number $\rho(X/Z)=1$, and
\item $-(K_X+\Delta)$ is $f$-ample.
\end{enumerate}

The \textbf{flip} $f^+\colon\map X^+.Z.$ of a flipping contraction $f\colon\map X.Z.$ is a
small birational projective morphism of normal varieties $f^+\colon\map X^+.Z.$ such that
$K_{X^+}+\Delta^+$ is $f^+$-ample, where $\Delta^+$ is the strict transform of $\Delta$.
\end{definition}

\begin{lemma}\label{l_flip} Let $(X,\Delta)$ be a log canonical pair, let 
$f\colon\map X.Z.$ be a flipping contraction and let $f^+\colon\map X^+.Z.$ be the flip of
$f$.  Let $\phi\colon\rmap X.X^+.$ be the induced birational map.

Then 
\begin{enumerate} 
\item $\phi$ is the log canonical model of $(X,\Delta)$ over $Z$,  
\item $(X^+,\Delta^+=\phi _*\Delta)$ is log canonical, 
\item $X^+$ is $\mathbb{Q}$-factorial, and
\item $\rho(X^+/Z)=1$. 
\end{enumerate} 
\end{lemma}
\begin{proof} As $\phi$ is small, \eqref{l_negativity} implies that $\phi$ is
$(K_X+\Delta)$-negative and this implies (1) and (2).

Let $B$ be a divisor on $X^+$ and let $C\subset X$ be the strict transform of $B$.  Since
$\rho(X/Z)=1$, we may find $\lambda\in \mathbb{R}$ such $C+\lambda(K_X+\Delta)$ is
numerically trivial over $U$.  But then \eqref{t_base} implies that
$C+\lambda(K_X+\Delta)=f^*D$, for some $\mathbb{R}$-Cartier divisor $D$ on $Z$.  Therefore
$B+\lambda(K_{X^+}+\Delta^+)=f^{+*}D$, and this implies (3) and (4).
\end{proof} 

In terms of our induction, we will need to work with a more restrictive notion of flipping
contraction.
\begin{definition}\label{d_pl} Let $(X,\Delta)$ be a purely log terminal
pair and let $f\colon\map X.Z.$ be a projective morphism of normal varieties.  Then $f$ is
a \textbf{pl-flipping contraction} if
\begin{enumerate}
\item $X$ is $\mathbb{Q}$-factorial and $\Delta$ is an $\mathbb{R}$-divisor,
\item $f$ is a small birational morphism of relative Picard number $\rho (X/Z)=1$,
\item $-(K_X+\Delta)$ is $f$-ample, and
\item $S=\rdown \Delta.$ is irreducible and $-S$ is $f$-ample.
\end{enumerate}

A \textbf{pl-flip} is the flip of a pl-flipping contraction. 
\end{definition}

\begin{definition}\label{d_contraction} Let $(X,\Delta)$ be a log canonical pair and 
$f\colon\map X.Z.$ be a projective morphism of normal varieties.  Then $f$ is a
\textbf{divisorial contraction} if
\begin{enumerate}
\item $X$ is $\mathbb{Q}$-factorial and $\Delta$ is an $\mathbb{R}$-divisor,
\item $f$ is a birational morphism of relative Picard number $\rho(X/Z)=1$ with
exceptional locus a divisor, and
\item $-(K_X+\Delta)$ is $f$-ample.
\end{enumerate}
\end{definition}
\begin{remark}\label{r_div}If $f\colon\map X.Z.$ is a divisorial contraction, then 
an argument similar to \eqref{l_flip} shows that $(Z,f_*\Delta)$ is log canonical and $Z$
is $\mathbb{Q}$-factorial.
\end{remark}

\begin{definition}\label{d_mori-fibre-space} Let $(X,\Delta)$ be a log canonical pair and 
$f\colon\map X.Z.$ be a projective morphism of normal varieties. Then $f$
is a \textbf{Mori fibre
  space} if
\begin{enumerate}
\item $X$ is $\mathbb{Q}$-factorial and $\Delta$ is an $\mathbb{R}$-divisor,
\item $f$ is a contraction morphism, $\rho(X/Z)=1$ and $\dim Z< \dim X$, and
\item $-(K_X+\Delta)$ is $f$-ample.
\end{enumerate}
\end{definition}

The objective of the MMP is to produce either a log terminal model or a Mori fibre space.
Note that if $K_X+\Delta$ has a log terminal model then $K_X+\Delta$ is pseudo-effective
and if $K_X+\Delta$ has a Mori fibre space then $K_X+\Delta$ is not pseudo-effective, so
these two cases are mutually exclusive.  

There are several versions of the MMP, depending on the singularities that are allowed
(typically, one restricts to kawamata log terminal singularities or divisorially log
terminal singularities or terminal singularities with $\Delta=0$) and depending on the
choices of negative extremal rays that are allowed (traditionally any choice of an
extremal ray is acceptable).

In this paper, we will run the MMP with scaling for divisorially log terminal pairs
satisfying certain technical assumptions.  We will need the following key result.

\begin{lemma}\label{l_touch} Let $\pi\colon\map X.U.$ be a projective morphism of 
normal quasi-projective varieties.  Suppose that the pair $(X,\Delta=A+B)$ has kawamata
log terminal singularities, where $A\geq 0$ is big over $U$, $B\geq 0$, $D$ is a nef
$\mathbb{R}$-Cartier divisor over $U$, but $K_X+\Delta$ is not nef over $U$.  Set
$$
\lambda=\sup \,\{\,\mu\, | \, D+\mu (K_X+\Delta) \text{ is nef over $U$}\,\}.
$$ 

Then there is a $(K_X+\Delta)$-negative extremal ray $R$ over $U$, such that
$$
(D+\lambda(K_X+\Delta))\cdot R=0.
$$
\end{lemma}
\begin{proof} By \eqref{l_big-to-ample} we may assume that $A$ is ample over $U$.    

By \eqref{c_low} there are only finitely many $(K_X+\Delta)$-negative extremal rays
$\llist R.k.$ over $U$.  For each $(K_X+\Delta)$-negative extremal ray $R_i$, pick a curve
$\Sigma_i$ which generates $R_i$.  Let
$$
\mu=\min _i \frac {D\cdot \Sigma_i}{-(K_X+\Delta)\cdot \Sigma_i}.
$$
Then $D+\mu (K_X+\Delta)$ is nef over $U$, since it is non-negative on each $R_i$, but it
is zero on one of the extremal rays $R=R_i$.  Thus $\lambda=\mu$.  \end{proof}

\begin{lemma}\label{l_mmp} Let $\pi\colon\map X.U.$ be a projective morphism of 
normal quasi-projective varieties and let $(X,\Delta_0)$ be a kawamata log terminal pair.
Suppose that $(X,\Delta=A+B)$ is a log canonical pair, where $A\geq 0$ is big over $U$,
$B\geq 0$, $\mathbf{B}_+(A/U)$ contains no non kawamata log terminal centres of $(X,\Delta)$ and $C$
is an $\mathbb{R}$-Cartier divisor such that $K_X+\Delta$ is not nef over $U$, whilst
$K_X+\Delta+C$ is nef over $U$.

Then there is a $(K_X+\Delta)$-negative extremal ray $R$ and a real number $0<\lambda\leq
1$ such that $K_X+\Delta+\lambda C$ is nef over $U$ but trivial on $R$.
\end{lemma}
\begin{proof} By \eqref{l_big-to-ample} we may assume $K_X+\Delta$ is kawamata log
terminal and that $A$ is ample over $U$.  Apply \eqref{l_touch} to
$D=K_X+\Delta+C$. \end{proof}

\begin{remark}\label{r_special} Assuming existence and termination of the 
relevant flips, we may use \eqref{l_mmp} to define a special minimal model program, which
we will refer to as the {\bf $(K_X+\Delta)$-MMP with scaling of $C$}.

Let $\pi\colon\map X. U.$ be a projective morphism of normal quasi-projective varieties,
where $X$ is $\mathbb{Q}$-factorial, $(X,\Delta +C=S+A+B+C)$ is a divisorially log
terminal pair, such that $\rdown \Delta .=S$, $A\geq 0$ is big over $U$,
$\mathbf{B}_+(A/U)$ does not contain any non kawamata log terminal centres of $(X,\Delta +C)$, and
$B\geq 0$, $C\geq 0$.  We pick $\lambda\geq 0$ and a $(K_X+\Delta)$-negative extremal ray
$R$ over $U$ as in \eqref{l_mmp} above.  If $K_X+\Delta$ is nef over $U$ we stop.
Otherwise $\lambda>0$ and we let $f\colon\map X.Z.$ be the extremal contraction over $U$
defined by $R$.  If $f$ is not birational, we have a Mori fibre space over $U$ and we
stop.  If $f$ is birational then either $f$ is divisorial and we replace $X$ by $Z$ or $f$
is small and assuming the existence of the flip $f^+\colon\map X^+.Z.$, we replace $X$ by
$X^+$. In either case $K_X+\Delta+\lambda C$ is nef over $U$ and $K_X+\Delta$ is
divisorially log terminal and so we may repeat the process.

In this way, we obtain a sequence $\phi_i\colon\rmap X_i.X_{i+1}.$ of $K_X+\Delta$ flips
and divisorial contractions over $U$ and real numbers $1\geq \lambda _1\geq \lambda _2
\geq \cdots$ such that $K_{X_{i}}+\Delta _{i}+\lambda _{i}C_{i}$ is nef over $U$ where
$\Delta _{i}=(\phi _{i-1})_{*}\Delta _{i-1}$ and $C_{i}=(\phi _{i-1})_{*}C _{i-1}$.

Note that by \eqref{l_centre}, each step of this MMP preserves the condition that
$(X,\Delta)$ is divisorially log terminal and $\mathbf{B}_+(A/U)$ does not contain any non
kawamata log terminal centres of $(X,\Delta)$ so that we may apply \eqref{l_mmp}.  By
\eqref{l_flip} and \eqref{r_div}, $(X,\Delta+\lambda C)$ is log canonical.  However the
condition that $\mathbf{B}_+(A/U)$ does not contain any non kawamata log terminal centres
of $(X,\Delta +\lambda C)$ is not necessarily preserved.
\end{remark}

\begin{lemma}\label{l_centre} Let $\pi\colon\map X.U.$ be a projective morphism
of normal quasi-projective varieties.  Suppose that $K_X+\Delta$ is divisorially log
terminal, $X$ is $\mathbb{Q}$-factorial and let $\phi\colon\rmap X.Y.$ be a sequence of
steps of the $(K_X+\Delta)$-MMP over $U$.

If $\Gamma=\phi_*\Delta$ then
\begin{enumerate} 
\item $\phi$ is an isomorphism at the generic point of every non kawamata log terminal centre of
$K_Y+\Gamma$.  In particular $(Y,\Gamma)$ is divisorially log terminal.
\item If $\Delta=S+A+B$, where $S=\rdown \Delta.$, $A\geq 0$ is big over $U$,
$\mathbf{B}_+(A/U)$ does not contain any non kawamata log terminal centres of $K_X+\Delta$, and $B\geq
0$ then $\phi_*S=\rdown\Gamma.$, $\phi_*A$ is big over $U$ and $\mathbf{B}_+(\phi_*A/U)$
does not contain any non kawamata log terminal centres of $K_Y+\Gamma$.

In particular $\Gamma\sim_{\mathbb{R},U}\Gamma'$ where $(Y,\Gamma')$ is kawamata log
terminal and $\Gamma '$ is big over $U$.
\end{enumerate} 
\end{lemma}
\begin{proof} We may assume that $\phi$ is either a flip or a divisorial contraction over
$U$. We first prove (1).  Let $p\colon\map W.X.$ and $q\colon\map W.Y.$ be common log
resolutions, which resolve the indeterminacy of $\phi$.  We may write
$$
p^*(K_X+\Delta)=q^*(K_Y+\Gamma)+E,
$$
where $E$ is exceptional and contains every exceptional divisor over the locus where
$\phi^{-1}$ is not an isomorphism.  In particular the log discrepancy of every valuation
with centre on $Y$ contained in the locus where $\phi^{-1}$ is not an isomorphism with
respect to $K_Y+\Gamma$ is strictly greater than the log discrepancy with respect to
$K_X+\Delta$.  Hence (1).

Now suppose that $\Delta=S+A+B$, $A$ is big over $U$ and no non kawamata log terminal centre of
$(X,\Delta)$ is contained in $\mathbf{B}_+(A/U)$.  Pick a divisor $C$ on $Y$ which is a
general ample $\mathbb{Q}$-divisor over $U$.  Possibly replacing $C$ by a smaller
multiple, we may assume that $\mathbf{B}((A-\phi^{-1}_*C)/U)$ does not contain any log
canonical centres of $(X,\Delta)$.  Thus $A-\phi^{-1}_*C \sim_{\mathbb{R},U} D\geq 0$,
where $D$ does not contain any non kawamata log terminal centres of $(X,\Delta)$.  If $\epsilon>0$ is
sufficiently small, then $(X,\Delta'=\Delta+\epsilon D)$ is divisorially log terminal,
$(X,\Delta)$ and $(X,\Delta')$ have the same non kawamata log terminal centres and $\phi$ is a
$K_X+\Delta'$ flip or a divisorial contraction over $U$.  (1) implies that
$(Y,\Gamma'=\phi_*\Delta')$ is divisorially log terminal and hence $\phi_*D$ does not
contain any non kawamata log terminal centres of $(Y,\Gamma)$.  (2) follows as $\phi_*A-C
\sim_{\mathbb{R},U} \phi_*D\geq 0$.  \end{proof}

\begin{lemma}\label{l_valuation} Let $\pi\colon\map X.U.$ be a projective morphism 
of quasi-projective varieties.  Suppose that $(X,\Delta)$ is a divisorially log terminal
pair.  Let $S$ be a prime divisor.  

Let $f_i\colon\rmap X_i.X_{i+1}.$ a sequence of flips and divisorial contractions over
$U$, starting with $X_1:=X$, for the $(K_X+\Delta)$-MMP, which does not contract $S$.  If
$f_i$ is not an isomorphism in a neighbourhood of the strict transform $S_i$ of $S$, then
neither is the induced birational map $\rmap X_i.X_j.$, $j>i$.
\end{lemma}
\begin{proof} Since the map $f_i$ is $(K_{X_i}+\Delta_i)$-negative and $\rmap
X_{i+1}.X_j.$ is $(K_{X_{i+1}}+\Delta_{i+1})$-negative, there is some valuation $\nu$
whose centre intersects $S_i$, such that
$$
a(\nu,X_i,\Delta_i)<a(\nu,X_{i+1},\Delta_{i+1}) \quad \text{and} \quad a(\nu,X_{i+1},\Delta_{i+1})\leq a(\nu,X_j,\Delta_j).
$$
But then $a(\nu,X_i,\Delta_i) < a(\nu,X_j,\Delta_j)$, so that $\rmap X_i.X_j.$ is not an
isomorphism in a neighbourhood of $S_i$.
\end{proof}

\subsection{Shokurov's polytopes}

We will need some results from \cite{Shokurov96}.  First some notation.  Given a ray
$R\subset \ccone X.$, let
$$
R^{\perp}=\{\, \Delta\in \mathcal{L}(V) \,|\, (K_X+\Delta)\cdot R=0\,\}.
$$

\begin{theorem}\label{t_polytopal} Let $\pi\colon\map X.U.$ be a projective morphism 
of normal quasi-projective varieties.  Let $V$ be a finite dimensional affine subspace of
$\WDiv _{\mathbb{R}}(X)$, which is defined over the rationals.  Fix an ample
$\mathbb{Q}$-divisor $A$ over $U$.  Suppose that there is a kawamata log terminal pair
$(X,\Delta_0)$.

Then the set of hyperplanes $R^{\perp}$ is finite in $\mathcal{L}_A(V)$, as $R$ ranges
over the set of extremal rays of $\ccone X/U.$.  In particular, $\mathcal{N}_{A,\pi}(V)$
is a rational polytope.
\end{theorem}

\begin{corollary}\label{c_polytopal} Let $\pi\colon\map X.U.$ be a projective morphism 
of normal quasi-projective varieties.  Let $V$ be a finite dimensional affine subspace of
$\WDiv _{\mathbb{R}}(X)$, which is defined over the rationals.  Fix a general ample
$\mathbb{Q}$-divisor $A$ over $U$.  Suppose that there is a kawamata log terminal pair
$(X,\Delta_0)$.  Let $\phi\colon\rmap X.Y.$ be any birational contraction over $U$.

Then $\mathcal{W}_{\phi,A,\pi}(V)$ is a rational polytope.  Moreover there are finitely many
morphisms $f_i\colon\map Y.Z_i.$ over $U$, $1\leq i\leq k$, such that if $f\colon\map Y.Z.$
is any contraction morphism over $U$ and there is a divisor $D$ on $Z$ such that
$K_Y+\Gamma=\phi_*(K_X+\Delta)\sim_{\mathbb{R},U} f^*D$ for some $\Delta\in
\mathcal{W}_{\phi,A,\pi}(V)$, then there is an index $1\leq i\leq k$ and an isomorphism
$\eta\colon\map Z_i.Z.$ such that $f=\eta\circ f_i$.
\end{corollary}

\begin{corollary}\label{c_relative} Let $\pi\colon\map X.U.$ be a projective morphism 
of normal quasi-projective varieties.  Let $V$ be a finite dimensional affine subspace of
$\WDiv _{\mathbb{R}}(X)$, which is defined over the rationals.  Fix a general ample
$\mathbb{Q}$-divisor $A$ over $U$.  Let $(X,\Delta_0)$ be a kawamata log terminal pair,
let $f\colon\map X.Z.$ be a morphism over $U$ such that $\Delta_0\in \mathcal L _A(V)$ and
$K_X+\Delta_0\sim _{\mathbb R ,U}f^*H$, where $H$ is an ample divisor over $U$.  Let $\phi\colon\rmap X.Y.$
be a birational map over $Z$.

Then there is a neighbourhood $P_0$ of $\Delta_0$ in $\mathcal{L}_A(V)$ such that for all
$\Delta\in P_0$ if $\phi$ is a log terminal model for $K_X+\Delta$ over $Z$ then $\phi$ is
a log terminal model for $K_X+\Delta$ over $U$.
\end{corollary}

\begin{proof}[Proof of \eqref{t_polytopal}] Since $\mathcal{L}_A(V)$ is compact it
suffices to prove this locally about any point $\Delta\in \mathcal{L}_A(V)$.  By
\eqref{l_linear-dlt} we may assume that $K_X+\Delta$ is kawamata log terminal.  Fix
$\epsilon>0$ such that if $\Delta'\in \mathcal{L}_A(V)$ and $\|\Delta'-\Delta\|<\epsilon$,
then $\Delta'-\Delta+A/2$ is ample over $U$.  Let $R$ be an extremal ray over $U$ such
that $(K_X+\Delta')\cdot R=0$, where $\Delta'\in \mathcal{L}_A(V)$ and
$\|\Delta'-\Delta\|<\epsilon$.  We have
$$
(K_X+\Delta-A/2)\cdot R=(K_X+\Delta')\cdot R-(\Delta'-\Delta+A/2)\cdot R<0.
$$
Finiteness then follows from \eqref{c_low}.  

$\mathcal{N}_{A,\pi}(V)$ is surely a closed subset of $\mathcal{L}_A(V)$.  If $K_X+\Delta$
is not nef over $U$ then \eqref{t_low} implies that $K_X+\Delta$ is negative on a rational
curve $\Sigma$ which generates an extremal ray $R$ of $\ccone X/U.$.  Thus
$\mathcal{N}_{A,\pi}(V)$ is the intersection of $\mathcal{L}_A(V)$ with the half-spaces
determined by finitely many of the extremal rays of $\ccone X/U.$.  \end{proof}

\begin{proof}[Proof of \eqref{c_polytopal}] Since $\mathcal{L}_A(V)$ is a rational
polytope its span is an affine subspace of $V_A$, which is defined over the rationals.
Possibly replacing $V$, we may therefore assume that $\mathcal{L}_A(V)$ spans $V_A$.  By
compactness, to prove that $\mathcal{W}_{\phi,A,\pi}(V)$ is a rational polytope, we may
work locally about a divisor $\Delta\in \mathcal{W}_{\phi,A,\pi}(V)$.  By
\eqref{l_linear-dlt} we may assume that $K_X+\Delta$ is kawamata log terminal in which
case $K_Y+\Gamma$ is kawamata log terminal as well.  Let $W\subset \WDiv_{\mathbb R}(Y)$
be the image of $V$.  If $C=\phi_*A$ then $C$ is big over $U$ and by
\eqref{l_linear-dlt}, we may find a rational affine linear isomorphism $L\colon\map
W.W.'$ and an ample $\mathbb{Q}$-divisor $C'$ over $U$ such that $L(\Gamma)$ belongs to
the interior of $\mathcal{L}_{C'}(W')$ and $L(\Gamma) \sim_{\mathbb{Q},U} \Gamma$ for any
$\Gamma\in W$.  \eqref{t_polytopal} implies that $\mathcal{N}_{C',\psi}(W')$ is a rational
polytope, where $\psi\colon\map Y.U.$ is the structure morphism, and so
$\mathcal{N}_{C,\psi}(W)$ is a rational polytope locally about $\Gamma$.

Let $p\colon\map Z.X.$ be a log resolution of $(X,\Delta)$ which resolves the
indeterminacy locus of $\phi$, via a birational map $q\colon\map Z.Y.$.  We may write
\begin{align*} 
K_Z+\Psi &= p^*(K_X+\Delta)\\ 
K_Z+\Phi &= q^*(K_Y+\Gamma).
\end{align*} 
Note that $\Delta\in\mathcal{W}_{\phi,A,\pi}(V)$ if and only if $\Gamma=\phi_*\Delta\in
\mathcal{N}_{C,\psi}(W)$ and $\Psi-\Phi\geq 0$.  Since the map $L\colon\map V.W.$ given by
$\map \Delta.\Gamma=\phi_*\Delta.$ is rational linear, the first statement is clear.

Note that if $f\colon\map Y.Z.$ and $f'\colon\map Y.Z'.$ are two contraction morphisms
over $U$ then there is an isomorphism $\eta\colon\map Z.Z'.$ such that $f'=\eta\circ f$ if
and only if the curves contracted by $f$ and $f'$ coincide. 

Let $f\colon\map Y.Z.$ be a contraction morphism over $U$, such that
$$
K_Y+\Gamma=K_Y+\phi_*\Delta \sim_{\mathbb{R},U} f^*D,
$$ 
where $\Delta\in\mathcal{W}_{\phi,A,\pi}(V)$.  $\Gamma$ belongs to the interior of a unique
face $G$ of $\mathcal{N}_{C,\psi}$ and the curves contracted by $f$ are determined by $G$.
Now $\Delta$ belongs to the interior of a unique face $F$ of $\mathcal{W}_{\phi,A,\pi}(V)$
and $G$ is determined by $F$.  But as $\mathcal{W}_{\phi,A,\pi}(V)$ is a rational polytope
it has only finitely many faces $F$.  \end{proof}

\begin{proof}[Proof of \eqref{c_relative}] By \eqref{t_polytopal} we may find finitely
many extremal rays $\llist R.k.$ of $Y$ over $U$ such that if
$K_Y+\Gamma=K_Y+\phi_*\Delta$ is not nef over $U$, then it is negative on one of these
rays.  If $\Gamma_0=\phi_*\Delta_0$ then we may write
$$
K_Y+\Gamma=K_Y+\Gamma_0+(\Gamma-\Gamma_0) \sim_{\mathbb{R},U} g^*H+\phi_*(\Delta-\Delta_0),
$$ 
where $g\colon\map Y.Z.$ is the structure morphism.  Therefore there is a neighbourhood
$P_0$ of $\Delta_0$ in $\mathcal{L}_A(V)$ such that if $K_Y+\Gamma$ is not nef over $U$
then it is negative on an extremal ray $R_i$, which is extremal over $Z$.  In particular
if $\phi$ is a log terminal model of $K_X+\Delta$ over $Z$ then it is a log terminal model
over $U$.  The other direction is clear.  
\end{proof}

\makeatletter
\renewcommand{\thetheorem}{\thesection.\arabic{theorem}}
\@addtoreset{theorem}{section} 
\makeatother
\section{Special finiteness}\label{s_special}

\begin{lemma}\label{l_terminal-negative} Let $\pi\colon\map X.U.$ be a projective morphism 
of quasi-projective varieties.  Let $(X,\Delta=S+A+B)$ be a log smooth pair, where
$S=\rdown\Delta.$ is a prime divisor, $A$ is a general ample $\mathbb{Q}$-divisor over $U$
and $B\geq 0$.  Let $C=A|_S$ and $\Theta=(\Delta-S)|_S$.  Let $\phi\colon\rmap X.Y.$ be a
birational map over $U$ which does not contract $S$, let $T$ be the strict transform of
$S$ and let $\tau\colon\rmap S.T.$ be the induced birational map.  Let
$\Gamma=\phi_*\Delta$ and define $\Psi$ on $T$ by adjunction,
$$
(K_Y+\Gamma)|_T=K_T+\Psi.
$$
If $(S,\Theta)$ is terminal and $\phi$ is a weak log canonical model of $K_X+\Delta$ over
$U$ then there is a divisor $C\leq\Xi\leq \Theta$ such that $\tau$ is a weak log
canonical model of $K_S+\Xi$ over
$U$, where $\tau_*\Xi=\Psi$.  
\end{lemma}
\begin{proof} Let $p\colon\map W.X.$ and $q\colon\map W.Y.$ be log resolutions of
$(X,\Delta)$ and $(Y,\Gamma)$, which resolve the indeterminacy of $\phi$, where
$\Gamma=\phi_*\Delta$.  Then we may write
$$
K_W+\Delta'=p^*(K_X+\Delta)+E \quad \text{and} \quad K_W+\Gamma'=q^*(K_Y+\Gamma)+F,
$$
where $\Delta'\geq 0$ and $E\geq 0$ have no common components, $\Gamma'\geq 0$ and $F\geq
0$ have no common components, $p_*\Delta'=\Delta$, $q_*\Gamma'=\Gamma$, $p_*E=0$ and
$q_*F=0$.  Since $\phi$ is $(K_X+\Delta)$-non-positive, we have $F-E+\Delta'-\Gamma'\geq 0$
and so
$$
F\geq E \qquad \text{and} \qquad \Gamma'\leq \Delta'.
$$
If $R$ is the strict transform of $S$ on $W$ then there are two birational morphisms
$f=p|_R\colon\map R.S.$ and $g=q|_R\colon\map R.T.$.  Since $(X,\Delta)$ is purely log
terminal, $(Y,\Gamma)$ is purely log terminal.  In particular $(T,\Psi)$ is kawamata log
terminal.  It follows that
$$
K_R+\Theta'=f^*(K_S+\Theta)+E' \quad \text{and} \quad K_R+\Psi'=g^*(K_T+\Psi)+F',
$$
where $\Theta'=(\Delta'-R)|_R$, $E'=E|_R$, $\Psi'=(\Gamma'-R)|_R$ and $F'=F|_R$.
Moreover, every component of $F'$ is $g$-exceptional, by (4) of \eqref{d_adjunction}.  As
$p$ and $q$ are log resolutions
$$
F'\geq E' \qquad \text{and} \qquad \Psi'\leq \Theta'.
$$

Suppose that $B$ is a prime divisor on $R$ which is $f$-exceptional but not
$g$-exceptional.  Then $B$ is not a component of $F'$, and so it is not a component of
$E'$.  But then the log discrepancy of $B$ with respect to $(S,\Theta)$ is at most one,
which contradicts the fact that $(S,\Theta)$ is terminal.  Thus $\tau$ is a birational
contraction.  

As $g_*E'=0$, 
$$
\tau_*\Theta=g_*\Theta'\geq g_*\Psi'=\Psi.
$$  
Let $\Xi\leq\Theta$ be the biggest divisor on $S$ such that $\tau_*\Xi=\Psi$.  In other
words, if a component of $\Theta$ is not $\tau$-exceptional, we replace its coefficient by
the corresponding coefficient of $\Psi$ and if a component of $\Theta$ is
$\tau$-exceptional, then we do not change its coefficient.  As $A$ is a general ample
$\mathbb{Q}$-divisor over $U$, $C$ is not contained in the locus where $\phi$ is not an
isomorphism.  It follows that $C\leq\Xi$.  We have
$$
f^*(K_S+\Xi)=g^*(K_T+\Psi)+L, 
$$
where $g_*L=0$.  Contrast this with
$$
f^*(K_S+\Theta)=g^*(K_T+\Psi)+M.
$$
We have already seen that $M=F'-E'+\Theta '-\Psi '\geq 0$.  By definition of $\Xi$, $f_*L$ and $f_*M$
agree on the $\tau$-exceptional divisors.  Since $f_*L$ is $\tau$-exceptional, it follows
that $f_*L\geq 0$.  But then \eqref{l_enough} implies that $L\geq 0$ and so $\tau$ is a
weak log canonical model of $K_S+\Xi$ over $U$.  \end{proof}

\begin{lemma}\label{l_linear-system} Let $\map X.U.$ and $\map Y.U.$ be two 
projective morphisms of normal quasi-projective varieties, where $X$ and $Y$ are
$\mathbb{Q}$-factorial.  Let $f\colon\rmap X.Y.$ be a small birational map over $U$.  Let
$S$ be a prime divisor on $X$, let $T$ be its strict transform on $Y$, let $g\colon\rmap
S.T.$ be the induced birational map and let $H$ be an ample $\mathbb{R}$-divisor on $Y$
over $U$.  Let $D$ be the strict transform of $H$.

If $g$ is an isomorphism and $D|_S=g^*(H|_T)$ then $f$ is an isomorphism in a
neighbourhood of $S$ and $T$.
\end{lemma}
\begin{proof} If $p\colon\map W.X.$ is the normalisation of the graph of $f$ and
$q\colon\map W.Y.$ be the induced birational morphism, then we may write
$$
p^*D=q^*H+E,
$$
where $E$ is $p$-exceptional.  Let $R$ be the strict transform of $S$ in $W$ and let
$p_1=p|_R\colon\map R.S.$ and $q_1=q|_R\colon\map R.T.$.  Since
$$
p_1^*(D|_S)=(p^*D)|_R=(q^*H+E)|_R=q_1^*(H|_T)+E|_R,
$$
we have $E|_R=0$.  

Since $X$ is $\mathbb{Q}$-factorial the exceptional locus of $p$ is a divisor by (3) of
\eqref{l_negativity}.  Notice also that if $F$ is a $p$-exceptional divisor, then it is
covered by a family of $p$-exceptional curves. These curves are not contracted by $q$ and
so the image of $F$ is contained in the indeterminacy locus of $f$.

Let $V\subset X$ be a subvariety contained in the indeterminancy locus of $f$, which is
maximal with respect to inclusion.  If $F$ is a $p$-exceptional divisor such that
$V\subset p(F)$ then $V=p(F)$ and $F$ is contained in the support of $E$ by (2) of
\eqref{l_negativity}.  It follows that $V$ does not intersect $S$ and so $f$ is a morphism
in a neighbourhood of $S$.  Since $Y$ is $\mathbb{Q}$-factorial, it follows that $f$ is an
isomorphism in neighbourhood of $S$ and $T$.
\end{proof}

\begin{lemma}\label{l_determine-neighbourhood} Let $\pi\colon\map X.U.$ be a projective
morphism of normal quasi-projective varieties.  Let $(X,\Delta_i)$ be two $i=1$, $2$, purely log
terminal pairs, where $S=\rdown \Delta_i.$ is a prime divisor, which is independent of
$i$.  Let $\phi_i\colon\rmap X.Y_i.$ be $\mathbb{Q}$-factorial log canonical models of
$K_X+\Delta_i$ over $U$, where $\phi_i$ are birational and do not contract $S$ for $i=1$, $2$.  Let $T_i$ be the strict transform
of $S$ and let $\tau_i\colon\rmap S.T_i.$ be the induced birational maps.  Let $\Phi_i$ be
the different
$$
(K_{Y_i}+T_i)|_{T_i}=K_{T_i}+\Phi_i.
$$

If 
\begin{enumerate}
\item the induced birational map $\chi\colon\rmap Y_1.Y_2.$ is small,  
\item the induced birational map $\sigma\colon\rmap T_1.T_2.$ is an isomorphism, 
\item $\sigma^*\Phi_2=\Phi_1$, and 
\item for every component $B$ of the support of $(\Delta_2-S)$, we have
$$
(\phi_{1*}B)|_{T_1}=\sigma^*((\phi_{2*}B)|_{T_2}), 
$$
\end{enumerate} 
then $\chi$ is an isomorphism in a neighbourhood of $T_1$ and $T_2$.
\end{lemma}
\begin{proof} Let $\Gamma_i=\phi_{i*}\Delta_i$ and define $\Psi_i$ by adjunction,
$$
(K_{Y_i}+\Gamma_i)|_{T_i}=K_{T_i}+\Psi_i.
$$
We have 
\begin{align*} 
K_{T_i}+\Psi_i &= (K_{Y_i}+T_i+\Gamma_i-T_i)|_{T_i}\\ 
              &= K_{T_i}+\Phi_i+(\Gamma_i-T_i)|_{T_i},
\end{align*} 
so that $\Psi_i=\Phi_i+(\phi_{i*}(\Delta_i-S))|_{T_i}$.  Conditions (3) and (4) then 
imply that 
$$
\sigma^*\Psi_2=\sigma^*(\Phi_2+(\phi_{2*}(\Delta_2-S))|_{T_2})=\Phi_1+(\phi_{1*}(\Delta_2-S))|_{T_1}.  
$$
It follows that 
\begin{align*} 
(\chi^*(K_{Y_2}+\Gamma_2))|_{T_1} &=(\chi^*\phi_{2*}(K_X+\Delta_2))|_{T_1} \\
                                &=(\phi_{1*}(K_X+\Delta_2))|_{T_1} \\
                                &=(K_{Y_1}+T_1+\phi_{1*}(\Delta_2-S))|_{T_1} \\
                                &=K_{T_1}+\Phi_1+(\phi_{1*}(\Delta_2-S))|_{T_1} \\
                                &=\sigma^*(K_{T_2}+\Psi_2).
\end{align*} 
Thus \eqref{l_linear-system} implies
that $\chi$ is an isomorphism in a neighbourhood of $T_1$ and $T_2$.  \end{proof}

\begin{lemma}\label{l_A-to-B} Theorem~\ref{t_finite}$_{n-1}$ implies
Theorem~\ref{t_special}$_n$.  
\end{lemma}
\begin{proof} Suppose not.  Then there would be an infinite sequence of
$\mathbb{Q}$-factorial weak log canonical models over $U$, $\phi_i\colon\rmap X.Y_i.$ for
$K_X+\Delta_i$ where $\Delta_i\in \mathcal{L}_{S+A}(V)$, which only contract components of
$\mathfrak{E}$ and which do not contract every component of $S$, such that if the induced
birational map $f_{ij}\colon\rmap Y_i.Y_j.$ is an isomorphism in a neighbourhood of the
strict transforms $S_i$ and $S_j$ of $S$ then $i=j$.  Since $\mathfrak{E}$ is finite,
possibly passing to a subsequence, we may assume that $f_{ij}$ is small.  Possibly passing
to a further subsequence, we may also assume that there is a fixed component $T$ of $S$
such that $\phi_i$ does not contract $T$ and $f_{ij}$ is not an isomorphism in a
neighbourhood of the strict transforms of $T$.  Replacing $S$ by $T$, we may therefore
assume that $S$ is irreducible.

Pick $\llist H.k.$ general ample over $U$ $\mathbb{Q}$-divisors which span $\WDiv_{\mathbb{R}}(X)$
modulo numerical equivalence over $U$ and let $H=\alist H.+.k.$ be their sum.  We may
replace $V$ by the subspace of $\WDiv_{\mathbb R}(X)$ generated by $V$ and $\llist H.k.$.
Passing to a subsequence, we may assume that $\lim \Delta_i=\Delta_{\infty}\in
\mathcal{L}_{S+A}(V)$.  By \eqref{l_linear-dlt} we may assume that $K_X+\Delta_{\infty}$
is purely log terminal and $\Delta_{\infty}$ is in the interior of $\mathcal L_{S+A}(V)$.
Possibly passing to a subsequence we may therefore assume that $\Delta_i$ and
$\Delta_{\infty}$ have the same support.  We may therefore assume that $K_X+\Delta_i$ is
purely log terminal and $\Delta_i$ contains the support of $H$.  By \eqref{l_ample-trick}
we may therefore assume that $\phi_i$ is the ample model of $K_X+\Delta_i$ over $U$.  In
particular $\Delta_i=\Delta_j$ implies $i=j$.

Pick $\Delta\in \mathcal{L}_{S+A}(V)$ such that $K_X+\Delta$ is purely log terminal and
$\Delta_i\leq \Delta$ for all $i>0$.  Let $f\colon\map Y.X.$ be a log resolution of
$(X,\Delta)$.  Then we may write
$$
K_Y+\Gamma=f^*(K_X+\Delta)+E,
$$
where $\Gamma\geq 0$ and $E\geq 0$ have no common components, $f_*\Gamma=\Delta$ and
$f_*E=0$.  Let $T$ be the strict transform of $S$.  Possibly blowing up more, we may
assume that $(T,\Theta)$ is terminal, where $\Theta=(\Gamma-T)|_T$.  We may find $F\geq 0$
exceptional so that $f^*A-F$ is ample over $U$, $(Y,\Gamma+F)$ is purely log terminal and
$(T,\Theta+F|_T)$ is terminal.  Let $A' \sim_{\mathbb Q,U} f^*A-F$ be a general ample
$\mathbb{Q}$-divisor over $U$.  For every $i$, we may write
$$
K_Y+\Gamma_i=f^*(K_X+\Delta_i)+E_i,
$$
where $\Gamma_i\geq 0 $ and $E_i\geq 0$ have no common components, $f_*\Gamma_i=\Delta_i$
and $f_*E_i=0$.  Let 
$$
\Gamma_i'=\Gamma_i-f^*A+F+A'\sim_{\mathbb{Q},U}\Gamma_i.
$$  
Then $K_Y+\Gamma_i'$ is purely log terminal and \eqref{l_equiv-models} implies that
$f\circ\phi_i$ is both a weak log canonical model, and the ample model, of $K_Y+\Gamma_i$ over
$U$.

Replacing $X$ by $Y$ we may therefore assume that $(X,\Delta)$ is log smooth and
$(S,\Theta)$ is terminal, where $\Theta=(\Delta-S)|_S$.  Let
$$
C=A|_S\leq \Xi_i\leq (\Delta_i-S)|_S\leq \Theta,
$$ 
be the divisors whose existence is guaranteed by \eqref{l_terminal-negative} so that
$\tau_i=\phi_i|_S\colon\rmap S.S_i.$ is a weak log canonical model (and the ample model)
of $K_S+\Xi_i$ over $U$, where $S_i$ is the strict transform of $S$.  As we are assuming
Theorem~\ref{t_finite}$_{n-1}$, possibly passing to a subsequence we may assume that the
restriction $g_{ij}=f_{ij}|_{S_i}\colon\rmap S_i.S_j.$ is an isomorphism.

Since $Y_i$ is $\mathbb{Q}$-factorial, we may define two $\mathbb R$-divisors on $S_i$ by adjunction,
$$
(K_{Y_i}+\Gamma_i)|_{S_i}=K_{S_i}+\Psi_i \qquad \text{and} \qquad (K_{Y_i}+S_i)|_{S_i}=K_{S_i}+\Phi_i,
$$
where $\Gamma_i=\phi_{i*}\Delta_i$.  Then 
$$
0\leq \Phi_i\leq \Psi_i=\tau_{i*}\Xi_i\leq \tau_{i*}\Theta.
$$
As the coefficients of $\Phi_i$ have the form $\frac{(r-1)}r$ for some integer $r>0$, by
(3) of \eqref{d_adjunction}, and the coefficients of $\Theta$ (and hence also of
$\tau_{i*}\Theta$) are all less than one, it follows that there are only finitely many
possibilities for $\Phi_i$.  Let $B$ be a component of the support of $\Delta_i-S$ and let
$G=(\phi_{i*}(B))|_{S_i}$.  By (3) of \eqref{d_adjunction}, the coefficients of $G$ are
integer multiples of $1/r$.  Pick a real number $\beta >0$ such that the coefficients of
$B$ in $\Delta_i-S$ are at least $\beta$.  As $0\leq \phi_{i*}(\Delta_i-S)|_{S_{i}}\leq
\tau_{i*}(\Theta)$, there are only finitely many possibilities for $G$.  Possibly passing
to a subsequence, we may assume $\Phi_i=f_{ij}^*\Phi_j$ and that if $B$ is a component of
the support of $\Delta_j-S$ then $(\phi_{i*}(B))|_{S_i}=g_{ij}^*((\phi_{j*}B)|_{S_j})$.
\eqref{l_determine-neighbourhood} implies that $f_{ij}$ is an isomorphism in a
neighbourhood of $S_i$ and $S_j$, a contradiction.  \end{proof}


\section{Log terminal models}
\label{s_models}

\begin{lemma}\label{l_termination} Assume Theorem~\ref{t_special}$_n$.

Let $\pi\colon\map X.U.$ be a projective morphism of normal quasi-projective varieties,
where $X$ is $\mathbb{Q}$-factorial of dimension $n$.  Suppose that
$$
K_X+\Delta+C=K_X+S+A+B+C,
$$ 
is divisorially log terminal and nef over $U$, where $S$ is a sum of prime divisors,
$\mathbf{B}_{+}(A/U)$ does not contain any non kawamata log terminal centres of $(X,\Delta+C)$ and
$B\geq 0$, $C\geq 0$.

Then any sequence of flips and divisorial contractions for the $(K_X+\Delta)$-MMP over $U$
with scaling of $C$, which does not contract $S$, is eventually disjoint from $S$.
\end{lemma}
\begin{proof} We may assume that $S$ is irreducible and by \eqref{l_big-to-ample} we may
assume that $A$ is a general ample $\mathbb{Q}$-divisor over $U$.  Let $f_i\colon\rmap
X_i.X_{i+1}.$ be a sequence of flips and divisorial contractions over $U$, starting with
$X_1:=X$, for the $(K_X+\Delta)$-MMP with scaling of $C$.

Let $\mathfrak{E}$ be the set of prime divisors in $X$ which are contracted by any of the
induced birational maps $\phi_i\colon\rmap X.X_i.$.  Then the cardinality of
$\mathfrak{E}$ is less than the relative Picard number of $X$ over $U$.  In particular
$\mathfrak{E}$ is finite.

By assumption there is a non-increasing sequence of real numbers $\ilist \lambda. \in
[0,1]$, such that $K_{X_i}+\Delta_i+\lambda_i C_i$ is nef over $U$, where $\Delta_i$ is
the strict transform of $\Delta$ and $C_i$ is the strict transform of $C$.  Further the
birational map $f_i$ is $(K_{X_i}+\Delta_i)$-negative and $(K_{X_i}+\Delta_i+\lambda_i
C_i)$-trivial, so that $f_i$ is $(K_{X_i}+\Delta_i+\lambda C_i)$-non-positive if and only
if $\lambda\leq \lambda_i$.  By induction, the birational map $\phi_i$ is therefore
$(K_X+\Delta+\lambda_iC)$-non-positive.  In particular $\phi_i$ is a
$\mathbb{Q}$-factorial weak log canonical model over $U$
of $(X,\Delta+\lambda_iC)$.

Let $V$ be the smallest affine subspace of $\WDiv_{\mathbb R}(X)$ containing $\Delta-S$
and $C$, which is defined over the rationals.  As $\Delta+\lambda_iC\in
\mathcal{L}_{S+A}(V)$, Theorem~\ref{t_special}$_n$ implies that there is an index $k$ and
infinitely many indices $l$ such that the induced birational map $\rmap X_k.X_l.$ is an
isomorphism in a neighbourhood of $S_k$ and $S_l$.  But then \eqref{l_valuation} implies
that $f_i$ is an isomorphism in a neighbourhood of $S_i$ for all $i\geq k$.
\end{proof}

We use \eqref{l_termination} to run a special MMP:
\begin{lemma}\label{l_subtract} Assume Theorem~\ref{t_existence}$_{n}$ and Theorem~\ref{t_special}$_n$.

Let $\pi\map\colon X.U.$ be a projective morphism of normal quasi-projective varieties,
where $X$ is $\mathbb{Q}$-factorial of dimension $n$.  Suppose that $(X,\Delta+C=S+A+B+C)$
is a divisorially log terminal pair, such that $\rdown\Delta.=S$, $A\geq 0$ is big over
$U$, $\mathbf{B}_{+}(A/U)$ does not contain any non kawamata log terminal centres of $(X,\Delta+C)$,
and $B\geq 0$, $C\geq 0$.  Suppose that there is an $\mathbb{R}$-divisor $D\geq 0$ whose
support is contained in $S$ and a real number $\alpha\geq 0$, such that
\[
\label{e_ast} K_X+\Delta \sim_{\mathbb{R},U} D+\alpha C. \tag{$\ast$}   
\]
If $K_X+\Delta+C$ is nef over $U$, then there is a log terminal model $\phi\colon\rmap
X.Y.$ for $K_X+\Delta$ over $U$, where $\mathbf{B}_+(\phi_*A/U)$ does not contain any log
canonical centres of $(Y,\Gamma=\phi_*\Delta)$.
\end{lemma}
\begin{proof} By \eqref{l_mmp} and \eqref{l_centre} we may run the $(K_X+\Delta)$-MMP with
scaling of $C$ over $U$, and this will preserve the condition that $\mathbf{B}_{+}(A/U)$
does not contain any non kawamata log terminal centres of $(X,\Delta)$.  Pick $t\in [0,1]$ minimal
such that $K_X+\Delta+tC$ is nef over $U$.  If $t=0$ we are done.  Otherwise we may find a
$(K_X+\Delta)$-negative extremal ray $R$ over $U$, such that $(K_X+\Delta+tC)\cdot R=0$.
Let $f\colon\map X.Z.$ be the associated contraction over $U$.  As $t>0$, $C\cdot R>0$ and
so $D\cdot R<0$.  In particular $f$ is always birational.

If $f$ is divisorial, then we can replace $X$, $S$, $A$, $B$, $C$ and $D$ by their images
in $Z$.  Note that \eqref{e_ast} continues to hold.  

Otherwise $f$ is small.  As $D\cdot R<0$, $R$ is spanned by a curve $\Sigma$ which is
contained in a component $T$ of $S$, where $T\cdot \Sigma<0$.  Note that
$K_X+S+A+B-\epsilon(S-T)$ is purely log terminal for any $\epsilon\in (0,1)$, and so $f$
is a pl-flip.  As we are assuming Theorem~\ref{t_existence}$_{n}$, the flip $f'\colon \map
X'.Z.$ of $f\colon\map X.Z.$ exists.  Again, if we replace $X$, $S$, $A$, $B$, $C$ and $D$
by their images in $X'$, then \eqref{e_ast} continues to hold.

On the other hand this MMP is certainly not an isomorphism in a neighbourhood of $S$ and
so the MMP terminates by \eqref{l_termination}.
\end{proof}

\begin{definition}\label{d_neutral} Let $\pi\colon\map X.U.$ be a projective morphism of normal quasi-projective 
varieties.  Let $(X,\Delta=A+B)$ be a $\mathbb{Q}$-factorial divisorially log terminal
pair and let $D$ be an $\mathbb{R}$-divisor, where $A\geq 0$, $B\geq 0$ and $D\geq 0$.  A
\textbf{neutral model} over $U$ for $(X,\Delta)$, with respect to $A$ and $D$, is any
birational map $f\colon\rmap X.Y.$ over $U$, such that
\begin{itemize} 
\item $f$ is a birational contraction, 
\item the only divisors contracted by $f$ are components of $D$,
\item $Y$ is $\mathbb{Q}$-factorial and projective over $U$, 
\item $\mathbf{B}_+(f_*A/U)$ does not contain any non kawamata log terminal centres of
$(Y,\Gamma=f_{*}\Delta )$, and
\item $K_Y+\Gamma$ is divisorially log terminal and nef over $U$.
\end{itemize}
\end{definition}

\begin{lemma}\label{l_smooth} Assume Theorem~\ref{t_existence}$_{n}$ and Theorem~ \ref{t_special}$_n$.

Let $\pi\colon\map X.U.$ be a projective morphism of normal quasi-projective varieties,
where $X$ has dimension $n$.  Let $(X,\Delta=A+B)$ be a divisorially log terminal log pair
and let $D$ be an $\mathbb{R}$-divisor, where $A\geq 0$ is big over $U$, $B\geq 0$ and
$D\geq 0$.

If 
\begin{itemize} 
\item[(i)] $K_X+\Delta\sim_{\mathbb{R},U} D$,
\item[(ii)] $(X,G)$ is log smooth, where $G$ is the support of $\Delta+D$, and
\item[(iii)] $\mathbf{B}_+(A/U)$ does not contain any non kawamata log terminal centres of $(X,G)$
\end{itemize} 
then $(X,\Delta)$ has a neutral model over $U$, with respect to $A$ and $D$.  
\end{lemma}
\begin{proof} We may write $D=D_1+D_2$, where every component of $D_1$ is a component of
$\rdown\Delta.$ and no component of $D_2$ is a component of $\rdown\Delta.$.  We proceed
by induction on the number of components of $D_2$.

Suppose $D_2=0$.  If $H$ is any general ample $\mathbb{Q}$-divisor over $U$, which is
sufficiently ample then $K_X+\Delta+H$ is divisorially log terminal and ample over $U$.
As the support of $D$ is contained in $\rdown\Delta.$, \eqref{l_subtract} implies that
$(X,\Delta)$ has a neutral model $f\colon\rmap X.Y.$ over $U$, with respect to $A$ and
$D$.

Now suppose that $D_2\neq 0$. Let
$$
\lambda =\sup \{\, t\geq 0 \,|\, \text{$(X,\Delta+tD_2)$ is log canonical}\,\},
$$
be the log canonical threshold of $D_2$.  Then $\lambda >0$ and $(X,\Theta=\Delta+\lambda
D_2)$ is divisorially log terminal and log smooth, $K_X+\Theta\sim_{\mathbb{R},U}
D+\lambda D_2$ and the number of components of $D+\lambda D_2$ that are not components of
$\rdown\Theta.$ is smaller than the number of components of $D_2$.  By induction there is
a neutral model $f\colon\rmap X.Y.$ over $U$ for $(X,\Theta)$, with respect to $A$ and
$D$.  

Now
\begin{align*} 
K_Y+f_*\Delta&\sim_{\mathbb{R},U} f_*D_1+f_*D_2, \\  
K_Y+f_*\Theta&=K_Y+f_*\Delta+\lambda f_* D_2,
\end{align*} 
where $K_Y+f_*\Theta$ is divisorially log terminal and nef over $U$, and the support of
$f_*D_1$ is contained in $\rdown f_*\Delta.$.  Since $\mathbf{B}_+(f_*A/U)$ does not
contain any non kawamata log terminal centres of $(Y,f_*\Theta)$, \eqref{l_subtract} implies that
$(Y,f_*\Delta)$ has a neutral model $g\colon\rmap Y.Z.$ over $U$, with respect to $f_*A$
and $f_*D$.  The composition $g\circ f\colon\rmap X.Z.$ is then a neutral model over $U$
for $(X,\Delta)$, with respect to $A$ and $D$.  \end{proof}

\begin{lemma}\label{l_contract} Let $\pi\colon\map X.U.$ be a projective morphism of normal
quasi-projective varieties.  Let $(X,\Delta=A+B)$ be a $\mathbb{Q}$-factorial divisorially
log terminal log pair and let $D$ be an $\mathbb{R}$-divisor, where $A\geq 0$ is big over
$U$, $B\geq 0$ and $D\geq 0$.

If every component of $D$ is either semiample over $U$ or a component of
$\mathbf{B}((K_X+\Delta)/U)$ and $f\colon\rmap X.Y.$ is a neutral model over $U$ for
$(X,\Delta)$, with respect to $A$ and $D$, then $f$ is a log terminal model for
$(X,\Delta)$ over $U$.
\end{lemma}
\begin{proof} By hypothesis the only divisors contracted by $f$ are components of
$\mathbf{B}((K_X+\Delta)/U)$.  Since the question is local over $U$, we may assume that
$U$ is affine.  Since $\mathbf{B}_+(f_*A/U)$ does not contain any non kawamata log terminal centres of
$(Y,\Gamma=f_*\Delta )$, \eqref{l_linear-ample} implies that we may find
$K_Y+\Gamma'\sim_{\mathbb{R},U}K_Y+\Gamma$ where $K_Y+\Gamma'$ is kawamata log terminal
and $\Gamma'$ is big over $U$.  \eqref{c_base} implies that $K_Y+\Gamma$ is semiample over
$U$.  

If $p\colon\map W.X.$ and $q\colon\map W.Y.$ resolve the indeterminacy of $f$ then we may
write
$$
p^*(K_X+\Delta)+E=q^*(K_Y+\Gamma)+F,
$$
where $E\geq 0$ and $F\geq 0$ have no common components, and both $E$ and $F$ are
exceptional for $q$.

As $K_Y+\Gamma$ is semiample, $\mathbf{B}((q^*(K_Y+\Gamma)+F)/U)$ and $F$ have the same
support.  On the other hand, every component of $E$ is a component of
$\mathbf{B}((p^*(K_X+\Delta)+E)/U)$.  Thus $E=0$ and any divisor contracted by $f$ is
contained in the support of $F$, and so $f$ is a log terminal model of
$(X,\Delta)$ over $U$.  \end{proof}

\begin{lemma}\label{l_B-to-C} Theorem~\ref{t_existence}$_{n}$ and
Theorem~\ref{t_special}$_n$ imply Theorem~\ref{t_model}$_n$.
\end{lemma}
\begin{proof} By \eqref{l_big-to-ample} we may assume that $\Delta=A+B$, where $A$ is a
general ample $\mathbb{Q}$-divisor over $U$ and $B\geq 0$.  By \eqref{p_mobile} we may
assume that $D=M+F$, where every component of $F$ is a component of $\mathbf{B}(D/U)$ and
there is a positive integer $m$ such that if $L$ is a component of $M$ then $mL$ is
mobile.

Pick a log resolution $f\colon\map Y.X.$ of the support of $D$ and $\Delta$, which
resolves the base locus of each linear system $|mL|$, for every component $L$ of $M$.  If
$\Phi$ is the divisor defined in \eqref{l_high} then every component of the exceptional
locus belongs to $\mathbf{B}((K_Y+\Phi)/U)$ and replacing $\Phi$ by an
$\mathbb{R}$-linearly equivalent divisor, we may assume that $\Phi$ contains an ample
divisor over $U$.  In particular, replacing $m\pi^*L$ by a general element of the linear
system $|m\pi^*L|$, we may assume that $K_Y+\Phi \sim_{\mathbb{R},U} N+G$, where every
component of $N$ is semiample, every component of $G$ is a component of
$\mathbf{B}((K_Y+\Phi)/U)$, and $(Y,\Phi+N+G)$ is log smooth.  By \eqref{l_high}, we
may replace $X$ by $Y$ and the result follows by \eqref{l_smooth} and
\eqref{l_contract}.  \end{proof}
\section{Non-vanishing}
\label{s_effective}

We follow the general lines of the proof of the non-vanishing theorem, see for example
Chapter 3, \S 5 of \cite{KM98}.  In particular there are two cases:

\begin{lemma}\label{l_bounded} Assume Theorem~\ref{t_model}$_n$.  Let $(X,\Delta)$ be a 
projective, log smooth pair of dimension $n$, where $\rdown\Delta.=0$, such that
$K_X+\Delta$ is pseudo-effective and 
$\Delta-A\geq 0$ for an ample $\mathbb{Q}$-divisor
$A$.  Suppose that for every positive integer $k$ such that $kA$ is integral,
$$
h^0(X,\ring X.(\rdown mk (K_X+\Delta).+ kA)),
$$
is a bounded function of $m$.  

 Then there is an $\mathbb{R}$-divisor $D$ such that $K_X+\Delta\sim_{\mathbb{R}}D\geq 0$.  
\end{lemma}
\begin{proof} By \eqref{p_linear} it follows that $K_X+\Delta$ is numerically equivalent
to $N_{\sigma}(K_X+\Delta)$.  Since $N_{\sigma}(K_X+\Delta)-(K_X+\Delta)$ is numerically
trivial and ampleness is a numerical condition, it follows that
$$
A''=A+N_{\sigma}(K_X+\Delta)-(K_X+\Delta),
$$
is ample and numerically equivalent to $A$.  Thus, since $A''$ is $\mathbb{R}$-linearly
equivalent to a positive linear combination of ample $\mathbb{Q}$-divisors, there exists
$0\le A'\sim_{\mathbb{R}}A''$ such that
$$
K_X+\Delta'=K_X+A'+(\Delta-A),
$$
is kawamata log terminal and numerically equivalent to $K_X+\Delta$, and
$$
K_X+\Delta'\sim_{\mathbb{R}} N_{\sigma}(K_X+\Delta)\geq 0.
$$
Thus by Theorem~\ref{t_model}$_n$, $K_X+\Delta'$ has a log terminal model $\phi\colon\rmap
X.Y.$, which, by \eqref{l_equiv-models}, is also a log terminal model for $K_X+\Delta$.  Replacing
$(X,\Delta)$ by $(Y,\Gamma)$ we may therefore assume that $K_X+\Delta$ is nef and the
result follows by the base point free theorem, cf. \eqref {c_base}.  \end{proof}

\begin{lemma}\label{l_unbounded} Let $(X,\Delta=A+B)$ be a projective, log smooth pair 
where $A$ is a general ample $\mathbb{Q}$-divisor and $\rdown B.=0$.  Suppose that there
is a positive integer $k$ such that $kA$ is integral and
$$
h^0(X,\ring X.(\rdown mk (K_X+\Delta).+ kA)),
$$
is an unbounded function of $m$.   

Then we may find a projective, log smooth pair $(Y,\Gamma)$ and a general ample
$\mathbb{Q}$-divisor $C$ on $Y$, where
\begin{itemize} 
\item $Y$ is birational to $X$, 
\item $\Gamma-C\geq 0$, 
\item $T=\rdown \Gamma.$ is an irreducible divisor, and
\item $\Gamma$ and $N_{\sigma}(K_Y+\Gamma)$ have no common components.  
\end{itemize} 

Moreover the pair $(Y,\Gamma)$ has the property that $K_X+\Delta \sim_{\mathbb{R}} D\geq
0$ for some $\mathbb{R}$-divisor $D$ if and only if $K_Y+\Gamma \sim_{\mathbb{R}} G\geq 0$ for some
$\mathbb{R}$-divisor $G$.
\end{lemma}
\begin{proof} Pick $m$ large enough so that
$$
h^0(X,\ring X.(\rdown mk(K_X+\Delta).+ kA)) > \binom{kn+n}{n}.
$$
By standard arguments, given any point $x\in X$, we may find a divisor $H'\geq 0$ which is
$\mathbb{R}$-linearly equivalent to
$$
\rdown mk(K_X+\Delta).+kA,
$$
of multiplicity greater than $kn$ at $x$.  In particular, we may find an
$\mathbb{R}$-Cartier divisor
$$
0\leq H\sim_{\mathbb{R}} m(K_X+\Delta)+A,
$$
of multiplicity greater than $n$ at $x$.  Given $t\in [0,m]$, consider
\begin{align*} 
(t+1)(K_X+\Delta) &= K_X+\frac {m-t}mA+B+t(K_X+\Delta+\frac 1mA) \\ 
                  &\sim_{\mathbb{R}}K_X+\frac {m-t}mA+B+\frac t mH \\ 
                  &=K_X+\Delta_t. 
\end{align*} 

Fix $0<\epsilon\ll 1$, let $A'=(\epsilon/m) A$ and $u=m-\epsilon$.  We have:
\begin{enumerate} 
\item $K_X+\Delta_0$ is kawamata log terminal,
\item $\Delta_t \ge A'$,  for any $t\in[0,u]$ and
\item the non kawamata log terminal locus of $(X,\Delta_u)$ contains a very general point
$x$ of $X$.
\end{enumerate} 

Let $\pi\colon\map Y.X.$ be a log resolution of $(X,\Delta+H)$.  We may write
$$
K_Y+\Psi_t=\pi^*(K_X+\Delta_t)+E_t,
$$
where $E_t\geq 0$ and $\Psi_t\geq 0$ have no common components, $\pi_*\Psi_t=\Delta_t$ and
$E_t$ is exceptional.  Pick an exceptional divisor $F\geq 0$ such that $\pi^*A'-F$ is
ample and let $C\sim_{\mathbb{Q}}\pi^*A'-F$ be a general ample $\mathbb{Q}$-divisor.  For
any $t\in [0,u]$, let
$$
\Phi_t=\Psi_t-\pi^*A'+C+F\sim_\mathbb{R} \Psi_t \qquad \text{and} \qquad \Gamma_t=\Phi_t-\Phi_t\wedge N_{\sigma}(K_Y+\Phi_t).
$$
Then properties (1-3) above become
\begin{enumerate} 
\item $K_Y+\Gamma_0$ is kawamata log terminal,
\item $\Gamma_t\ge C$,  for any $t\in[0,u]$, and 
\item $(Y,\Gamma_u)$ is not kawamata log terminal.   
\end{enumerate} 
Moreover
\begin{enumerate}
\setcounter{enumi}{3}
\item $(Y,\Gamma_t)$ is log smooth, for any $t\in [0,u]$, and 
\item $\Gamma_t$ and $N_{\sigma}(K_Y+\Gamma_t)$ do not have any common components. 
\end{enumerate} 

Let
$$
s=\sup \{\, t\in [0,u] \,|\, \text{$K_Y+\Gamma_t$ is log canonical} \,\}. 
$$
Note that 
\begin{align*} 
N_{\sigma}(K_Y+\Phi_t) &= N_{\sigma}(K_Y+\Psi_t) \\ 
                       &= N_{\sigma}(\pi^*(K_X+\Delta_t)) +E_{t}\\ 
                       &= N_{\sigma}((t+1)\pi^*(K_X+\Delta))+E_{t} \\ 
                       &= (t+1)N_{\sigma}(\pi^*(K_X+\Delta))+E_{t}.
\end{align*} 
Thus $K_Y+\Gamma_t$ is an affine linear function of $t$.  Setting $\Gamma=\Gamma_s$, we
may write
$$
\Gamma=T+C+B',
$$
where $T=\rdown\Gamma.\neq 0$, $C$ is ample and $B'\geq 0$.  Possibly perturbing $\Gamma$,
we may assume that $T$ is irreducible, so that $K_Y+\Gamma$ is purely log
terminal.  \end{proof}

We will need the following consequence of Kawamata-Viehweg vanishing:

\begin{lemma}\label{l_pseudo} Let $(X,\Delta=S+A+B)$ be a $\mathbb{Q}$-factorial projective
purely log terminal pair and let $m>1$ be an integer.  Suppose that
\begin{enumerate}
\item $S=\rdown \Delta.$ is irreducible,
\item $m(K_X+\Delta)$ is integral,
\item $m(K_X+\Delta)$ is Cartier in a neighbourhood of $S$,
\item $h^0(S,\ring S.(m(K_X+\Delta)))> 0$,
\item $K_X+G+B$ is kawamata log terminal, where $G\geq 0$,
\item $A\sim_{\mathbb{Q}}(m-1)tH+G$ for some $t$, and 
\item $K_X+\Delta+tH$ is big and nef.
\end{enumerate}

Then $h^0(X,\ring X.(m(K_X+\Delta)))>0$.
\end{lemma}
\begin{proof} Considering the long exact sequence associated to the restriction exact
sequence,
$$
\ses {\ring X.(m(K_X+\Delta)-S)}.{\ring X.(m(K_X+\Delta))}.{\ring S.(m(K_X+\Delta))}.,
$$
it suffices to observe that
$$
H^1(X,\ring X.(m(K_X+\Delta)-S))=0,
$$
by Kawamata-Viehweg vanishing, since
\begin{align*}
m(K_X+\Delta)-S &= (m-1)(K_X+\Delta)+K_X+A+B\\
                &\sim_{\mathbb{Q}}K_X+G+B+(m-1)(K_X+\Delta+tH),
\end{align*}
and $K_X+\Delta+tH$ is big and nef.  \end{proof}

\begin{lemma}\label{l_convex} Let $X$ be a normal projective variety, let $S$ be a prime 
divisor and let $D_i$ be three, $i=1$, $2$, $3$, $\mathbb{R}$-Cartier divisors on $X$.
Suppose that $f_i\colon\rmap X.Z_i.$ are ample models of $D_i$, $i=1$, $2$, $3$, where
$f_i$ is birational and $f_i$ does not contract $S$.  Suppose that $Z_i$ is
$\mathbb{Q}$-factorial, $i=1$, $2$, and the induced birational map $\rmap Z_1.Z_2.$ is an 
isomorphism in a neighbourhood of the strict transforms of $S$.

If $D_3$ is a positive linear combination of $D_1$ and $D_2$ then the induced birational
map $\rmap Z_i.Z_3.$ is an isomorphism in a neighbourhood of the strict transforms of $S$.
\end{lemma}
\begin{proof} By assumption $D_3=\lambda_1D_1+\lambda_2D_2$, where $\lambda_i>0$.  If
$g_i\colon\map Z.Z_i.$ is the normalisation of the graph of $\rmap Z_1.Z_2.$ then $g_i$ is
by assumption an isomorphism in a neighbourhood of the strict transforms of $S$, for
$i=1$, $2$.  Let $f\colon\rmap X.Z.$ be the induced birational map.  Let $g\colon\map
W.X.$ resolve the indeterminacy of $f$ and $f_3$.  Replacing $X$ by $W$ and $D_i$ by
$g^*D_i$, we may assume that $f$ and $f_3$ are morphisms.

Let $H_i=f_{i*}D_i$, $i=1$, $2$, $3$.  As the ample model is a semiample model, by (4) of
\eqref{l_ample}, $E_i=D_i-f_i^*H_i\geq 0$ is $f_i$-exceptional.  As $H_i$ is ample, $i=1$
and $2$, $H=\lambda_1g_1^*H_1+\lambda _2g_2^*H_2$ is semiample.  Let $\Sigma\subset Z$ be
a curve.  As $\Sigma$ is not contracted by both $g_1$ and $g_2$, $g_i^*H_i\cdot \Sigma\geq
0$, with equality for at most one $i$.  Thus $H$ is ample.  Now
\begin{align*} 
f_3^*H_3+E_3&=D_3  \\ 
           &=\lambda_1D_1+\lambda_2D_2 \\ 
           &=\lambda_1f_1^*H_1+\lambda_2f_2^*H_2+\lambda_1 E_1+\lambda_2E_2 \\ 
           &=f^*(\lambda_1g_1^*H_1+\lambda_2g_2^*H_2)+\lambda_1 E_1+\lambda_2E_2 \\ 
           &=f^*H+\lambda_1 E_1+\lambda_2E_2.
\end{align*} 
Note that $E_3\leq \lambda_1E_1+\lambda_2E_2$, as $f^*H$ has no stable base locus.  Let
$T=f(S)$.  We may write
$$
f^*H=f_3^*H_3+E_3-(\lambda_1 E_1+\lambda_2E_2).
$$
Since $\lambda_1 E_1+\lambda_2E_2$ is $f$-exceptional in a neighbourhood of $f^{-1}(T)$,
$E_3=\lambda_1 E_1+\lambda_2E_2$ in a neighbourhood of $f^{-1}(T)$, by
\eqref{l_negativity}.  But then $f$ and $f_3$ contract precisely the same curves in the
same neighbourhood and so $Z$ is isomorphic to $Z_3$ in a neighbourhood of the strict
transforms of $S$.  \end{proof}

\begin{lemma}\label{l_special} Assume Theorem~\ref{t_special}$_{n}$ and
Theorem~\ref{t_model}$_n$.  

Let $(X,\Delta_0=S+A+B_0)$ be a log smooth projective pair of dimension $n$, where $A\geq
0$ is a general ample $\mathbb{Q}$-divisor, $\rdown\Delta_0.=S$ is a prime divisor and
$B_0\geq 0$.  Suppose that $K_X+\Delta_0$ is pseudo-effective and $S$ is not a component
of $N_{\sigma}(K_X+\Delta_0)$.  Let $V_0$ be a finite dimensional affine subspace of
$\WDiv_{\mathbb{R}}(X)$ containing $B_0$, which is defined over the rationals.

Then we may find a general ample $\mathbb{Q}$-divisor $H\geq 0$, a log terminal model
$\phi\colon\rmap X.Y.$ for $K_X+\Delta_0+H$ and a positive constant $\alpha$, such that if
$B\in V_0$ and
$$
\|B-B_0\|<\alpha t,
$$
for some $t\in (0,1]$, then there is a log terminal model $\psi\colon\rmap X.Z.$ of
$K_X+\Delta+tH=K_X+S+A+B+tH$, which does not contract $S$, such that the induced
birational map $\chi\colon\rmap Y.Z.$ is an isomorphism in a neighbourhood of the strict
transforms of $S$.
\end{lemma}
\begin{proof} Pick general ample $\mathbb{Q}$-Cartier divisors $\llist H.k.$, which span
$\WDiv_{\mathbb{R}}(X)$ modulo numerical equivalence, and let $V$ be the affine subspace
of $\WDiv_{\mathbb{R}}(X)$ spanned by $V_0$ and $\llist H.k.$.  Let $\mathcal{C}\subset
\mathcal{L}_{S+A}(V)$ be a convex subset spanning $V$, which does not contain $\Delta_0$,
but whose closure contains $\Delta_0$, such that if $\Delta\in \mathcal{C}$ then
$K_X+\Delta$ is purely log terminal, $\Delta-\Delta_0$ is ample and the support of
$\Delta-\Delta_0$ contains the support of the sum $\alist H.+.k.$.  In particular if
$\Delta\in \mathcal{C}$ then $N_{\sigma}(K_X+\Delta)\subset N_{\sigma}(K_X+\Delta_0)$.  As
the coefficients $\sigma_C$ of $N_{\sigma}$ are continuous on the big cone,
\eqref{d_sigma}, possibly replacing $\mathcal{C}$ by a subset, we may assume that if
$\Delta\in \mathcal{C}$ then $N_{\sigma}(K_X+\Delta)$ and $N_{\sigma}(K_X+\Delta_0)$ share
the same support.  Moreover, if $\Delta\in \mathcal{C}$ then $K_X+\Delta$ is big and so,
as we are assuming Theorem~\ref{t_model}$_n$, $K_X+\Delta$ has a log terminal model
$\phi\colon\rmap X.Y.$, whose exceptional divisors are given by the support of
$N_{\sigma}(K_X+\Delta_0)$.  In particular $\phi$ does not contract $S$.

Given $\phi\colon\rmap X.Y.$, define a subset $\mathcal{S}_{\phi}\subset\mathcal{C}$ as
follows: $\Delta\in\mathcal{S}_{\phi}$ if and only if there is a log terminal model
$\psi\colon\rmap X.Z.$ of $K_X+\Delta$ such that
\begin{itemize}
\item the induced rational map $\chi\colon\rmap Y.Z.$ is an isomorphism in a neighbourhood
of the strict transforms of $S$.  
\end{itemize}
Define a subset $\mathcal{S}'_{\phi}\subset \mathcal{S}_{\phi}$ by requiring in addition 
that 
\begin{itemize}
\item $\psi$ is isomorphic to the ample model of $K_X+\Delta$ in a neighbourhood 
of the strict transforms of $S$. 
\end{itemize}

As we are assuming Theorem~\ref{t_special}$_{n}$, there are finitely many $1\leq j\leq l$
birational maps $\phi_j\colon\rmap X.Y_i.$ such that 
$$
\mathcal{C}=\bigcup_{j=1}^l \mathcal{S}_{\phi_j}.
$$
On the other hand, \eqref{l_ample-trick} implies that 
$$
\bigcup_{j=1}^l \mathcal{S}'_{\phi_j},
$$
is a dense open subset of $\mathcal{C}$.  

We will now show that the sets $\mathcal{S}'_{\phi}$ are convex.  Suppose that
$\Delta_i\in \mathcal{S}'_{\phi}$, $i=1$, $2$ and $\Delta$ is a convex linear combination
of $\Delta_1$ and $\Delta_2$.  Then $\Delta\in \mathcal{C}$ and so $K_X+\Delta$ has a log
terminal model $\psi\colon\rmap X.Z.$.  By (3) of \eqref{l_ample} there is a birational
morphism $h\colon\map Z.Z'.$ to the ample model of $K_X+\Delta$.  Let $\psi_i\colon\rmap
X.Z_i.$ be the ample model of $K_X+\Delta_i$.  \eqref{l_convex} implies that the induced
birational map $\rmap Z_i.Z'.$ is an isomorphism in a neighbourhood of the strict
transforms of $S$.  As $Z$ and $Z_i$ are isomorphic in codimension $1$, $\map Z.Z'.$ is
small in a neighbourhood of the strict transforms of $S$.  As $Z$ and $Z'$ are
$\mathbb{Q}$-factorial, the morphism $\map Z.Z'.$ is also an isomorphism in a
neighbourhood of the strict transforms of $S$.  Thus $\Delta\in \mathcal{S}'_{\phi}$ and
so $\mathcal{S}'_{\phi}$ is convex.

Shrinking $\mathcal{C}$, we may therefore assume that
$\mathcal{C}=\mathcal{S}'_{\phi_{j}}$ for some $1\leq j\leq l$.
  
Pick $\Delta=S+A+B\in \mathcal{C}$ and pick $H \sim_{\mathbb{Q}} \Delta-\Delta_0$ a
general ample $\mathbb{Q}$-divisor.  Pick a positive constant $\alpha$ such that if
$\|B-B_0\|<\alpha t$ for some $B\in V_0$ and $t\in (0,1]$ then
$S+A+B+t(\Delta-\Delta_0)\in \mathcal{C}$.  By \eqref{l_equiv-models},
$K_X+S+A+B+t(\Delta-\Delta_0)$ and $K_X+S+A+B+tH$ have the same log terminal models.  It
follows that $\alpha$ and $H$ have the required properties.
\end{proof}

\begin{lemma}\label{l_C-to-D} Theorem~\ref{t_effective}$_{n-1}$,
Theorem~\ref{t_special}$_{n}$ and Theorem~\ref{t_model}$_n$ imply
Theorem~\ref{t_effective}$_n$.
\end{lemma}
\begin{proof} By \eqref{l_ridiculous}, it suffices to prove this result for the generic
fibre of $U$.  Thus we may assume that $U$ is a point, so that $X$ is a projective
variety.

Let $f\colon\map Y.X.$ be a log resolution of $(X,\Delta)$.  We may write
$$
K_Y+\Gamma=f^*(K_X+\Delta)+E,
$$
where $\Gamma\geq 0$ and $E\geq 0$ have no common components, $f_*\Gamma=\Delta$ and
$f_*E=0$.  If $F\geq 0$ is an $\mathbb R$-divisor whose support equals the union of all
$f$-exceptional divisors then $\Gamma+F$ is big.  Pick $F$ so that $(X,\Gamma+F)$ is kawamata log
terminal.  Replacing $(X,\Delta)$ by $(Y,\Gamma+F)$ we may therefore assume that
$(X,\Delta)$ is log smooth.  By \eqref{l_big-to-ample} we may assume that $\Delta=A+B$,
where $A$ is a general ample $\mathbb{Q}$-divisor and $B\geq 0$.  By \eqref{l_bounded} and
\eqref{l_unbounded}, we may therefore assume that $\Delta=S+A+B$, where $(X,\Delta)$ is
log smooth, $A$ is a general ample $\mathbb{Q}$-divisor, and $\rdown\Delta.=S$ is a prime
divisor, which is not a component of $N_{\sigma}(K_X+\Delta)$.

Let $V$ be the subspace of $\WDiv_{\mathbb{R}}(X)$ spanned by the components of $B$.  By
\eqref{l_special} we may find a constant $\alpha>0$, a general ample $\mathbb{Q}$-divisor
$H$ on $X$, and a log terminal model $\phi\colon\rmap X.Y.$ of $K_X+\Delta+H$ 
such that if $B'\in V$, $t\in (0,1]$ and $\|B-B'\|<\alpha t$ then there is a log terminal
model $\phi'\colon\rmap X.Y'.$ of $K_X+S+A+B'+tH$ such that the induced rational map
$\chi\colon\rmap Y.Y'.$ is an isomorphism in a neighbourhood of the strict transforms of
$S$.  Pick $\epsilon>0$ such that $A-\epsilon H$ is ample.

Let $T$ be the strict transform of $S$ on $Y$ and define $\Phi_0$ on $T$ by adjunction
$$
(K_Y+T)|_T=K_T+\Phi_0.  
$$
Let $W$ be the subspace of $\WDiv_{\mathbb{R}}(T)$ spanned by the components of
$(\phi_*B)|_T+\Phi_0$ and let $L\colon\map V.W.$ be the rational affine linear map
$L(B')=\Phi_0+\phi_*B'|_T$.  Let $\Gamma=\phi_*\Delta$ and let $C=\phi_*A|_T$.  If we
define $\Psi$ on $T$ by adjunction
$$
(K_Y+\Gamma)|_T=K_T+\Psi,
$$
then $\Psi=C+L(B)$.  As $K_T+C+L(B)+t\phi_*H|_T$ is nef for any $t>0$, it follows that $K_T+\Psi$
is nef.  Since $C$ is big and $K_T+\Psi$ is kawamata log terminal, \eqref{l_linear-ample}
implies that there is a rational affine linear isomorphism $L'\colon\map
\WDiv_{\mathbb{R}}(T).\WDiv_{\mathbb{R}}(T).$ which preserves $\mathbb{Q}$-linear
equivalence, an ample $\mathbb{Q}$-divisor $G$ on $T$ and a rational affine linear
subspace $W'$ of $\WDiv_{\mathbb{R}}(T)$ such that $L'(U)\subset \mathcal{L}_G(W')$, where
$U\subset \mathcal{L}_C(W)$ is a neighbourhood of $\Psi$.  $\mathcal{N}_G(W')$ is a rational
polytope, by \eqref{t_polytopal}.  In particular we may find a rational polytope
$\mathcal{C}\subset V$ containing $B$ such that if $B'\in\mathcal{C}$ then
$K_X+\Delta'=K_X+S+A+B'$ is purely log terminal and
$$
(K_Y+\Gamma')|_T=K_T+\Psi',
$$
is nef, where $\Gamma'=\phi_*\Delta'$ and $\Psi'=C+L(B')$.  Pick a positive integer $k$
such that if $r(K_Y+\Gamma')$ is integral then $rk(K_Y+\Gamma')$ is Cartier in a
neighbourhood of $T$.  By \eqref{c_base} there is then a constant $m>0$ such that $mr(K_Y+\Gamma')$ is Cartier in a
neighbourhood of $T$ and
$mr(K_T+\Psi')$ is base point free.

\eqref{l_diophantine} implies that there are real numbers $r_i>0$ with $\sum r_i=1$,
positive integers $p_i>0$ and $\mathbb{Q}$-divisors $B_i\in\mathcal{C}$ such that
$$
p_i(K_X+\Delta_i),
$$
is integral, where $\Delta_i=S+A+B_i$, 
$$
K_X+\Delta=\sum r_i(K_X+\Delta_i),
$$
and
$$
\|B_i-B\|< \frac{\alpha\epsilon}{m_i},
$$
where $m_i=m p_i$ .  Let $\Psi_i=C+L(B_i)$.  By our choice of $p_i$, $m_i(K_T+\Psi_i)$
is base point free and so
$$
h^0(T,\ring T.(m_i(K_T+\Psi_i)))>0.
$$
Let $t_i=\epsilon /m_i$.  By \eqref{l_special} there is a log terminal model
$\phi_i\colon\rmap X.Y_i.$ of $K_X+\Delta_i+t_iH$, such that the induced birational map
$\chi\colon\rmap Y.Y_i.$ is an isomorphism in a neighbourhood of $T$ and the strict
transform $T_i$ of $S$.  In particular if $\Gamma_i=T_i+\phi_{i*}A+\phi_{i*}B_i$ and
$\tau\colon\map T.T_i.$ is the induced isomorphism then
$$
\tau^*(K_{Y_i}+\Gamma_i)|_{T_i}=K_T+\Psi_i,
$$
and so the pair $(Y_i,\Gamma_i)$ clearly satisfies conditions (1), (2), (4) and (7) of
\eqref{l_pseudo}.  As the induced birational map $\rmap Y_i.Y.$ is an isomorphism in a
neighbourhood of $T_i$, (3) of \eqref{l_pseudo} also holds.  As
$$
(m_i-1)t_i=\left (\frac{m_i-1}{m_i}\right ) \epsilon <\epsilon,
$$
$A-(m_i-1)t_iH$ is ample and so we may pick a general ample $\mathbb{Q}$-divisor $L_i$
such that $L_i\sim_{\mathbb{Q}} A-(m_i-1)t_iH$.  Then $K_X+S+L_i+B_i$ is purely log
terminal and as $\phi_i$ is $(K_X+\Delta_i+t_iH)$-negative and $H$ is ample, $\phi_i$ is
$(K_X+S+L_i+B_i)$-negative.  It follows that $K_{Y_i}+T_i+\phi_{i*}(L_i+B_i)$ is purely
log terminal so that $K_{Y_i}+\phi_{i*}(L_i+B_i)$ is kawamata log terminal, and conditions
(5) and (6) of \eqref{l_pseudo} hold.  Therefore, \eqref{l_pseudo} implies that
$$
h^0(Y_i,\ring Y_i.(m_i(K_{Y_i}+\Gamma_i)))>0.
$$

As $\phi_i$ is $(K_X+\Delta_i+t_iH)$-negative and $H$ is ample, it follows that $\phi_i$
is $(K_X+\Delta_i)$-negative.  But then
$$
h^0(X,\ring X.(m_i(K_X+\Delta_i)))=h^0(Y,\ring Y.(m_i(K_Y+\Gamma_i)))>0.
$$
In particular there is an $\mathbb{R}$-divisor $D$ such that
\[
K_X+\Delta=\sum r_i (K_X+\Delta_i)\sim_{\mathbb{R}} D\geq 0. \qedhere
\]
\end{proof}
\section{Finiteness of models}
\label{s_finite}

\begin{lemma}\label{l_cover} Assume Theorem~\ref{t_model}$_n$ and Theorem~\ref{t_effective}$_n$.
 
Let $\pi\colon\map X.U.$ be a projective morphism of normal quasi-projective varieties,
where $X$ has dimension $n$.  Let $V$ be a finite dimensional affine subspace of $\WDiv
_{\mathbb{R}}(X)$, which is defined over the rationals.  Fix a general ample
$\mathbb{Q}$-divisor $A$ over $U$.  Let $\mathcal{C}\subset \mathcal{L}_A(V)$ be a
rational polytope such that if $\Delta\in \mathcal{C}$ then $K_X+\Delta$ is kawamata log
terminal.

Then there are finitely many rational maps $\phi_i\colon\rmap X.Y_i.$ over $U$, $1\leq
i\leq k$, with the property that if $\Delta\in\mathcal{C}\cap \mathcal{E}_{A,\pi}(V)$ then
there is an index $1\leq i\leq k$ such that $\phi_i$ is a log terminal model of
$K_X+\Delta$ over $U$.
\end{lemma}
\begin{proof} Possibly replacing $V_A$ by the span of $\mathcal{C}$, we may assume that
$\mathcal{C}$ spans $V_A$.  We proceed by induction on the dimension of $\mathcal{C}$.

Suppose that $\Delta_0\in \mathcal{C}\cap \mathcal{E}_{A,\pi}(V)$.  As we are assuming
Theorem~\ref{t_effective}$_n$ there is an $\mathbb{R}$-divisor $D_0\geq 0$ such that
$K_X+\Delta_0 \sim_{\mathbb{R},U}D_0$ and so, as we are assuming
Theorem~\ref{t_model}$_n$, there is a log terminal model $\phi\colon\rmap X.Y.$ over $U$
for $K_X+\Delta_0$.  In particular we may assume that $\dim \mathcal{C}>0$.

First suppose that there is a divisor $\Delta_0\in\mathcal{C}$ such
that $K_X+\Delta_0\sim_{\mathbb{R},U}0$.  Pick $\Theta\in\mathcal{C}$, $\Theta\neq
\Delta_0$.  Then there is a divisor $\Delta$ on the boundary of $\mathcal{C}$ such that
$$
\Theta-\Delta_0=\lambda(\Delta-\Delta_0),
$$
for some $0<\lambda\leq 1$.  Now
\begin{align*}
K_X+\Theta=&\lambda(K_X+\Delta)+(1-\lambda)(K_X+\Delta_0)\\
         \sim_{\mathbb{R},U}&\lambda(K_X+\Delta).
\end{align*}
In particular $\Delta\in\mathcal{E}_{A,\pi}(V)$ if and only if
$\Theta\in\mathcal{E}_{A,\pi}(V)$ and \eqref{l_equiv-models} implies that $K_X+\Delta$ and
$K_X+\Theta$ have the same log terminal models over $U$.  On the other hand the boundary
of $\mathcal{C}$ is contained in finitely many affine hyperplanes defined over the
rationals, and we are done by induction on the dimension of $\mathcal{C}$.

We now prove the general case.  By \eqref{l_linear-dlt}, we may assume that $\mathcal{C}$
is contained in the interior of $\mathcal{L}_A(V)$.  Since $\mathcal{L}_A(V)$ is compact
and $\mathcal{C}\cap \mathcal{E}_{A,\pi}(V)$ is closed, it suffices to prove this result
locally about any divisor $\Delta_0\in \mathcal{C}\cap \mathcal{E}_{A,\pi}(V)$.  Let
$\phi\colon\rmap X.Y.$ be a log terminal model over $U$ for $K_X+\Delta_0$.  Let
$\Gamma_0=\phi_*\Delta_0$.  

Pick a neighbourhood $\mathcal{C}_0\subset \mathcal{L}_A(V)$ of $\Delta_0$, which is a
rational polytope.  As $\phi$ is $(K_X+\Delta_0)$-negative we may pick $\mathcal{C}_0$
such that for any $\Delta\in \mathcal{C}_0$, $a(F,K_X+\Delta)<a(F,K_Y+\Gamma)$ for all
$\phi$-exceptional divisors $F\subset X$, where $\Gamma=\phi_*\Delta$.  Since
$K_Y+\Gamma_0$ is kawamata log terminal and $Y$ is $\mathbb{Q}$-factorial, possibly
shrinking $\mathcal{C}_0$, we may assume that $K_Y+\Gamma$ is kawamata log terminal for
all $\Delta\in\mathcal{C}_0$.  In particular, replacing $\mathcal{C}$ by $\mathcal{C}_0$,
we may assume that the rational polytope $\mathcal{C}'=\phi_*(\mathcal{C})$ is contained
in $\mathcal{L}_{\phi_*A}(W)$, where $W=\phi_*(V)$.  By \eqref{l_linear-ample}, there is a
rational affine linear isomorphism $L\colon\map W.V'.$ and a general ample
$\mathbb{Q}$-divisor $A'$ over $U$ such that $L(\mathcal{C}')\subset\mathcal{L}_{A'}(V')$,
$L(\Gamma)\sim _{\mathbb{Q} , U} \Gamma$ for all $\Gamma \in \mathcal{C}'$ and
$K_Y+\Gamma$ is kawamata log terminal for any $\Gamma\in L(\mathcal{C}' )$.

Note that $\dim V'\leq \dim V$.  By \eqref{l_equiv-models} and \eqref{l_modelXY}, any log
terminal model of $(Y,L(\Gamma))$ over $U$ is a log terminal model of $(X,\Delta)$ over
$U$ for any $\Delta \in \mathcal{C}$. Replacing $X$ by $Y$ and $\mathcal{C}$ by
$L(\mathcal{C}')$, we may therefore assume that $K_X+\Delta_0$ is $\pi$-nef.

By \eqref{c_base} $K_X+\Delta_0$ has an ample model $\psi\colon\map X.Z.$ over $U$.  In
particular $K_X+\Delta_0\sim_{\mathbb{R},Z} 0$.  By what we have already proved there are
finitely many birational maps $\phi_i\colon\rmap X.Y_i.$ over $Z$, $1\leq i\leq k$, such
that for any $\Delta\in\mathcal{C}\cap \mathcal{E}_{A,\psi}(V)$, there is an index $i$
such that $\phi_i$ is a log terminal model of $K_X+\Delta$ over $Z$.  Since there are only
finitely many indices $1\leq i\leq k$, possibly shrinking $\mathcal{C}$,
\eqref{c_relative} implies that if $\Delta\in \mathcal{C}$ then $\phi_i$ is a log terminal
model for $K_X+\Delta$ over $Z$ if and only if it is log terminal model for $K_X+\Delta$
over $U$.

Suppose that $\Delta\in \mathcal{C}\cap \mathcal{E}_{A,\pi}(V)$.  Then $\Delta\in
\mathcal{C}\cap \mathcal{E}_{A,\psi}(V)$ and so there is an index $1\leq i\leq k$ such
that $\phi_i$ is a log terminal model for $K_X+\Delta$ over $Z$.  But then $\phi_i$ is a
log terminal model for $K_X+\Delta$ over $U$.  \end{proof}

\begin{lemma}\label{l_polyope-E} Assume Theorem~\ref{t_model}$_n$ and
Theorem~\ref{t_effective}$_n$.

Let $\pi\colon\map X.U.$ be a projective morphism of normal quasi-projective varieties,
where $X$ has dimension $n$. Suppose that there is a kawamata log terminal pair $(X,\Delta _{0})$.
Fix $A\geq 0$, a general ample $\mathbb{Q}$-divisor over
$U$.  Let $V$ be a finite dimensional affine subspace of $\WDiv_{\mathbb{R}}(X)$ which is
defined over the rationals.  Let $\mathcal{C}\subset \mathcal{L}_A(V)$ be a rational
polytope.

Then there are finitely many birational maps $\psi_j\colon\rmap X.Z_j.$ over $U$, $1\leq
j\leq l$ such that if $\psi\colon\rmap X.Z.$ is a weak log canonical model of $K_X+\Delta$
over $U$, for some $\Delta\in\mathcal{C}$ then there is an index $1\leq j\leq l$ and an
isomorphism $\xi\colon\map Z_j.Z.$ such that $\psi=\xi\circ\psi_j$.
\end{lemma}
\begin{proof} Suppose that $\Delta\in \mathcal{C}$ and $K_X+\Delta' \sim_{\mathbb{R},U}
K_X+\Delta$ is kawamata log terminal.  \eqref{l_equiv-models} implies that
$\psi\colon\rmap X.Z.$ is a weak log canonical model of $K_X+\Delta$ over $U$ if and only
if $\psi$ is a weak log canonical model of $K_X+\Delta'$ over $U$.  By
\eqref{l_linear-dlt} we may therefore assume that if $\Delta\in\mathcal{C}$ then
$K_X+\Delta$ is kawamata log terminal.

Let $G$ be any divisor which contains the support of every element of $V$ and let
$f\colon\map Y.X.$ be a log resolution of $(X,G)$.  Given $\Delta\in \mathcal{L}_A(V)$ we
may write 
$$
K_Y+\Gamma=f^*(K_X+\Delta)+E,
$$
where $\Gamma$ and $E\geq 0$ have no common components, $f_*\Gamma=\Delta$ and $f_*E=0$.
If $\psi\colon\rmap X.Z.$ is a weak log canonical model of $K_X+\Delta$ over $U$ then
$\psi\circ f\colon\rmap Y.Z.$ is a weak log canonical model of $K_Y+\Gamma$ over $U$.  If
$\mathcal{C}'$ denotes the image of $\mathcal{C}$ under the map $\map \Delta.\Gamma.$ then
$\mathcal{C}'$ is a rational polytope, cf. \eqref{l_polytopal}, and if $\Gamma\in
\mathcal{C}'$ then $K_Y+\Gamma$ is kawamata log terminal.  In particular
$\mathbf{B}_+(f^*A/U)$ does not contain any non kawamata log terminal centres of $K_Y+\Gamma$, for any
$\Gamma\in \mathcal{C}'$.  Let $W$ be the subspace of $\WDiv_{\mathbb{R}}(Y)$ spanned by
the components of the strict transform of $G$ and the exceptional locus of $f$.  By
\eqref{l_linear-ample} and \eqref{l_equiv-models}, we may assume that there is a general
ample $\mathbb{Q}$-divisor $A'$ on $Y$ over $U$ such that
$\mathcal{C}'\subset\mathcal{L}_{A'}(W)$.  Replacing $X$ by $Y$ and $\mathcal{C}$ by
$\mathcal{C}'$, we may therefore assume that $X$ is smooth.

Pick general ample $\mathbb{Q}$-Cartier divisors $\llist H.p.$ over $U$, which generate
$\WDiv_{\mathbb{R}}(X)$ modulo relative numerical equivalence and let $H=\alist H.+.k.$ be
their sum.  By \eqref{l_linear-dlt} we may assume that if $\Delta\in \mathcal{C}$ then
$\Delta$ contains the support of $H$.  Let $W$ be the affine subspace of
$\WDiv_{\mathbb{R}}(X)$ spanned by $V$ and the divisors $\llist H.p.$.  Pick
$\mathcal{C}'$ a rational polytope $\mathcal{L}_A(W)$ containing $\mathcal{C}$ in its
interior such that if $\Delta\in \mathcal{C}'$ then $K_X+\Delta$ is kawamata log
terminal. 

By \eqref{l_cover} there are finitely many $1\leq i\leq k$ rational maps
$\phi_i\colon\rmap X.Y_i.$ over $U$, such that given any $\Delta'\in \mathcal{C}'\cap
\mathcal{E}_{A,\pi}(W)$, we may find an index $1\leq i\leq k$ such that $\phi_i$ is a log
terminal model of $K_X+\Delta'$ over $U$.  By \eqref{c_polytopal} for each index $1\leq
i\leq k$ there are finitely many contraction morphisms $f_{i,m}\colon\map Y_i.Z_{i,m}.$
over $U$ such that if $\Delta'\in \mathcal{W}_{\phi_i,A,\pi}(W)$ and there is a contraction
morphism $f\colon\map Y_i.Z.$ over $U$, with
$$
K_{Y_i}+\Gamma_i=K_{Y_i}+{\phi_i}_*\Delta'\sim_{\mathbb{R},U}f^*D,
$$
for some $\mathbb{R}$-divisor $D$ on $Z$, then there is an index $(i,m)$ and an
isomorphism $\xi\colon\map Z_{i,m}.Z.$ such that $f=\xi\circ f_{i,m}$.  Let
$\psi_j\colon\rmap X.Z_j.$, $1\leq j\leq l$ be the finitely many rational maps obtained by
composing every $\phi_i$ with every $f_{i,j}$.

Pick $\Delta\in \mathcal{C}$ and let $\psi\colon\rmap X.Z.$ be a weak log canonical model
of $K_X+\Delta$ over $U$.  Since the components of $\Delta$ span $\WDiv_{\mathbb{R}}(X)$
modulo numerical equivalence over $U$, by \eqref{l_ample-trick}, we may find
$\Delta'\in\mathcal{C}'\cap \mathcal{E}_{A,\pi}(W)$ such that $Z$ is the ample model for
$K_X+\Delta'$ over $U$.  Since $\Delta'\in\mathcal{C}'\cap \mathcal{E}_{A,\pi}(W)$ there
is an index $1\leq i\leq k$ such that $\phi_i$ is a log terminal model of $K_X+\Delta'$
over $U$.  By \eqref{l_weak} there is a morphism $f\colon\map Y_i.Z.$, where
$$
K_{Y_i}+\Gamma_i=f^*(K_Z+\psi_*\Delta').
$$  
It follows that there is an index $m$ and an isomorphism $\xi\colon \map Z_{i,m}.Z.$ such
that $f=\xi\circ f_{i,m}$.  But then
$$
\psi=f\circ\phi_i=\xi\circ f_{i,m}\circ \phi_i =\xi\circ\psi_j,
$$ 
for some index $1\leq j\leq l$. 
\end{proof}

\begin{lemma}\label{l_D-to-E} Theorem~\ref{t_model}$_n$ and Theorem~\ref{t_effective}$_n$ imply
Theorem~\ref{t_finite}$_n$.
\end{lemma}
\begin{proof} Since $\mathcal{L}_A(V)$ is itself a rational polytope by
\eqref{l_polytopal}, this is immediate from \eqref{l_polyope-E}.
\end{proof}

\section{Finite generation}

\begin{lemma}\label{l_E-to-F} Theorem~\ref{t_model}$_n$, Theorem~\ref{t_effective}$_n$ and 
Theorem~\ref{t_finite}$_n$ imply Theorem~\ref{t_ezd}$_n$.
\end{lemma}
\begin{proof} Theorem~\ref{t_model}$_n$ and Theorem~\ref{t_effective}$_n$ imply that there
is is a log terminal model $\mu\colon\rmap X.Y.$ of $K_X+\Delta$, and
$K_Y+\Gamma=\mu_*(K_X+\Delta)$ is semiample by (1) of \eqref{l_weak}.  (1) follows, as
$$
R(X,K_X+\Delta)\simeq R(Y,K_Y+\Gamma).
$$

As $K_Y+\Gamma$ is semiample the prime divisors contained in the stable base locus of
$K_X+\Delta$ are precisely the exceptional divisors of $\mu$.  But there is a constant
$\delta>0$ such that if $\Xi\in V$ and $\|\Xi-\Delta\|<\delta$ then the exceptional
divisors of $\mu$ are also $(K_X+\Xi)$-negative.  Hence (2).

Note that by \eqref{c_finite}, there is a constant $\eta>0$ such that if $\Xi\in W$ and
$\|\Xi-\Delta\|<\eta$ then $\mu$ is also a log terminal model of $K_X+\Xi$.
\eqref{c_base} implies that there is a constant $r>0$ such that if $m(K_Y+\Gamma)$ is
integral and nef, then $mr(K_Y+\Gamma)$ is base point free.  It follows that if
$k(K_X+\Xi)/r$ is Cartier then every component of $\fix (k(K_X+\Xi))$ is contracted by
$\mu$ and so every such component is in the stable base locus of $K_X+\Delta$.  This is
(3).  \end{proof}
\section{Proof of Theorems}
\label{s_main}

\begin{proof}[Proof of Theorems~\ref{t_existence}, \ref{t_special}, ~\ref{t_model},
~\ref{t_effective}, ~\ref{t_finite} and ~\ref{t_ezd}] Immediate from the main result of \cite{HM07},
\eqref{l_A-to-B}, \eqref{l_B-to-C}, \eqref{l_C-to-D}, \eqref{l_D-to-E} and \eqref{l_E-to-F}.  \end{proof}

\begin{proof}[Proof of \eqref{t_main}] Suppose $K_X+\Delta$ is $\pi$-big.  Then we may
write $K_X+\Delta \sim_{\mathbb{R},U} B\geq 0$.  If $\epsilon>0$ is sufficiently small then
$K_X+\Delta+\epsilon B$ is kawamata log terminal.  \eqref{l_equiv-models} implies that
$K_X+\Delta$ and $K_X+\Delta+\epsilon B$ have the same log terminal models over $U$.
Replacing $\Delta$ by $\Delta+\epsilon B$ we may therefore assume that $\Delta$ is big
over $U$.  (1) follows by Theorem \ref{t_model} and Theorem~\ref{t_effective}.

 (2) and (3) follow from (1) and \eqref{c_base}.  \end{proof}

\section{Proof of Corollaries}
\label{s_corollaries}

\begin{proof}[Proof of \eqref{c_general}] (1), (2) and (3) are immediate from \eqref{t_main}.
(4) is Theorem D of \cite{EGZ06}.  \end{proof}

\begin{proof}[Proof of \eqref{c_generated}] Immediate by Theorem~5.2 of \cite{FM02} 
and (3) of \eqref{t_main}.  \end{proof}

\begin{proof}[Proof of \eqref{c_flops}] Note that $Y_1$ and $Y_2$ are isomorphic in
codimension one.  Replacing $U$ by the common ample model of $(X,\Delta)$, we may assume
that $K_{Y_i}+\Gamma_i$ is numerically trivial over $U$.  Let $H_2$ be a divisor on $Y_2$,
which is ample over $U$.  Let $H_1$ be the strict transform on $Y_1$.  Possibly replacing
$H_2$ by a small multiple, we may assume that $K_{Y_i}+\Gamma_i+H_i$ is kawamata log
terminal.
 
Suppose that $K_{Y_1}+\Gamma_1+H_1$ is not nef over $U$.  Then there is a
$(K_{Y_1}+\Gamma_1+H_1)$-flip over $U$ which is automatically a $(K_{Y_1}+\Gamma_1)$-flop
over $U$.  By finiteness of log terminal models for $K_X+\Delta$ over $U$, this
$(K_{Y_1}+\Gamma_1+H_1)$-MMP terminates.  Thus we may assume that $K_{Y_1}+\Gamma_1+H_1$
is nef over $U$.  But then $K_{Y_2}+\Gamma_2+H_2$ is the corresponding ample model, and so
there is a small birational morphism $f\colon\map Y_1.Y_2.$.  As $Y_2$ is
$\mathbb{Q}$-factorial, $f$ is an isomorphism.  \end{proof}

\begin{proof}[Proof of \eqref{c_finite}] We first prove (1) and (2).  By
Theorem~\ref{t_finite} and \eqref{c_polytopal}, and since ample models are unique by (1)
of \eqref{l_ample}, it suffices to prove that if $\Delta\in\mathcal{E}_{A,\pi }(V)$ then
$K_X+\Delta$ has both a log terminal model over $U$ and an ample model over $U$.

By \eqref{l_big-to-ample} and \eqref{l_equiv-models} we may assume that $K_X+\Delta$ is
kawamata log terminal.  \eqref{t_main} implies the existence of a log terminal model over
$U$ and the existence of an ample model then follows from \eqref{l_weak}.

(3) follows as in the proof of \eqref{c_polytopal}.  \end{proof}

\begin{proof}[Proof of \eqref{c_kodaira}] Immediate consequence of \eqref{c_finite}.
\end{proof}

\begin{proof}[Proof of \eqref{c_cox}] Let $V_A$ be the affine subspace of
$\WDiv_{\mathbb{R}}(X)$ generated by $\llist \Delta.k.$.  \eqref{c_finite} implies that
there are finitely many $1\leq p\leq q$ rational maps $\phi_p\colon\rmap X.Y_p.$ over $U$
such that if $\Delta\in \mathcal{E}_{A,\pi}(V)$ then there is an index $1\leq p\leq q$
such that $\phi_p$ is a log terminal model of $K_X+\Delta$ over $U$.  Let
$\mathcal{C}\subset \mathcal{L}_A(V)$ be the polytope spanned by $\llist \Delta.k.$ and
let
$$
\mathcal{C}_p=\mathcal{W}_{\phi_p,A,\pi}(V)\cap \mathcal{C}.
$$
Then $\mathcal{C}_p$ is a rational polytope.  Replacing $\llist\Delta.k.$ by the vertices
of $\mathcal{C}_p$, we may assume that $\mathcal{C}=\mathcal{C}_p$, and we will drop the
index $p$.  Let $\pi'\colon\map Y.U.$ be the induced morphism.  Let
$\Gamma_i=\phi_*\Delta_i$.  If we pick a positive integer $m$ so that both
$G_i=m(K_Y+\Gamma_i)$ and $D_i=m(K_X+\Delta_i)$ are Cartier for every $1\leq i\leq k$, 
then
$$
\mathfrak{R}(\pi,D^{\bullet})\simeq \mathfrak{R}(\pi',G^{\bullet}).
$$ 
Replacing $X$ by $Y$, we may therefore assume that $K_X+\Delta_i$ is nef over $U$.  By
\eqref{c_base} $K_X+\Delta_i$ is semiample over $U$ and so the Cox ring is finitely
generated.

\noindent\textbf{Aliter:} By \eqref{l_big-to-ample}, we may assume that each $K_X+\Delta_i$ is
kawamata log terminal.  Pick a log resolution $f\colon\map Y.X.$ of $(X,\Delta)$, where
$\Delta$ is the support of the sum $\Delta_1+\Delta_2+\dots+\Delta_k$.  Then we may write
$$
K_Y+\Gamma_i=\pi^*(K_X+\Delta_i)+E_i,
$$
where $\Gamma_i$ and $E_i$ have no common components, $f_*\Gamma_i=\Delta_i$ and
$f_*E_i=0$.  We may assume that there is an exceptional $\mathbb Q$-divisor $F\geq 0$ such
that $f^*A-F$ is ample over $U$ and $(Y,\Gamma_i+F)$ is kawamata log terminal.  Pick $A'
\sim_{\mathbb{Q},U}f^*A-F$ a general ample $\mathbb{Q}$-divisor over $U$.  Then
$$
G_i=K_Y+\Gamma_i+F-f^*A+A' \sim_{\mathbb{Q},U} K_Y+\Gamma_i,
$$
is kawamata log terminal and $\mathfrak{R}(\pi,D^{\bullet})$ if finitely generated if and
only if $\mathfrak{R}(\pi\circ f,G^{\bullet})$ if finitely generated, since they have
isomorphic truncations.  Replacing $X$ by $Y$ we may therefore assume that $(X,\Delta)$ is
log smooth.

Pick a positive integer $m$ such that $D_i=m(K_X+\Delta_i)$ is integral for $1\leq i\leq
k$.  Let 
$$
E=\bigoplus_{i=1}^k \ring X.(m\Delta_i).
$$ 
and let $Y=\gproj X.E.$ with projection map $f\colon\map Y.X.$.  Pick $\sigma_i\in \ring
X.(m\Delta_i)$ with zero locus $m\Delta_i$ and let $\sigma=(\llist \sigma.k.)\in
H^0(X,E)$.  Let $S$ be the divisor corresponding to $\sigma$ in $Y$.  Let $\llist T.k.$ be
the divisors on $Y$, given by the summands of $E$, let $T$ be their sum, and let
$\Gamma=T+S/m$.  Note that $\ring Y.(m(K_Y+\Gamma))$ is the tautological bundle associated
to $E(mK_X)$; indeed if $k=1$ this is clear and if $k>1$ then it has degree one on the
fibres and $\ring Y.(m(K_Y+\Gamma))$ restricts to the tautological line bundle on each
$T_i$ by adjunction and induction on $k$.  Thus
$$
\mathfrak{R}(\pi ,D^{\bullet})\simeq \mathfrak{R}(\pi \circ f ,m(K_Y+\Gamma)).
$$
On the other hand $(Y,\Gamma)$ is log smooth outside $T$ and restricts to a divisorially
log terminal pair on each $T_i$ by adjunction and induction.  Therefore $(Y,\Gamma)$ is
divisorially log terminal by inversion of adjunction.  Now $f^*A\leq S/m\leq \Gamma$.
Since $T$ is ample over $X$, there is a positive rational number  $\epsilon>0$ such that
$f^*A+\epsilon T$ is ample over $U$.  Let $A' \sim_{\mathbb{Q},U} f^*A+\epsilon T$ be a
general ample $\mathbb{Q}$-divisor over $U$.  Then
$$
K_Y+\Gamma'=K_Y+\Gamma-\epsilon T-f^*A+A' \sim_{\mathbb{Q},U} K_Y+\Gamma,
$$
is kawamata log terminal.  Thus $\mathfrak{R}(\pi \circ f ,m(K_Y+\Gamma))$ is finitely generated over $U$ by 
\eqref{c_generated}.  \end{proof}

We will need a well-known result on the geometry of the moduli spaces of $n$-pointed
curves of genus $g$:

\begin{lemma}\label{l_moduli} Let $X=\mgn g.n.$ and let $D$ be the sum of the boundary 
divisors.  

 Then 
\begin{itemize} 
\item $X$ is $\mathbb{Q}$-factorial, 
\item $K_X$ is kawamata log terminal and 
\item $K_X+D$ is log canonical and ample.  
\end{itemize} 
\end{lemma}
\begin{proof} $X$ is $\mathbb{Q}$-factorial, $K_X$ is kawamata log terminal and $K_X+D$ is 
log canonical as the pair $(X,D)$ is locally a quotient of a normal crossings pair.

It is proved in \cite{Mumford77} that $K_X+D$ is ample when $n=0$.  To prove the general
case, consider the natural map
$$
\pi\colon\map{\mgn g.n+1.}.{\mgn g.n.}.,
$$
which drops the last point.  Let $Y=\mgn g.n+1.$ and $G$ be the sum of the boundary
divisors.  If $\overline{\mathcal{M}}_{g,n}$ is the moduli stack of stable curves of genus
$g$, then $\mgn g.n.$ is the coarse moduli space and so there is a representable morphism
$$
f\colon\map \overline{\mathcal{M}}_{g,n}.{\mgn g.n.}.,
$$
which only ramifies over the locus of stable curves with automorphisms.  If $g=1$ and
$n\leq 1$ then $K_Y+G$ is obviously $\pi$-ample.  Otherwise the locus of smooth curves
with extra automorphisms has codimension at least two and in this case
$K_{\overline{\mathcal{M}}_{g,n}}+\Delta=f^*(K_X+D)$, where $\Delta$ is the sum of the
boundary divisors.  On the other hand there is a fibre square
$$
\begin{diagram}
\overline{\mathcal{M}}_{g,n+1}   &   \rTo    &  \mgn g.n+1.  \\
\dTo^{\psi} &  & \dTo_{\pi} \\
\overline{\mathcal{M}}_{g,n}   &   \rTo^f    &  \mgn g.n.,
\end{diagram}
$$  
where $\psi$ is the universal morphism.  Since the stack is a fine moduli space, $\psi$ is
the universal curve.  If $\Gamma$ is the sum of the boundary divisors then
$K_{\overline{\mathcal{M}}_{g,n+1}}+\Gamma$ has positive degree on the fibres of $\psi$,
by adjunction and the definition of a stable pair.  In particular $K_Y+G$ is ample on the
fibres of $\pi$.

On the other hand, we may write
$$
K_Y+G=\pi^*(K_X+D)+\psi,
$$
for some $\mathbb{Q}$-divisor $\psi$.  It is proved in \cite{GKM00} that $\psi$ is nef. 
We may assume that $K_X+D$ is ample, by induction on $n$.  If $\epsilon>0$ is sufficiently 
small, it follows that $\epsilon (K_Y+G)+(1-\epsilon)\pi^*(K_X+D)$ is ample.  But then 
\begin{align*} 
K_Y+G &= \epsilon(K_Y+G)+(1-\epsilon)(K_Y+G)\\ 
      &= \epsilon(K_Y+G)+(1-\epsilon)\pi^*(K_X+D)+(1-\epsilon)\psi,
\end{align*} 
is also ample and so the result follows by induction on $n$. \end{proof}

\begin{proof}[Proof of \eqref{c_curves}] By \eqref{l_moduli} $K_X+D$ is ample and log canonical,
where $D$ is the sum of the boundary divisors.  In particular $K_X+\Delta$ is kawamata log
terminal, provided none of the $a_i$ are equal to one.

Pick a general ample $\mathbb{Q}$-divisor $A \sim_{\mathbb{Q}} \delta(K_X+D)$.  Note that
\begin{align*} 
(1+\delta)(K_X+\Delta)  &= K_X+\delta(K_X+D)+(1+\delta)\Delta-\delta D \\ 
                        &\sim_{\mathbb{Q}} K_X+A+B. 
\end{align*} 
Now 
$$
0\leq (\Delta-\delta D)+\delta\Delta=B=\Delta+\delta(\Delta-D)\leq D.
$$
Thus the result is an immediate consequence of Theorem~\ref{t_model},
Theorem~\ref{t_effective} and Theorem~\ref{t_finite}. \end{proof}

\begin{proof}[Proof of \eqref{c_dream}] Immediate by \eqref{c_cox} and the main result of
\cite{HK00}.  \end{proof}

\begin{proof}[Proof of \eqref{c_fibre}] Pick any $\pi$-ample divisor $A$ such that
$K_X+\Delta+A$ is $\pi$-ample and kawamata log terminal.  We may find $\epsilon>0$ such
that $K_X+\Delta+\epsilon A$ is not $\pi$-pseudo-effective.  We run the
$(K_X+\Delta+\epsilon A)$-MMP over $U$ with scaling of $A$.  Since every step of this MMP
is $A$-positive, it is automatically a $(K_X+\Delta)$-MMP as well.  But as
$K_X+\Delta+\epsilon A$ is not $\pi$-pseudo-effective, this MMP must terminate with a Mori
fibre space $g\colon\map Y.W.$ over $U$.  \end{proof}

\begin{proof}[Proof of \eqref{c_batyrev}] As $X$ is a Mori dream space, the cone of
pseudo-effective divisors is a rational polyhedron.  It follows that $\cnone X.$ is also a
rational polyhedron, as it is the dual of the cone of pseudo-effective divisors.

Now suppose that $F$ is a co-extremal ray.  As $\cnone X.$ is polyhedral, we may pick a
pseudo-effective divisor $D$ which supports $F$.  Pick $\epsilon>0$ such that $D-\epsilon
(K_X+\Delta)$ is ample.  Pick a general ample $\mathbb{Q}$-divisor $A \sim_{\mathbb{Q}}
\frac 1 \epsilon D-(K_X+\Delta)$ such that $K_X+\Delta+kA$ is ample and kawamata log
terminal for some $k>1$.  Then $(X,\Theta=\Delta+A)$ is kawamata log terminal and
$K_X+\Theta\sim_{\mathbb{Q}}\frac 1 \epsilon D$ supports $F$.  As in the proof of
\eqref{c_fibre}, the $(K_X+\Delta)$-MMP with scaling of $A$ ends with a
$(K_X+\Theta)$-trivial Mori fibre space $f\colon\map Y.Z.$, and it is easy to see that the
pullback to $X$ of a general curve in the fibre of $f$ generates $F$.  \end{proof}

\begin{proof}[Proof of \eqref{c_existence}] The flip of $\pi$ is precisely the log
canonical model, so that this result follows from Theorem~\ref{t_model}.  \end{proof}

\begin{proof}[Proof of \eqref{c_scaling}] Immediate from \eqref{l_valuation} and
Theorem~\ref{t_finite}.
\end{proof}

\begin{proof}[Proof of \eqref{c_extract}] Pick an ample $\mathbb{Q}$-divisor $A\geq 0$
which contains the centre of every element of $\mathfrak{E}$ of log discrepancy one, but
no non kawamata log terminal centres.  If $\epsilon>0$ is sufficiently small, then
$(X,\Delta+\epsilon A)$ is log canonical and so replacing $\Delta$ by $\Delta+\epsilon A$,
we may assume that $\mathfrak{E}$ contains no valuations of log discrepancy one.
Replacing $\Delta$ by $(1-\eta)\Delta+\eta\Delta_0$, where $\eta>0$ is sufficiently small,
we may assume that $K_X+\Delta$ is kawamata log terminal.

We may write
$$
K_W+\Psi=f^*(K_X+\Delta)+E,
$$
where $\Psi\geq 0$ and $E\geq 0$ have no common components, $f_*\Psi=\Delta$ and $E$ is
exceptional.  Let $F$ be the sum of all the exceptional divisors which are neither
components of $E$ nor correspond to elements of $\mathfrak{E}$.

Pick $\epsilon>0$ such that $K_W+\Phi=K_Y+\Psi+\epsilon F$ is kawamata log terminal.  As
$f$ is birational, $\Phi$ is big over $X$ and so by \eqref{t_main} we may find a log
terminal model $g\colon\rmap W.Y.$ for $K_W+\Phi$ over $X$.  Let $\pi\colon \map Y.X.$ be
the induced morphism.  If $\Gamma=g_*\Phi$ and $E'=g_*(E+\epsilon F)$ then
$$
K_Y+\Gamma=\pi^*(K_X+\Delta)+E',
$$
where $K_Y+\Gamma$ is nef over $X$.  Negativity of contraction implies that $E'=0$ so that
$$
K_Y+\Gamma=\pi^*(K_X+\Delta).
$$
But then we must have contracted every exceptional divisor which does correspond to an
element of $\mathfrak{E}$.  \end{proof}

\begin{proof}[Proof of \eqref{c_model}] We may write
$$
K_Z+\Psi=f^*(K_X+\Delta)+E_0,
$$
where $\Psi\geq 0$ and $E_0\geq 0$ have no common components, $E_0$ is $f$-exceptional and
$f_*\Psi=\Delta$.  We may write $\Psi=\Psi_1+\Psi_2$ where every component of $\Psi_1$ has
coefficient less than one and every component of $\Psi_2$ has coefficient at least one.
Let $E_1$ be the sum of the components of $\Psi_1$ which are exceptional.  Pick $\delta>0$
such that $K_Z+\Psi_3$ is kawamata log terminal, where $\Psi_3=\Psi_1+\delta E_1$.  Let
$g\colon\rmap Z.Y.$ be a log terminal model of $K_Z+\Psi_3$ over $X$, whose existence is
guaranteed by \eqref{t_main}.  If $\pi\colon\map Y.X.$ is the induced birational morphism
then we may write
$$
K_Y+\Gamma=\pi^*(K_X+\Delta)+E,
$$
where $\Gamma=g_*(\Psi+\delta E_1)\geq 0$, $E=g_*(E_0+\delta E_1)\geq 0$ and
$K_Y+\Gamma_3$ is nef over $X$, where $\Gamma_3=g_*\Psi_3$.  If $E\neq 0$ then by
\eqref{l_negativity} there is a family of curves $\Sigma\subset Y$ contracted by $\pi$,
which sweeps out a component of $E$ such that $E\cdot \Sigma<0$.  Since $E$ and
$\Gamma-\Gamma_3$ have no common components, $(\Gamma-\Gamma_3)\cdot \Sigma\geq 0$, so
that $(K_Y+\Gamma_3)\cdot \Sigma<0$, a contradiction.  Thus $E=0$, in which case
$\Gamma_3=\Gamma_1$ and no component of $\Gamma_1$ is exceptional.  \end{proof}

\begin{lemma}\label{l_special-model} Let $(X,\Delta)$ be a quasi-projective 
divisorially log terminal pair and let $S$ be a component of the support of $\Delta$.

Then there is a small projective birational morphism $\pi\colon\map Y.X.$, where $Y$ is
$\mathbb{Q}$-factorial and $-T$ is nef over $X$, where $T$ is the strict transform of $S$.

In particular if $\Sigma$ is a curve in $Y$ which is contracted by $\pi$ and $\Sigma$
intersects $T$ then $\Sigma$ is contained in $T$.
\end{lemma}
\begin{proof} By (2.43) of \cite{KM98} we may assume that $(X,\Delta)$ is kawamata log
terminal.  Let $f\colon\map Z.X.$ be a log terminal model of $(X,\Delta)$ over $X$, let
$\Psi$ be the strict transform of $\Delta$ and let $R$ be the strict transform of $S$.
Let $g\colon\rmap Z.Y.$ be a log terminal model of $K_Z+\Psi-\epsilon R$ over $X$, for any
$\epsilon>0$ sufficiently small.  As $\Psi-\epsilon R$ is big over $X$ and $K_Z+\Psi$ is
kawamata log terminal, \eqref{l_linear-ample} and \eqref{c_finite} imply that we may
assume that $Y$ is independent of $\epsilon$.  In particular $-T$ is $\pi$-nef, where
$\pi\colon \map Y.X.$ is the induced morphism.
\end{proof}

\begin{proof}[Proof of \eqref{c_inversion}] We work locally about $S$.  Let $a$ be the log
discrepancy of $K_S+\Theta$ and let $b$ be the minimum of the log discrepancy with respect
to $K_X+\Delta$ of any valuation whose centre on $X$ is of codimension at least two.  It
is straightforward to prove that $b\leq a$, cf. (17.2) of \cite{Kollaretal}.  If $B$ is a
prime divisor on $S$ which not a component of $\Theta$ then the log discrepancy of $B$
with respect to $(S,\Theta)$ is one, so that $a\leq 1$.  In particular we may assume that
$b<1$.

By the main theorem of \cite{Kawakita05}, $(X,\Delta)$ is log canonical near the image of
$S$ if and only if $(S,\Theta)$ is log canonical and so we may assume that $(X,\Delta)$ is
log canonical and hence $b\geq 0$.

Suppose that $(X,\Delta)$ is not purely log terminal.  By \eqref{c_model} we may find a
birational projective morphism $\pi\colon\map Y.X.$ which only extracts divisors of log
discrepancy zero and if we write $K_Y+\Gamma=\pi^*(K_X+\Delta)$, then $Y$ is
$\mathbb{Q}$-factorial and $K_Y$ is kawamata log terminal.  By connectedness, see (17.4)
of \cite{Kollaretal}, the non kawamata log terminal locus of $K_Y+\Gamma$ is connected
(since we work locally about $S$) and contains the strict transform of $S$.
Thus we are free to replace $X$ by $Y$ and so we may assume that $X$ is
$\mathbb{Q}$-factorial and $K_X$ is kawamata log terminal.

Suppose that $(X,\Delta)$ is purely log terminal.  By \eqref{l_special-model} there is a
birational projective morphism $\pi\colon\map Y.X.$ such that $Y$ is
$\mathbb{Q}$-factorial and the exceptional locus is contained in the strict transform of
$S$.  Replacing $X$ by $Y$, we may assume that $X$ is $\mathbb{Q}$-factorial.

We may therefore assume that $X$ is $\mathbb{Q}$-factorial and $K_X$ is kawamata log
terminal.  Let $\nu$ be any valuation of log discrepancy $b$.  Suppose that the centre of
$\nu$ is not a divisor.  By \eqref{c_extract} there is a birational projective morphism
$\pi\colon\map Y.X.$ which extracts a single exceptional divisor $E$ corresponding to
$\nu$.  Since $X$ is $\mathbb{Q}$-factorial, the exceptional locus of $\pi$ is equal to
the support of $E$ and so $E$ intersects the strict transform of $S$.  Let
$\Gamma=\Delta'+(1-b)E$, where $\Delta'$ is the strict transform of $\Delta$.  Then
$K_Y+\Gamma=\pi^*(K_X+\Delta)$.

Replacing $X$ by $Y$ we may therefore assume there is a prime divisor $D\neq S$ on $X$,
whose coefficient $\Delta$ is $1-b$.  By \eqref{d_adjunction} some component $B$ of
$\Theta$ has coefficient at least $1-b$, that is log discrepancy $b$ and so $a\leq
b$.  \end{proof}

\begin{proof}[Proof of \eqref{c_moishezon}] Since this result is local in the \'etale
topology, we may assume that $U$ is affine.  Let $f\colon\map Y.X.$ be a log resolution of
$(X,\Delta)$, so that the composition $\psi\colon\map Y.U.$ of $f$ and $\pi$ is
projective.  We may write
$$
K_Y+\Gamma'=f^*(K_X+\Delta)+E,
$$
where $\Gamma'\geq 0$ and $E\geq 0$ have no common components and $E$ is $f$-exceptional.
If $F\geq 0$ is a $\mathbb{Q}$-divisor whose support equals the exceptional locus of $f$,
then for any $0<\epsilon \ll 1$, $(Y,\Gamma =\Gamma '+\epsilon F)$ is kawamata log
terminal.

Pick a $\psi$-ample divisor $H$ such that $K_Y+\Gamma+H$ is $\psi$-ample.  We run the
$(K_Y+\Gamma)$-MMP over $U$ with scaling of $H$.  Since $X$ contains no rational curves
contracted by $\pi$, this MMP is automatically a MMP over $X$.  Since $\Gamma$ is big over
$X$ and termination is local in the \'etale topology, this MMP terminates by
\eqref{c_scaling}. 

Thus we may assume that $K_Y+\Gamma$ is $f$-nef, so that $E+\epsilon F$ is empty, and
hence $f$ is small.  As $X$ is analytically $\mathbb{Q}$-factorial it follows that $f$ is
an isomorphism.  But then $\pi$ is a log terminal model.  \end{proof}

\bibliographystyle{hamsplain}
\bibliography{/home/mckernan/Jewel/Tex/math}

\end{document}